Università degli Studi "Roma Tre"

Facoltà di Scienze Matematiche, Fisiche e Naturali

Ph-D Thesis in Mathematics
(XXIV Ciclo)

# Syzygies, Pluricanonical Maps and the Birational Geometry of Irregular Varieties

Candidate:

Sofia Tirabassi

Advisor:

Prof. Giuseppe Pareschi (Università di Roma "Tor Vergata")

Head of the Doctoral School:
Prof. Luigi Chierchia

ACCADEMIC YEAR 2011-2012
December 18, 2011

**Abstract**

In this thesis we looked into three different problems which share, as a common factor, the exstensive use of the Fourier–Mukai transform as research tool.

In the first Part we investigated the syzygies of Kummer varieties (i.e. quotients of abelian varieties by the $\mathbb{Z}/2\mathbb{Z}$-action induced by the group operation), extending to higher syzygies results on projective normality and degree of equations of Sasaki ([68]), Kempf ([44]) and Khaled ([45, 46]).

The second Part of this Thesis is dedicated to the study of pluricanonical linear systems on varieties of maximal Albanese dimension. More precisely, in Chapter 3 we prove that the 4-canonical map of a smooth variety of general type and maximal Albanese dimension is always birational into its image, the content of this section can also be found in [72]. Chapter 4 is based on a joint work with Z. Jiang and M. Lahoz ([39]) in which we prove that, in any Kodaira dimension, the 4-canonical map of a smooth variety of maximal Albanese dimension induces the Iitaka fibration, while, in the case of varieties of general type, the 3-canonical map is sufficient (and hence the 3-canonical map of these varieties is always birational). We remark that these last results are both sharp.

Finally, in the last part of this thesis we consider the problem of classification of varieties with small invariants. The final goal of our investigation is to provide a complete cohomological charaterization of products of theta divisors by proving that every smooth projective variety $X$, of maximal Albanese dimension, with Euler characteristic equal to 1, and whose Albanese image is not fibered by tori is birational to a product of theta divisors. Under these hypothesis we show that the Albanese map has degree one. Furthermore, we present a new characterization of $\Theta$-divisor in principally polarized abelian varieties.

**Aknowledgements**

This thesis would never had been born without the apport and help of many people.

First of all I would like to thank the Università degli Studi Roma TRE (and the Ph. D grant commission L. Chierchia, M. Pontecorvo and R. Spigler) for the opportunity given to me and for awarding me a Ph. D grant.

Secondly I am deeply indebted to my advisor, G. Pareschi, that in these three years taught me a lot of mathematics, proposed to me many engaging problems and guided me to their solution.

Afterward, I owe very much to my doctoral commission, C. Ciliberto, A. F. Lopez, M. Popa and R. Pardini. I am especially grateful R. Pardini, who gave me many suggestions and comments, and to M. Popa for the time and efforts they dedicated to the reading of this work.

A whole chapter of this thesis would not be there if not for the help of my two coauthors, Z. Jiang, and M. Lahoz, to whom I am really indebted.

During these three years in Rome I had the occasion to meet and talk to many mathematicians and every single one of them taught me many things. Hence here I thank for the many mathematical conversations E. Arbarello, L. Caporaso, C. Ciliberto, F. Flamini, A. F. Lopez, G. Pareschi, A. Rapagnetta, E. Sernesi, S. Verra, F. Viviani.

Finally I have to thank all my fellow Ph. D students (both former and new, Lorenzo and Fulvio obviously included!) for the friendly environment they provided at the mathematics department of Università degli Studi Roma TRE.

Last (but not least) it comes my family (my parents and my granny) and friends (especially the Argelato's RPG group!) who supported me during my Ph. D. In particular I am deeply grateful to my future husband, Lorenzo, and his family. Without you I would never have had the strenght to accomplish this, especially during these last few crazy months before deadline!

I hope I did not forget anyone, if that is the case I am deeply sorry for that. Thanks to you all!!!

# INTRODUCTION

Abelian varieties are among the most studied objects in algebraic geometry. In 1981 Mukai ([53]) developed a tool (called Fourier-Mukai transform) in order to study moduli spaces of deformation of sheaves on abelian varieties. In this thesis we use the Fourier-Mukai functors in order to investigate geometric objects that are closely related to abelian varieties: Kummer varieties (i.e. quotients of abelian varieties by the natural $\mathbb{Z}/2\mathbb{Z}$ action induced by the group operation) and maximal Albanese dimension varieties (varieties which admit a generically finite morphism into an abelian variety). In particular we analyzed three different problems:

- we investigated the syzygies of Kummer varieties unifying and enhancing results of Sasaki ([68]), Kempf ([44]) and Khaled ([45, 46]);
- we studied pluricanonical systems on varieties of maximal Albanese dimension, improving the work of Chen–Hacon ([9, 10]), Pareschi–Popa ([60]) and Jiang ([36]);
- we employed the Fourier-Mukai transform to the classification problem of irregular varieties with given invariants.

The problem of projective normality and degree of defining equations of a curve $C$ embedded in the projective space by a very ample linear system $|D|$ is very classical. In [24] Green realized that both questions are different faces of a wider problem about the syzygies of $C$, and that this issue could be addressed by computing the cohomology of the Koszul complex of certain vector bundles on $C$. More precisely the question is to calculate the minimal degree of generators of the section algebra $R_{\mathscr{O}_C(D)} := \oplus_n H^0(C, \mathscr{O}_C(nD))$ over $S_{\mathscr{O}_C(D)} := \mathrm{Sym}^\bullet(H^0(C, \mathscr{O}_C(D))$. The techniques presented in [24] work in any dimension and in the last three decades a lot of work has been done to extend Green's work (see, for example [25, 26, 28, 27, 15] and the very recent preprint [17]). For what it concerns abelian varieties, the results obtained in the 70s by Koizumi ([47]) and Mumford ([55]) were recently generalized to higher syzygies by Kempf ([41]), Pareschi [58], Pareschi–Popa ([62]) and Lazarsfeld–Pareschi–Popa [51].

A Kummer variety $\mathcal{K}_A$ is quotient of abelian variety $A$ by the action of the involution $(-1)_A : A \to A$ defined by $a \mapsto -a$. Since ample line bundles on Kummer varieties could be easily described in terms of (even powers of) line bundels on the associate abelian varieties, Kummer varieties present a deep affinity with abelian varieties. In the second Chapter of this thesis we use this special kinship in order to extend to higher syzygies the results on projective normality and degree of defining equations of Kummer varieties obtained in the 90s by Sasaki, Kempf and Khaled. For a more precise account of these achievements we invite the reader to look over the introduction to Chapter 2; in fact the statements are somewhat thechnical and we deemed it better to postpone their presentation after having explained all the terminology involved. The main ideas behind the proofs is to "pull back" the problem on the Kummer to a problem on the abelian variety $A$ and then use the machinery granted by the Fourier-Mukai transform in order to solve it.

While in Chapter 2 we brought in to play integral transforms in order to study a quotient of an abelian variety, in the second and third Part of this thesis, we considered the triple given by a smooth complex projective variety $X$, its Albanese variety $\mathrm{Alb}(X)$, i.e. the dual torus to $H^1(X, \mathscr{O}_X)/H^1(X, \mathbb{Z})$, and its Albanese morphism

$$\mathrm{alb}_X : X \to \mathrm{Alb}(X).$$

Thus if previously we used the Fourier-Mukai transform in order to investigate the properties of sheaves pulled back from a variety $K_A$, now we will concentrate on sheaves pushed forward to an abelian variety. In both second and third Part of this thesis we focus our attention on varieties of *maximal Albanese dimension* (i.e. varieties whose Albanese map is generically finite into its image), heeding, in particular, to those of general type, whose canonical line bundle is big by definition. Thus, when $m$ is an integer large enough and divisible (meaning that $H^0(X, \omega_X^{\otimes m}) \neq 0$), the rational map induced by the $m$-canonical linear system is birational. In this case it is usually said, by a slight abuse, that the $m$-canonical linear system itself is birational.

One of the main issues about varieties of general type is an effectiveness problem that arises every time we have to deal with quantities that are "big enough". In fact it is very natural to ask onself if a bound can be found for such numbers. Hacon and McKernan in [31], and indipendently Takayama [71], proved that in any dimension there exists a bound $m_0$ depending only on the dimension of $X$, for which the $m$-pluricanonical map is birational for every divisible $m \geq m_0$. For what it concerns curves this bound is 3, and this result is an easy consequence of Riemann–Roch Theorem. Bombieri [6] showed that $m_0 = 5$ is the optimal bound for surfaces. As the dimension grows, the situation gets more and more complicated and still there are many open problems. For example, we cite the work of Chen–Chen who proved that 73 is a bound for threefolds, but it is not know if this is optimal. Some results on 4-folds were discovered by Di Biagio in his Ph-D thesis [14, 13]

The study of $m$-canonical map of irregular varieties was started by Chen and Hacon who in [10] proved that, in case of varieties of maximal Albanese dimension, the dependecy on the dimension of Hacon–McKernan was linear. However, in a later paper, [9] they realized that it was non existent, that is some bound can be found (3 for varieties with positive Euler characteristic, 6 for every other varieties) that worked in any dimension. The same results were later found again by Pareschi and Popa ([60]) as a consequence of the Fourier-Mukai based thecniques ideated in [61]. Jiang in [36] proved that the 5-canoncial map for varieties of maximal Albanese dimension is always birational, lowering by one the bound of Chen–Hacon. In addition he demonstrated that, even when the variety $Z$ is of intermediate Kodaira dimension the pluricanonical linear system $|mK_Z|$ induces the Iitaka fibration for every $m \geq 5$, shifting the attention from varieties of general type to varieties of any Kodaira dimension.

In Chapter 3 we studied the tetracanonical map and showed that, in the case of varieties of general type, it is always birational. However this result is not sharp, and in Chapter 4 we presents an improvement of its obtained in collaboration with Z. Jiang and M. Lahoz ([39]). Namely we were able to prove that the tricanonical map of varieties of general type and maximal Albanese dimension is birational; furthermore, indipendently from the Kodaira dimension, the 4-canonical map of varieties of maximal Albanese dimension induces the Iitaka fibration. We remark that both achievements are sharp: infact the bicanonical map of desingularizations of irreducible principal polarizations is not birational and, in addition, we were able to produce an example of variety of maximal Albanese dimension and intermediated Kodaira dimension whose tricanonical map could not induce the Iitaka fibration.

The second issue related to irregular varieties is the classification problem. It is well known that the holomorphic Euler characteristic of varieties of maximal Albanese dimension is non-negative. Ein–Lazarsfeld ([16]) proved that if the Euler characteristic of $X$, $\chi(X)$, is zero, then the Albanese image of $X$ is fibered in tori. Augmenting the Euler characteristic by one, we find the first examples of varieties of maximal Albanese dimension whose Albanese image is not fibered in translates of abelian subvarieties of $\mathrm{Alb}(X)$: smooth models of (irreducible) theta divisor in principally polarized abelian varieties. An interesting line of research started by Ein–Lazarsfeld ([16]) and later pursuit by Hacon ([30]), Hacon–Pardini ([32]), Lazarsfeld–Popa ([52]), Barja–Lahoz–Naranjo–Pareschi ([2]) and Pareschi ([59]), consists in characterizing smooth models of theta divisors by their birational invariants. The next step along the road of better comprehension of varieties of maximal Albanese dimension and Euler characteristic one would be to provide a similar description for *products* of irreducible theta divisors. In this setting Pareschi conjectured the following

> *A complex smooth projective variety $X$ of maximal Albanese dimension such that $\chi(X, \omega_X) = 1$ and its Albanese image is not fibered in tori is birational to a product of theta divisors.*

The above statement holds for surfaces, thank to the work of Beauville ([3]), who proved that a surface $S$ with irregularity $q(S) = 1$ and $\chi(S) = 1$ is a product of genus 2 curves, Pirola ([66]) and Hacon–Pardini ([33]) who studied surfaces $S$ with $\chi(S) = 1$ and $q(S) = 3$. The conjecture is also known to be true for highly irregular varieties thank to [34] where the authors proved that the irregularity $q(X) := h^1(X, \mathscr{O}_X)$ of a variety $X$ as above satisfies the inequality

$$q(X) \leq 2 \dim X$$

and equality holds if and only if $X$ is birational to a product of curves of genus 2. If Pareschi's conjecture were proved to be true it will reduce the problem of (birational) classification of varieties of maximal Albanese dimension and $\chi = 1$ to the study of varieties whose Albanese image is not of general type. Moreover, it could lead to a better understandig of those varieties whose bicanonical map is not birational, completing the work of Barja–Lahoz–Naranjo–Pareschi ([2]) and Lahoz ([48]).

In the last Part of this thesis we study smooth complex projective varieties $X$ of maximal Albanese dimension with Albanese image not fibered in subtori of $\mathrm{Alb}(X)$ and whose Euler charecteristic is equal to one, proving some partial results that hopefully will lead to the resolution of Pareschi's conjecture. In particular, we find that the Albanese morphism of such varieties is always birational and we used this fact to give a new cohomological characterization of $\Theta$ divisors under the further hypothesis of the Albanese image of $X$ being normal (again, since this last result is a bit technical we invite the reader to the introduction to the Chapter 5 for the complete statement).

**Notation:**

Through out this thesis we work on an algebraically closed field $k$; restrictions to the characteristic or to the field itself will be announced when needed. Unless otherwise stated, with the word "variety" we will mean a projective variety over $k$ and a "sheaf" $\mathscr{F}$ on $X$ will always stand for a coherent sheaf. The abelian category of coherent sheaves on $X$ will be denoted with $\mathbf{Coh}(X)$, while $\mathbf{D}(X)$ will be the bounded derived category of complexes of coherent sheaves on $X$, i. e. $\mathbf{D}(X) := \mathbf{D}^b(\mathbf{Coh}(X))$.

Given $\mathscr{F}$ a sheaf on $X$, its cohomology groups will be denoted by $H^i(X, \mathscr{F})$, or simply $H^i(\mathscr{F})$ when there is no chance of mistaking the variety $X$. By $h^i(X, \mathscr{F})$ (or simply $h^i(\mathscr{F})$) we will mean the dimension of $H^i(X, \mathscr{F})$ as $k$-vector space.

Let $x \in X$, by $k(x)$ we denote the skyscraper sheaf at $x$. Given $\mathscr{F}$ a coherent sheaf on a variety $X$, and $V$ a subspace of $H^0(X, \mathscr{F})$, we will denote by $\mathrm{Bs}(V)$ the *base locus of* $V$ i.e. the locus of points $x \in X$ where the map

$$V \otimes k(x) \to \mathscr{F} \otimes k(x)$$

fails to be surjective. If $V \simeq H^0(X, \mathscr{F})$ we will call $\mathrm{Bs}(V)$ *the base locus of* $\mathscr{F}$ and we will denote it by $\mathrm{Bs}(\mathscr{F})$. The reader could find a more accurate list of the symbols used in this thesis at the end of this document, just before the index.

# CONTENTS

# CHAPTER 1

# GENERIC VANISHING: BACKGROUND MATERIAL

We expose here some basic results in use throughout the thesis. Moreover, at the beginning of each Chapter, the reader will find more preliminary material.

## 1.1 Fourier-Mukai Functors

One of the main technical tools applied in this thesis is the *Fourier-Mukai* functor introduced by Mukai in [53] in order to study moduli space of deformations of Picard sheaves. It is constructed as follows: given $A$ an abelian variety of dimension $q$ and $\widehat{A}$ its dual, one can consider $\mathscr{P}_A$ (or simply $\mathscr{P}$ when there is no chance of confusion) the Poincaré line bundle on the product $A \times \widehat{A}$ and build the exact functor:

$$\mathbf{R}S_A := \mathbf{R}q_*(p^*(\cdot) \otimes \mathscr{P}) : \mathbf{D}(A) \longrightarrow \mathbf{D}(\widehat{A}) \tag{1.1.1}$$

where $p$ and $q$ are respectively the left and right projection from $A \times \widehat{A}$. Mukai's inversion thoeorem [53, Theorem 2.2] tells us that this functor is an equivalence of triangulated categories with a quasi-inverse given by his twin functor,

$$\mathbf{R}\widehat{S}_A := \mathbf{R}p_*(q^*(\cdot) \otimes \mathscr{P}),$$

composed with the exact equivalence $(-1_A)^*[+q]$, where the map $-1_A : A \to A$ is the "multiplication by -1" in the abelian variety, and $[\cdot]$ stands, as usual, for the shift functor in a triangulated category. In fact the followings composition formulas hold:

$$\mathbf{R}S_A \circ \mathbf{R}\widehat{S}_A \simeq (-1_A)^* \circ [-q], \quad \mathbf{R}\widehat{S}_A \circ \mathbf{R}S_A \simeq (-1_{\widehat{A}})^* \circ [-q]; \tag{1.1.2}$$

A key point of the proof is the remark that $\mathbf{R}S_A(\mathscr{O}_A) \simeq k(\widehat{0})[-q]$.

Often in the sequel, when there is no chance of confusion about the variety $A$, we will write just $\mathbf{R}S$ or $\mathbf{R}\widehat{S}$ instead of $\mathbf{R}S_A$ or $\mathbf{R}\widehat{S}_A$.

Denote with $R^iS(\mathscr{F})$ (respectively $R^i\widehat{S}(\mathscr{F})$) the $i$-th cohomology group of the complex $\mathbf{R}S(\mathscr{F})$ (respectively $\mathbf{R}\widehat{S}(\mathscr{F})$). Then we can define the followings

**Definition 1.1.1** ([53, Definition 2.3]). We say that *Weak Index Theorem* (in brief W.I.T.) holds for an object $\mathscr{F}$ in $\mathbf{D}(A)$ if the $R^iS_A(\mathscr{F})$ vanish for all but one $i$; this $i$ is denoted by $i(\mathscr{F})$ and called *index* of $\mathscr{F}$, the coherent sheaf $R^{i(\mathscr{F})}S_A(\mathscr{F})$ is denoted by $\widehat{\mathscr{F}}$ and is called the *Fourier-Mukai transform of* $\mathscr{F}$. By W.I.T.$(j)$ we denote the class of objects in $\mathbf{D}(A)$ that satisfies W.I.T. with index $j$.

We say that *Index Theorem* (in brief I.T.) holds for a coherent sheaf $\mathscr{F}$ on $A$ if for any $\alpha \in \widehat{A}$ and all but one $i$ we have

$$h^i(A, \mathscr{F} \otimes \alpha) = 0.$$

By I.T.$(j)$ we will mean the set of coherent sheaves on $A$ that satisfies I.T. with index $j$.

It can be proved using base change that a sheaf $\mathscr{F}$ satisfies I.T. with index $i$ if and only if it satisfies W.I.T. with index $i$ and its transform $\widehat{\mathscr{F}}$ is a vector bundle.

A very nice example of the use of the Fourier-Mukai functor in order to study sheaves on abelian varieties is the following cohomological characterization of principal polarizations, due to Hacon, tha we will be needing afterwards.

**Proposition 1.1.2** ([30, Proposition 2.2]). *Let $A$ be an abelian variety and $\mathscr{F}$ a torsion free sheaf of rank 1 (i. e. $\mathscr{F} \simeq \mathscr{I} \otimes L$ with $L$ a line bundle and $\mathscr{I}$ and ideal sheaf). Suppose that $\mathscr{F}$ satisfies I.T. with index* 0 *and that its Euler characteristic is 1. Then $\mathscr{F}$ is a line bundle with $h^0(X, \mathscr{F}) = 1$ hence it is a principal polarization.*

### 1.1.1 Relations with Other Functors

In this paragraph we report some results of [53, Section 3] that show the behavior of $\mathbf{R}S$ and $\mathbf{R}\widehat{S}$ with respect to other classical functors of algebraic geometry. Before going further we introduce a piece of notation: from now on, given a topologically trivial line bundle $\alpha$ on $A$, the symbol $[\alpha]$ will stand for the point of $\widehat{A}$ parametrizing $\alpha$ (via $\mathscr{P}$.), i. e. $\alpha \simeq \mathscr{P}_{A\times[\alpha]}$. Let $p \in A$ be a point, by $t_p : A \to A$ we will denote the morphism "translation by $p$" defined by $x \mapsto x + p$. In a similar way we will define the morphism $t^*_{[\alpha]}$ with $[\alpha] \in \widehat{A}$. The topologically trivial line bundle $\mathscr{P}_{\{p\}\times\widehat{A}}$ will be denoted by $P_p$.

**Proposition 1.1.3** (Exchange of translation and tensor product). *Given $p \in A$, there is the following isomorphism of functors*

$$\mathbf{R}S \circ t_p^* \simeq (- \otimes P_{-p}) \circ \mathbf{R}S.$$

*Conversely if* $[\alpha] \in \widehat{A}$ *we get*

$$\mathbf{R}S \circ (- \otimes \alpha) \simeq t^*_{[\alpha]} \circ \mathbf{R}S$$

### 1.1.2 Generalized Fourier-Mukai Functors

After the first appearence of the Fourier-Mukai transform in 1981 a lot of work has been done to study the behavior of this functor and to better understand his usage. The first step in this investigation has been to extend this tool to a more general setting. In fact, given two varieties $X$ and $Y$, any object in $\mathscr{E} \in \mathbf{D}(X \times Y)$ can substitute the role of the Poincaré bundle $\mathscr{P}$ in Definition (1.1.1). What we get is an exact functor $\mathbf{R}\Phi_{\mathscr{E}} : \mathbf{D}(X) \to \mathbf{D}(Y)$ that is called *Integral transform with kernel* $\mathscr{E}$. Usually, in literature, this functors are said to be Fourier-Mukai functors whenever they yield an equivalence of categories.

In a sequence of articles (see for example [60, 63, 64, 65]) Pareschi and Popa studied a very special integral transform: given a variety $X$ of dimension $n$ with a non trivial morphism to an abelian variety

$$a : X \longrightarrow A,$$

$\dim A = q$ one can consider the product $X \times \widehat{A}$. The line bundle $\mathscr{P}_a := (a \times \mathrm{id})^*\mathscr{P}$ on $X \times \widehat{A}$ is given and we may constuct the functor $\mathbf{R}\Phi_{\mathscr{P}_a}$. It is an easy consequence of projection formula and base change that

$$\mathbf{R}\Phi_{\mathscr{P}_a} \simeq \mathbf{R}S_A \circ \mathbf{R}a_*. \tag{1.1.3}$$

In what follows it will be necessary to consider the integral transform with kernel $\mathscr{P}_a^{-1}$. Since by the "See-saw" Principle (see [55][Cor. 6, pg. 54]) it is not difficult to show that

$$\mathscr{P}^{-1} \simeq (1_A \times -1_{\widehat{A}})^*\mathscr{P}$$

we get

$$\mathbf{R}\Phi_{\mathscr{P}_a^{-1}} \simeq (-1_{\widehat{A}})^*\mathbf{R}\Phi_{\mathscr{P}_a}. \tag{1.1.4}$$

Now for any smooth variety $Z$ we may consider the *dualizing functor*

$$\mathbf{R}\Delta_Z := \mathscr{H}om_{\mathscr{O}_Z}(-, \omega_Z).$$

Again, when there will not be any chance of confusion, we will omit the subscript $Z$.

A key result of Grothendieck and Verdier explains the behaviour of $\mathbf{R}\Delta_Z$ with respect with the derived direct image functors. In the sequel, however, when we will refer to *Grothendieck duality* we will mean the following statement that explains the mutual relation between integral transforms and duality functors.

**Theorem 1.1.4** (Grothendieck duality, [65, Lemma 2.2]). *Let $X$ be a smooth variety of dimension $n$ with $a : X \to A$ a non trivial morphism to an abelian variety· There is an isomorphism of functors from $\mathbf{D}(X) \to \mathbf{D}(\widehat{A})$*

$$\mathbf{R}\Delta_A \mathbf{R}\Phi_{\mathscr{P}_a} \simeq (-1_{\widehat{A}})^* \mathbf{R}\Phi_{\mathscr{P}_a} \mathbf{R}\Delta_X[n].$$

Applying both sides of the equality above to a given object $\mathscr{F}$ in $\mathbf{D}(X)$ and taking cohomology sheaves of both complexes, we got the following isomorphisms of sheaves

$$\mathscr{E}xt^i(\mathbf{R}\Phi_{\mathscr{P}_a}\mathscr{F}, \mathscr{O}_{\widehat{A}}) \simeq R^{n+i}\Phi_{\mathscr{P}_a^{-1}}\mathbf{R}\Delta\mathscr{F} \tag{1.1.5}$$

## 1.2 Cohomological Support Loci and GV-sheaves

In the previous Section we introduced integral transforms and explained their realations with other functors; now we are able to present another fundamental tool we will be using: sheaves that satisfy generic vanishing.

**Definition 1.2.1** (Cohomological support loci). Given a sheaf $\mathscr{F}$ on $X$ its *$i$-th cohomological support locus with respect to $a$* is :

$$V^i_a(X, \mathscr{F}) := \{[\alpha] \in \widehat{A} | h^i(X, \mathscr{F} \otimes a^*\alpha) > 0\}.$$

As it happens for cohomology groups, when possible we will omit the variety $X$ in the notation above.

**Example 1.** If $\mathscr{F}$ is a sheaf on an abelian variety satisfying I.T. with index 0, then its cohomological support loci $V^i_{\mathrm{id}}(\mathscr{F})$ are empty for every $i \geq 1$.

An important invariant associated to this cohomological loci is the following:

**Definition 1.2.2** ([64, Definition 3.1]). Given a coherent sheaf $\mathscr{F}$ on $X$, the *generic vanishing index* of $\mathscr{F}$ with respect to $a$ is

$$\mathrm{gv}_a(\mathscr{F}) := \min_{i>0}\{\mathrm{codim}_{\widehat{A}}(V^i(\mathscr{F})) - i\}.$$

If the $V^i$'s are empty for $i > 0$ then we say by definition that the generic vanishing index of $\mathscr{F}$ is $+\infty$. When $\mathscr{F} = \omega_X$ then $\mathrm{gv}_a(\mathscr{F})$ is called the *generic vanishing index of* $X$ (with respect to $a$) and it is denoted by $\mathrm{gv}_a(X)$.

The class of sheaves whose generic vanishing index is greater or equal $k$ is usually denoted by $GV_k$. A sheaf $\mathscr{F}$ whose generic vanishing index is non-negative is called *$GV$-sheaf* (*generic vanishing sheaf*). In literature, when $\mathscr{F}$ is a sheaf on an abelian variety $A$ with $\mathrm{gv}_{\mathrm{id}}(\mathscr{F}) \geq 1$, it is often said that $\mathscr{F}$ is *Mukai regular* or, in biref, *$M$-regular.* A well known result of Green–Lazarsfeld provides us with many examples of generic vanishing sheaves:

**Theorem 1.2.3** ([28]). *Given $X$ a smooth variety with a non trivial morphism $a : X \to A$ to an abelian variety, then*

$$\mathrm{gv}_a(X) \geq \dim a(X) - \dim X.$$

*In particular, if $X$ has maximal Albanese dimension, then $\omega_X$ is a GV-sheaf.*

The integral transform defined in paragraph 1.1.2 is a perfect instrument to study the geometry of the generic vanishing loci. An example of such thing is the following result that relates the "size" of the $V^i_a(\mathscr{F})$ with the vanishing of the cohomology sheaves of the transform $\mathbf{R}\Delta_X\mathscr{F}$.

**Theorem 1.2.4** (W.I.T. criterion, [65, Theorem A]). *Let $\mathscr{F}$ be a sheaf on $X$ and suppose $n = \dim X$. Then the following are equivalent:*

*(i)* $\mathrm{gv}_a(\mathscr{F}) \geq -k$ *for* $k \geq 0$*;*

*(ii)* $R^i\Phi_{\mathscr{P}_a}(\mathbf{R}\Delta\mathscr{F}) = 0$ *for every* $i \neq n-k, \dots, n$.

Observe that, by the above result, the generic vanishing index of a sheaf $\mathscr{F}$ is non-negative if and only if the complex $\mathbf{R}a_*(\mathbf{R}\Delta\mathscr{F})$ satisfies W.I.T. with index $n$. Thus Theorem 1.2.4 provides us with a criterion (that we will call *W. I. T. criterion*) that will help us to understand whether a sheaf satisfies the generic vanishing.

### 1.2.1 GV-sheaves

Now we focus on some features of generic vanishing sheaves. The second equivalent condition of Theorem 1.2.4 tells us that the full transform of $\mathbf{R}\Delta_X\mathscr{F}$ is indeed a sheaf concentrated in degree $n = \dim(X)$. In particular we have that the object $\mathbf{R}a_*\mathbf{R}\Delta\mathscr{F}$ satisfies W.I.T. with index $n$. Another peculiar property of $GV$-sheaves is stated in the following Lemma.

**Lemma 1.2.5** ([59, Corollary 3.2]). *Given $\mathscr{F}$ a GV-sheaf on $X$ with respect to some map $a$. Then*

$$V^d_a(\mathscr{F}) \subseteq \cdots \subseteq V^1_a(\mathscr{F}) \subseteq V^0_a(\mathscr{F}).$$

The following Proposition provides two basic properties for its Fourier transform $\widehat{\mathbf{R}a_*\mathbf{R}\Delta\mathscr{F}}$.

**Proposition 1.2.6** ([59, Proposition 1.6]). *Let $\mathscr{F}$ be a GV-sheaf on $X$ with respect to $a$. Then*

*(i)* $\mathrm{rk}\widehat{\mathbf{R}a_*\mathbf{R}\Delta\mathscr{F}} = \chi(\mathscr{F})$

*(ii)* $\mathbf{R}\Delta_{\widehat{A}}(\widehat{\mathbf{R}a_*\mathbf{R}\Delta_X\mathscr{F}}) \simeq (-1_{\widehat{A}})^*\mathbf{R}\Phi_{\mathscr{P}_a}(\mathscr{F})$.

*Notation* 1.2.1. In literature one may find the notation $\widehat{\mathbf{R}\Delta\mathscr{F}}$ instead of the more cumbersome $\widehat{\mathbf{R}a_*\mathbf{R}\Delta\mathscr{F}}$ that we adopted. But in the sequel we will often need to track down which map to $A$ we are employing, hence we will stick to this heavier notation.

We conclude this section by stating another important property of $GV$-sheaves. We recall the following definition

**Definition 1.2.7.** Let $k$ be a non-negative integer. A coherent sheaf $\mathscr{F}$ on a smooth projective variety is called and $k$-syzygy sheaf if it exists an exact sequence

$$0 \to \mathscr{F} \to \mathscr{E}_k \to \cdots \to \mathscr{E}_1 \to \mathscr{G} \to 0$$

with $\mathscr{G}$ another coherent sheaf and $\mathscr{E}_j$, $j = 1, \ldots, k$ locally free sheaves.

**Example 2.** a) Any coherent sheaf is a 0-syzygy sheaf.
b) $\mathscr{F}$ is 1-syzygy if and only if it is torsion free.
c) $\mathscr{F}$ is 2-syzygy if and only if it it reflexive.

In [64] Pareschi and Popa used $k$-syzygy sheaves to describe the classes $GV_k$ when $k$ is positive, proving the following.

**Theorem 1.2.8.** *The following are equivalent*

*(i)* $\mathrm{gv}_a(\mathscr{F}) \geq k$*;*

*(ii)* $\widehat{\mathbf{R}a_*\mathbf{R}\Delta\mathscr{F}}$ *is a $k$-syzygy sheaf.*

As an immediate consequence we get

**Theorem 1.2.9.** *Let $\mathscr{F}$ be a GV sheaf on $X$. The following are equivalent:*

*(i)* $\mathrm{gv}_a(\mathscr{F}) = 0$*;*

*(ii) the sheaf* $\widehat{\mathbf{R}a_*\mathbf{R}\Delta\mathscr{F}}$ *is not torsion free.*

Thank to Theorem 1.2.8, combined with Evans–Griffith syzygy theorem ([20]) Pareschi and Popa were able to prove the following relation between the Euler characteristic of a sheaf and its generic vanishing index.

**Proposition 1.2.10** ([63, Corollary 4.1])**.** *Let $X$ be a compact Khäler manifold of maximal Albanese dimension and $a : X \to A$ a generically finite morphism whose image generates $A$, then*

$$\chi(\omega_X) \geq \mathrm{gv}_a(\omega_X).$$

# Part I

# Syzygies of Projective Varieties

# CHAPTER 2

# SYZYGIES OF KUMMER VARIETIES

Let $X$ be an abelian variety. Its *associated Kummer variety* $\mathcal{K}_X$ is the quotient of $X$ by the natural $(\mathbb{Z}/2\mathbb{Z})$-action induced by the morphism $-1_X : X \to X$ defined by $x \mapsto -x$. Given a Kummer variety $\mathcal{K}_X$ and an ample line bundle $A$ on $\mathcal{K}_X$, a result of Sasaki ([68]) states that $A^{\otimes m}$ is very ample and the embedding it defines is projectively normal as soon as $m \geq 2$. Later Khaled ([45]) proved that, under the same conditions, the homogeneus ideal of $\mathcal{K}_X$ is generated by elements of degree 2 and 3, while, if $m \geq 3$ it is generated only by quadrics; if, furthermore, we assume the $A$ is a *general* very ample line bundle on $X$, then the homogeneus ideal of $\mathcal{K}_X$ is generated in degree less or equal 4 (for a complete exposition of existing results conerning syzygies of Kummer varieties see Section 2.1.2). In this Chapter we prove that these statements are particular cases of more general results on the syzygies of the variety $\mathcal{K}_X$.

More precisely, let $Z$ be an algebraic variety over an algebraically closed field $k$ and let $\mathscr{A}$ be an ample invertible sheaf on $Z$, generated by its global sections. With $R_{\mathscr{A}}$ we will indicate the *sections ring associated to the sheaf* $\mathscr{A}$:

$$R_{\mathscr{A}} := \bigoplus_{n \in \mathbb{Z}} H^0(Z, \mathscr{A}^{\otimes n})$$

while $S_{\mathscr{A}}$ will be the symmetric algebra of $H^0(Z, \mathscr{A})$. The ring $R_{\mathscr{A}}$ is a finitely generated graded $S_{\mathscr{A}}$-algebra and as such it admits a *minimal free resolution* $E_\bullet$, i.e. an exact complex

$$E_\bullet = 0 \to \cdots \xrightarrow{f_{p+1}} E_p \xrightarrow{f_p} \cdots \xrightarrow{f_2} E_1 \xrightarrow{f_1} E_0 \xrightarrow{f_0} R_{\mathscr{A}} \to 0 \qquad (2.0.1)$$

where

(i) $E_0 = S_{\mathscr{A}} \oplus \bigoplus_j S_{\mathscr{A}}(-a_{0j})$, $a_{0j} \in \mathbb{Z}$, $a_{0j} \geq 2$ since $Z$ is embedded by a complete linear system,

(ii) $E_i = \oplus_j S_{\mathscr{A}}(-a_{ij})$, $a_{ij} \in \mathbb{Z}$, $a_{ij} \geq 0$,

(iii) and $\mathrm{Im}(f_p) \subset \mathfrak{m}E_{p-1}$ where $\mathfrak{m}$ is the maximal homogeneus ideal of $S_{\mathscr{A}}$.

This resolution is unique in the sense that given $E'_\bullet$ another minimal free resolution there exists a graded isomorphism of complexes $E_\bullet \to E'_\bullet$ inducing the identity map on $R_{\mathscr{A}}$ (see [18, Theorem 20.1]).

In order to extend classical results of Castelnuovo, Mattuck, Fujita and Saint-Donat on the projective embeddings of curves, Green ([24]) introduced the following

**Definition** (Property $N_p$ )**.** Let $p$ be a given integer. The line bundle $\mathscr{A}$ satisfies property $N_p$ if, in the notations above,

$$E_0 = S_{\mathscr{A}}$$

and

$$E_i = \oplus S_{\mathscr{A}}(-i-1) \quad 1 \leq i \leq p.$$

Pareschi ([58]) extended the above condition as follows : we say that, given a non negative integer $r$ property $N_0^r$ holds for $\mathscr{A}$ if, in the notation above, $a_{0j} \leq 1 + r$ for every $j$ (i. e. the embedded variety is $h$-normal for every $h \geq 2+r$). Inductively we say that $\mathscr{A}$ satisfies property $N_p^r$ if $N_{p-1}^r$ holds for $\mathscr{A}$ and $a_{pj} \leq p + 1 + r$ for every $j$.

Green in [24] proved that, if $Z$ is a smooth curve of genus $g$ and $\mathscr{A}$ a very ample line bundle on $Z$ then $\mathscr{A}$ satisfies $N_p$ if $\deg \mathscr{A} \geq 2g + 1 + p$. He also conjectured that, if $C$ is a smooth non-hyperelliptic curve and $K_C$ is its canonical divisor, then

$$\mathscr{O}_C(K_C) \text{ satisfies } N_p \text{ if and only if } p < \mathrm{Cliff}(C);$$

where $\mathrm{Cliff}(C)$ is the Clifford index of the curve(cfr. [19, Section 9A]). Green's conjecture was recently proved for the general curve by Voisin ([73, 74]), and for the general cover by Aprodu–Farkas ([1]). Farkas investigated syzygies of curves in order to evince geometrical properties of the moduli sapces $\mathcal{M}_g$ and $\overline{\mathcal{M}_g}$ (see for example [21, 22]).

Another line of research started by [24], is to see how to best extend Green's results in higher dimension. The case of surfaces has been challenged by Gallego and Purnaprajna in a long series of articles starting in 1996. For a survey in this matter please see [23]. The syzygies of the projective space were studied by Green in [25], where he proved that $\mathscr{O}_{\mathbb{P}^n}(d)$ satisfies $N_p$ for $d \geq p \geq 1$, by Ottaviani–Paoletti ([57]), who proved that $\mathscr{O}_{\mathbb{P}^2}(d)$ does not satisfy $N_p$ for $3d-2 < p$, and in the recent preprint [17], where the authors studied the asymptotic behaviour of syzygies of projective varieties in general, and those of $\mathbb{P}^n$ in particular, demonstrating some theorems conjectured in [57]. For arbitrary smooth varieties there is a general conjecture of Mukai and in [15] Ein–Lazarsfeld proved that, if $Z$ is of dimension $n$, denoting by $\omega_Z$ its canonical line bundle on $Z$, then for any $\mathscr{L}$ very ample on $Z$ the sheaf

$$\mathscr{A} := \omega_Z \otimes \mathscr{L}^{\otimes(n+1+d)}$$

satisfies $N_p$ for every $d \geq p \geq 1$.

Abelian varieties distinguish themselves among other smooth varieties since, at least for what it concerns their syzygies, they tend to behave in any dimension like elliptic curves. More precisely, Koizumi ([47]) proved that, given $\mathscr{A}$ ample on an abelian variety $X$, then for $m \geq 3$, $\mathscr{A}^{\otimes m}$ embeds $X$ in the projective space as a projectively normal variety. Furthermore, a classical theorem of Mumford ([55]), perfectionated by Kempf ([41]), states that the homogeneus ideal of $X$ is generated in degree 2 as long as $m \geq 4$. These results inspired Lazarsfeld to conjecture that $\mathscr{A}^{\otimes m}$ satisfies $N_p$ for every $m \geq p+3$. In [42] Kempf proved that $\mathscr{A}^{\otimes m}$ satisfies condition $N_p$ as soon as $m \geq \max\{3, 2p+2\}$. A generalized version of Lazarsfeld's conjecture, involving property $N_p^r$ rather than simply $N_p$, was proved in [58]; later in [62] Pareschi–Popa were able to recover and improve Pareschi's statements as a consequence of the powerful, Fourier-Mukai based, theory of $M$-regularity that they developed in [61].

Given the results on projective normality and degree of defining equations of Sasaki and Khaled, and the close relationship between abelian varieties and Kummer varieties, it was natural to conjecture that even for the latter could be found a bound $m_0(p,r)$, *independent of the dimension of* $\mathcal{K}_X$ such that $\mathscr{A}^{\otimes m}$ satisfies $N_p^r$ for every $m \geq m_0(p,r)$. In this Chapter we present some results in this direction. The main idea behind the proofs is that ample line bundles on Kummer variety $\mathcal{K}_X$ have a nice description in terms of ample line bundles on $X$. More precisely, denoting by $\pi_X : X \to \mathcal{K}_X$ the quotient map, then for every $A$ ample on $\mathcal{K}_X$ it exists $\mathscr{A}$ ample on $X$ such that $\pi_X^* A \simeq \mathscr{A}^{\otimes 2}$. Hence we can use Pareschi–Popa machinery to find some results on $\mathscr{A}^{\otimes 2m}$ and then study how the $\mathbb{Z}/2\mathbb{Z}$ action fits in the frame. Below we list the main achievements we obtained.

**Theorem 2.A.** *Fix two non negative integers $p$ and $r$ such that* $\mathrm{char}(k)$ *does not divide $p+1$, $p+2$. Let $A$ be an ample line bundle on a Kummer variety $\mathcal{K}_X$, then*

*(a) $A^{\otimes n}$ satisfies property $N_p$ for every $n \in \mathbb{Z}$ such that $n \geq p+2$.*

*(b) More generally $A^{\otimes n}$ satisfies property $N_p^r$ for every $n$ such that $(r+1)n \geq p+2$.*

Since it consists in an improvement of existings results on the degree of defining equations of Kummer varieties it is worth to emphasize individually the case $p = 1$ of the above statement. Thank the geometric meaning of property $N_p^r$ (Section 2.1.1) one can deduce the following:

**Particular Case 2.B.** *Let $A$ be a very ample line bundle on a Kummer variety $\mathcal{K}_X$. Then the ideal of the image $\varphi_A(\mathcal{K}_X)$ in $\mathbb{P}(H^0(X,\, A))$ is generated by forms of degree at most 4.*

This result was classically known to be true, thank to the work of Wirtinger, Andreotti–Mayer and Khaled, just when $A$ was a *general* very ample line bundle on $\mathcal{K}_X$.

Adding one hypothesis about the line bundle $A$ we can get a somewhat better result improving the work of Kempf and Khaled; namely:

**Theorem 2.C.** *Let $p$ and $r$ be two integers such that $p \geq 1$, $r \geq 0$ and* char$(k)$ *does not divide $p+1$, $p+2$. Let $A$ be an ample line bundle on a Kummer variety $\mathcal{K}_X$, such that its pullback $\pi_X^* A \simeq \mathscr{A}^{\otimes 2}$ with $\mathscr{A}$ an ample symmetric invertible sheaf on $X$ which does not have a base divisor. Then*

*(a) $A^{\otimes n}$ satisfies property $N_p$ for every $n \in \mathbb{Z}$ such that $n \geq p+1$.*

*(b) More generally $A^{\otimes n}$ satisfies property $N_p^r$ for every $n$ such that $(r+1)n \geq p+1$.*

Again, it is worth of single out the case $p = 1$ of the above Theorem, concerning the equations of the Kummer variety $\mathcal{K}_X$.

**Particular Case 2.D.** *Suppose that* char$(k)$ *does not divide 2 or 3 and let $A$ be an ample invertible sheaf on $\mathcal{K}_X$ such that $\pi_X^* A \simeq \mathscr{A}^{\otimes 2}$ with $\mathscr{A}$ without a base divisor. Then*

*(a) If $n \geq 2$ then the ideal $\mathscr{I}_{\mathcal{K}_X, A^{\otimes n}}$ of the embedding $\varphi_{A^{\otimes n}}$ is generated by quadrics.*

*(b) $\mathscr{I}_{\mathcal{K}_X, A}$ is generated by quadrics and cubics.*

The key point of the proofs of Theorems 2.A and 2.C will be to reduce the problem on the Kummer variety $\mathcal{K}_X$ to a different problem on the abelian variety $X$. Namely we will show that property $N_p^r$ on the Kummer is implied by the surjectivity of a map of the type:

$$(*) \qquad \bigoplus_{[\alpha]\in\widehat{U}} H^0(X,\, \mathscr{F}\otimes\alpha)\otimes H^0(X,\, \mathscr{H}\otimes\alpha) \xrightarrow{m_\alpha} H^0(X,\, \mathscr{F}\otimes\mathscr{H}\otimes\alpha)$$

where $\mathscr{F}$ and $\mathscr{H}$ are sheaves on $X$ and $\widehat{U}$ is a non empty open subset of $\widehat{X}$, the abelian variety dual to $X$ and $m_\alpha$ is just the multiplication of global sections. Criteria for the surjectiviy of such maps are implicit in Kempf's work ([42, 43]), for the case $\mathscr{F}$ a vector bundle and $\mathscr{H}$ a line bundle (for an explicit argument due to Lazarsfeld see [58]). These results had been improved and extended to general coherent sheaves by Pareschi–Popa in [62].

This Chapter is organized in the following manner: in the next section we expone some background material such us the relationship between proprerty $N_p^r$ and the cohomology of the Koszul complex and a useful criterion for the surjectivity of a map of type (*). In Section 3 we present some slightly modified version of results of Sasaki and Khaled. The last section is entirely devoted to the proof of the main theorems.

## 2.1 Background Material

### 2.1.1 Property $\mathbf{N_p^r}$ and Koszul Cohomology

In this first section we review some well known relations between property $N_p$, or more generally property $N_p^r$, and the surjectivity of certain multiplication maps of sections of vector bundles.

Let $Z$ be a projective variety and let $L$ be an ample invertible sheaf on $Z$. We begin by stating a well known result of homological algebra:

**Proposition 2.1.1** ([19, Proposition I.1.7]). *Let $S = k[x_0, \ldots, x_n]$ a polynomial ring and let $\mathbf{E} : \cdots \to E_1 \to E_0$ be the minimal free resolution of a finitely generated $S$-module $M$. Then if $\mathcal{S}$ is any minimal set of homogeneus generators of $E_i$, then the set*

$$\mathcal{S}_j := \{s \in S \mid s \text{ has degree } j\} \subseteq \mathcal{S}$$

*has cardinality* $\dim_k \operatorname{Tor}_i^S(k, M)_j$.

An immediate corollary of this is

**Corollary 2.1.2.** *Fix $p$ and $r$ non negative integers. Given $L$ an a very ample line bundle on a projective variety $Z$, it satisfies property $N_p$ ($N_p^r$) if*

*(i)* $\operatorname{Tor}_0^{S_L}(k, R_L)_j = 0$ *for every $j \geq 1$ (for every $j \geq r+3$),*

*(ii)* $\operatorname{Tor}_p^{S_L}(k, R_L)_j = 0$ *for every $j \geq p+2$ (for every $j \geq r+p+2$).*

Using this fact and computing the above Tor groups via the Koszul resolution of the field $k$, Green observed (see for example, [24, Thm. 1.2], [26, Thm. 1.2], or [49, p. 511]) that condition $N_p$ is equivalent to te exactness in the middle of the complex

$$\bigwedge^{p+1} H^0(L) \otimes H^0(L^{\otimes h}) \to \bigwedge^{p} H^0(L) \otimes H^0(L^{\otimes h+1}) \to \bigwedge^{p-1} H^0(L) \otimes H^0(L^{\otimes h+2}) \tag{2.1.1}$$

for any $h \geq 1$. More generally, condition $N_p^r$ is equivalent to exactness in the middle of (2.1.1) for every $h \geq r+1$. Suppose that $L$ generated by its global sections and consider the following exact sequence:

$$0 \to M_L \longrightarrow H^0(L) \otimes \mathscr{O}_Z \longrightarrow L \to 0. \tag{2.1.2}$$

Taking wedge product, for any $p$ one gets the following exact sequence

$$0 \to \bigwedge^{p+1} M_L \longrightarrow \bigwedge^{p+1} H^0(L) \otimes \mathscr{O}_Z \longrightarrow \bigwedge^{p} M_L \otimes L \to 0. \tag{2.1.3}$$

It follows (for details cfr. [49] or [17]) that property $N_p^r$ is implied by the surjectivity of

$$\bigwedge^{p+1} H^0(L) \otimes H^0(L^{\otimes h}) \longrightarrow H^0 \left( \bigwedge^{p} M_L \otimes L^{\otimes h+1} \right) \tag{2.1.4}$$

for any $h \geq r+1$, where (2.1.4) was obtained by twisting (2.1.3) by $L^{\otimes h}$ and taking cohomology. Thus from (2.1.3) it follows that if

$$H^1(Z, \bigwedge^{p+1} M_L \otimes L^{\otimes h}) = 0, \tag{2.1.5}$$

then for any $h \geq r+1$, then condition $N_p$ is satisfied. If $\mathrm{char}(k)$ does not divide $p$, $\bigwedge^p \mathscr{E}$ is a direct summands of $\mathscr{E}^{\otimes p}$ for any vector bundle $\mathscr{E}$. Therefore we are led to the following Lemma:

**Lemma 2.1.3.** *Assume that* $\mathrm{char}(k)$ *does not divide* $p$ *and* $p+1$.

*(a) If* $H^1(Z, M_L^{\otimes p+1} \otimes L^{\otimes h}) = 0$ *for any* $h \geq r+1$ *then* $L$ *satisfies* $N_p^r$.

*(b) Let* $W \subseteq H^0(Z, L)$ *be a free sublinear system and denote by* $M_W$ *the kernel of the evaluation map* $W \otimes \mathscr{O}_Z \to L$. *Assume that* $H^1(Z, M_W^{\otimes p} \otimes L^{\otimes h}) = 0$, *then* $H^1(Z, M_W^{\otimes p+1} \otimes L^{\otimes h}) = 0$ *if and only if the multiplication map*

$$W \otimes H^0(M_W^{\otimes p} \otimes L^{\otimes h}) \longrightarrow H^0(M_W^{\otimes p} \otimes L^{\otimes h+1})$$

*is surjective.*

*Proof.* The proof of (a) is straightforward, while (b) follows from the following exact sequence:

$$0 \to M_W^{\otimes p+1} \otimes L^{\otimes h} \longrightarrow W \otimes M_W^{\otimes p} \otimes L^{\otimes h} \longrightarrow M_W^{\otimes p} \otimes L^{\otimes h+1} \to 0.$$

□

### Property $\mathrm{N_p^r}$ for Small p's'

By definition , if a variety $Z$ is embedded in a projective space by a very ample line bundle $L$ satisfying property $N_0^r$, then the variety $Z$ is $h$-normal for every $h \geq r$. Hence property $N_0$ is equivalent to projective normality.

For $p = 1$ property $N_p^r$ carries information of geometric nature since it returns intelligence about the equations of the embedding of the variety $Z$. More specifically, we can prove the following result that will allow us to deduce Particular Case 2.B and Particular Case 2.D from the main theorems.

**Proposition 2.1.4.** *If $L$ is a very ample line bundle on an algebraic variety $Z$ satisfying $N_0^r$, then the homogeneus ideal of $Z$ is generated by homogeneus elements of degree at most $r+2$.*

*Proof.* Denote by $V$ the vector space $H^0(Z,\, L)$ and let $\mathrm{S}^k V$ be the component of degree $k$ of the symmetric algebra of $V$, $S$. Consider furthermore the two $S$-modules

$$I = \bigoplus H^0(\mathbb{P}(V),\, \mathscr{I}_{Z,L}(k))$$

and $R_L$ and take a look to the following commutative diagram where the middle column is given by the Koszul complex.

$$\begin{array}{ccccccccc}
 & & & & \mathrm{S}^{k-1}V \otimes \bigwedge^2 V & \xrightarrow{(3)} & H^0(Z,\, M_L \otimes L^{\otimes k}) & & \\
 & & & & \downarrow & & \downarrow & & \\
0 & \longrightarrow & I_k \otimes V & \longrightarrow & \mathrm{S}^k V \otimes V & \longrightarrow & \{R_L\}_k \otimes V & \longrightarrow & 0 \\
 & & \downarrow{\scriptstyle (1)} & & \downarrow{\scriptstyle (2)} & & \downarrow & & \\
0 & \longrightarrow & I_{k+1} & \longrightarrow & \mathrm{S}^{k+1} V & \longrightarrow & \{R_L\}_{k+1} & \longrightarrow & 0
\end{array}$$

Our aim is to see that the map (1) is surjective for every $k \geq r+2$. Suppose that $L$ satisfies property $N_1^r$, then in particular property $N_0^r$ holds for $L$ and the second and third row are exact for every $k \geq r+1$. Since (2) is surjective, by the Snake Lemma, for every $k \geq r+1$ the surjectivity of (1) is implied by the surjectivity of (3) for every $k \geq r+2$. Now we can factor (3) in the following way:

$$\begin{array}{ccc}
\mathrm{S}^{k-1}V \otimes \bigwedge^2 V & \xrightarrow{\quad (3) \quad} & H^0(Z,\, M_L \otimes L^{\otimes k}) \\
 & {\scriptstyle g}\searrow \quad \nearrow{\scriptstyle f} & \\
 & H^0(Z,\, L^{\otimes k-1}) \otimes \bigwedge^2 V &
\end{array}$$

where $g$ is the canonical mapping $H^0(\mathbb{P}(V),\, \mathscr{O}_{\mathbb{P}}(k-1)) \longrightarrow H^0(Z,\, L^{\otimes k-1})$ and $f$ is the map in (2.1.4). For every $k \geq r+2$, $g$ is surjective because $L$ satisfies $N_0^r$, while the surjectivity of $f$ is equivalent to property $N_1^r$; hence (3) is surjective. $\square$

### 2.1.2 Kummer Varieties: Definition, Projective Normality and Equations

Let $X$ be an abelian variety over a field $k$, with $\mathrm{char}(k) \neq 2$. As usual we denote by $-1_X : X \longrightarrow X$ the morphism given by $x \mapsto -x$. We recall that the *Kummer variety associated to* $X$, denoted with $\mathcal{K}_X$, is the quotient variety

$$\mathcal{K}_X := X / < -1_X > .$$

By $\pi_X : X \longrightarrow \mathcal{K}_X$ we will mean the canonical map to the quotient.

### 2.1.3 Line Bundles on Kummer Varieties

In this paragraph we recall some basic fact on symmetric sheaves and on line bundles on Kummer varieties. This part of the thesis, as well as any other passage in this section, is of expository nature: a complete treaty of the results here presented can be found in [54, pp. 303-305].

We recall that an invertible sheaf $\mathscr{L}$ on an abelian variety $X$ is called *symmetric* when

$$(-1_X)^*\mathscr{L} \simeq \mathscr{L}$$

Thus, take $\mathscr{L}$ a symmetric line bundle on $X$ and fix $\psi$ an isomorphism

$$\psi : \mathscr{L} \xrightarrow{\sim} (-1_X)^*\mathscr{L};$$

then for all $x \in X$ closed points, $\psi$ induces an isomorphism

$$\psi(x) : \mathscr{L}(x) \to \mathscr{L}(-x).$$

Therefore it is possible to *canonically normalize* $\psi$ by requesting that the map $\psi(0)$ is the identity.

**Definition 2.1.5.** The canonical isomorphism arised from the above construction is denoted by $\psi_{\mathscr{L}}$ and it is called *Mumford's normalized isomorphism of* $\mathscr{L}$

Remark that for every $x \in X$ point of order 2, $\psi_{\mathscr{L}}$ induce an involution $e(x)$ on $\mathscr{L}(x)$.

**Definition 2.1.6.** The line bundle $\mathscr{L}$ is *totally symmetric* if $e(x)$ is the identity for every 2-torsion point.

**Example 3.** An even power of a symmetric line bundle is always totally symmetric.

A converse to this example holds for ample line bundles:

**Proposition 2.1.7.** *Let* $\mathscr{L}$ *be a totally symmetric ample invertible sheaf on an abelian variety* $X$*, then* $\mathscr{L} \simeq \mathscr{A}^{\otimes 2}$ *with* $\mathscr{A}$ *an ample symmetric invertible sheaf on* $X$*.*

This statement is classical and well known. A proof of its is implicit in [68, Proof of Lemma 1.2]. We include our own proof for the reader's benefit. We will need the following two statements about totally symmetric line bundles.

**Lemma 2.1.8** (Properties of totally symmetric line bundles). *(i) Let $\mathscr{L}$ be an ample totally symmetric invertible sheaf of type $\delta = (d_1, \dots, d_g)$ on an abelian variety $X$ of dimension $g$. The group*

$$K(\mathscr{L}) := \operatorname{Ker}\left(x \mapsto t_x^*\mathscr{L} \otimes \mathscr{L}^{-1}\right)$$

*contains all the points of order 2 of $X$. Hence all the $d_i$'s are even.*

*(ii) If $\mathscr{L}_1$ and $\mathscr{L}_2$ are two totally symmetric line bundles such that they are algebraically equivalent, then $\mathscr{L}_1 \simeq \mathscr{L}_2$.*

*Proof.* Part (i) is Corollary 4 p. 310 in [54] while part (ii) is explained at p. 307 of the same paper. □

*Proof of Proposition 2.1.7.* By the Lemma 2.1.8(i) we have that $\mathscr{L} \simeq \mathscr{M}^{\otimes 2}$ with $\mathscr{M}$ an ample line bundle on $X$. If $\mathscr{M}$ is already symmetric there is nothing to prove. Thus we can suppose that $\mathscr{M}$ is not symmetric. Now consider $\alpha$ a topologically trivial line bundle such that $\alpha^{\otimes 2} \simeq (-1_X)^*\mathscr{M} \otimes \mathscr{M}^{\otimes -1}$ (it exists because $\operatorname{Pic}^0(X)$ is an abelian variety and therefore a divisible group ([55, (iv) p. 42])). The invertible sheaf $\mathscr{A} := \mathscr{M} \otimes \alpha$ is symmetric. In fact we have

$$\mathscr{A} \otimes (-1_X)^*\mathscr{A}^{\otimes -1} \simeq \mathscr{M} \otimes \alpha \otimes (-1_X)^*(\mathscr{M}^{\otimes -1} \otimes \alpha^{\otimes -1}) \simeq \alpha^{\otimes 2} \otimes \mathscr{M} \otimes (-1_X)^*\mathscr{M}^{\otimes -1} \simeq \mathscr{O}_X.$$

Therefore the sheaf $\mathscr{A}^{\otimes 2}$ is totally symmetric and algebraically equivalent to $\mathscr{L}$. Since also $\mathscr{L}$ is totally simmetric the statement follows for Lemma 2.1.8(ii). □

**Proposition 2.1.9** ([54, Proposition 1 page 305]). *Let $\mathscr{L}$ be an invertible sheaf on an abelian variety $X$ and consider the associated Kummer vairety $\pi_X : X \to \mathcal{K}_X$. Then $\mathscr{L}$ is of the form $\pi_X^* M$ with $M$ some line bunlde on $\mathcal{K}_X$ if and only if it is totally symmetric.*

Now take $\mathscr{L}$ a symmetric invertible sheaf on $X$. The $\mathbb{Z}/2\mathbb{Z}$ action on $X$ given by the involution $-1_X : X \longrightarrow X$ induces trough $\psi_{\mathscr{L}}$ a lifting of the action on $\mathscr{L}$. The composition

$$H^0(X,\, \mathscr{L}) \xrightarrow{(-1_X)^*} H^0(X,\, (-1_X)^*\mathscr{L}) \xrightarrow{(-1_X)^*(\psi_{\mathscr{L}})} H^0(X,\, \mathscr{L})$$

is denoted by $[-1]_{\mathscr{L}}$, or simply by $[-1]$ when there is no chance of misinterpretation, and it is an involution of $H^0(X,\, \mathscr{L})$ and admits just 1 and -1 as eigenvalues. We define

$$H^0(X,\, \mathscr{L})^{\pm} = \{s \in H^0(X,\, \mathscr{L}) \text{ such that } [-1]_{\mathscr{L}} s = \pm s\}$$

If $\mathscr{L}$ is totally symmetric, by Proposition 2.1.9 then it exists a line bundle $M$ on the Kummer variety $\mathcal{K}_X$ such that $\pi_X^* M \simeq \mathscr{L}$ and one can identify $H^0(\mathcal{K}_X, M)$ with $H^0(X, \mathscr{L})^+$.

*Notation* 2.1.1. From now till the end of this chapter, in order to lighten the notation, we will use the letter $i$ to denote the map $(-1)_X$.

### 2.1.4 Projective Normality and Equations

If $A$ is an ample invertible sheaf on $\mathcal{K}_X$ generated by its global sections, one may wonder about the "good" properties enjoyed by the morphism

$$\varphi_{|A^{\otimes n}|} : \mathcal{K}_X \longrightarrow \mathbf{P}(H^0(\mathcal{K}_X, A^n))$$

associated to the complete linear system $|A^n|$, $n \in \mathbb{N}^*$. In this paragraph we review some known results on very ampleness, projective normality and bound of the degree of equations of Kummer varieties.

For example, a well known fact (cfr. [5, Proposition 4.8.1]) says that, if $\pi_X^* A = 2\theta$ with $\theta$ a principal polarization on $X$, then $A$ is very ample, i.e. $\varphi_{|A|}$ is an embedding. For what it concerns projective normality, we have the following results.

**Theorem** (Sasaki, Khaled)**.** *In the above notation*

*(ii)* *$A^{\otimes n}$ is very ample and normally generated, (i. e. the embedding induced by $|A^{\otimes n}|$ is projectively normal) for every $n \geq 2$;*

*(i)* *If $A$ is very ample, then it is normally generated if and only if, after writing $\pi^* A \simeq \mathscr{A}^{\otimes 2}$ with $\mathscr{A}$ ample and symmetric on $X$, $0_X \notin \mathrm{Bs}(H^0(X, \mathscr{A} \otimes \alpha)^+)$ for every $[\alpha] \in \widehat{X}$ of order 2.*

Th first part of the above Theorem is due to Sasaki ([68]), while the second was proved by Khaled in [46]. The degree of generators of the homogeneous ideals of $\mathcal{K}_X$ was studied by Kempf and Khaled who proved the following statements.

**Theorem.** *In the notation above.*

*(i)* *Let $\theta$ and $A$ be respectively a principal polarization on the abelian variety $X$, and a line bundle on $\mathcal{K}_X$ such that $\pi_X^* A \simeq 2\theta$. Then $A^{\otimes 2}$ is normally generated and the homogeneous ideal associated to the embedding $\varphi_{\mathscr{A}^{\otimes 2}}$ is generated by its components of degree two and three.*

*(ii)* *If $n \geq 3$, then the ideal of $\mathcal{K}_X$ in the embedding given by $A^{\otimes n}$ is generated by forms of degree 2.*

*(iii) The image homogeneus ideal associated to* $\varphi_{A^{\otimes 2}}$ *is generated in degree two and three.*

*(iv) If A very ample and normally generated, then the ring*

$$\bigoplus_{n\in\mathbb{N}} H^0(\mathcal{K}_X,\, A^{\otimes n})$$

*is generated by* $H^0(\mathcal{K}_X,\, A)$ *modulo quadric, cubic, and quartic relations.*

*Proof.* The first part of the statement is proved in [44]. For the others please see [45]. □

In [4] Beauville proved that part (i) of the above Theorem is sharp:

**Proposition** ([4, Proposition 3.2])**.** *Let* $(X, \Theta)$ *be an indecomposable principally polarized abelin variety and consider* $\varphi_{|2\Theta|}$ *the morphism associated to the complete linear system* $|2\Theta|$. *The ideal* $\mathscr{I}_{|2\Theta|}$ *of* $\varphi_{|2\Theta|}(X)$ *in* $\mathbb{P}(H^0(X, \mathscr{O}_X(2\Theta))$ *cannot be genereted by its elements of degree* $\leq 3$.

### 2.1.5 M-regular Sheaves and Multiplication Maps

$M$-regular sheaves and $M$-regularity theory, introduced by Pareschi-Popa and reviewed in the first Chapter of this thesis, are crucial to our purpose thank to their application in determining whether a map of the form

$$\bigoplus_{[\alpha]\in U} H^0(X,\, \mathscr{F}\otimes\alpha)\otimes H^0(X,\, \mathscr{H}\otimes\alpha^\vee) \xrightarrow{m_\alpha} H^0(X,\, \mathscr{F}\otimes\mathscr{H}), \tag{2.1.6}$$

with $\mathscr{F}$ and $\mathscr{H}$ sheaves on an abelian variety $X$ and $U \subseteq \widehat{X}$ an open set, is surjective. We will list below all the results of such kind that we will be using troughout the paper. The first one is an extension of a theorem that had already appeared in the work of Kempf, Mumford and Lazarsfeld.

**Theorem 2.1.10** ([61], Theorem 2.5)**.** *Let* $\mathscr{F}$ *and* $\mathscr{H}$ *be sheaves on* $X$ *such that* $\mathscr{F}$ *is* $M$*-regular and* $\mathscr{H}$ *is locally free satisfying I.T. with index 0. Then* (2.1.6) *is surjective for any non empty Zariski open set* $U \subseteq \widehat{X}$.

An easy corollary of the above result is stated below.

**Corollary 2.1.11.** *If* $\mathscr{F}$ *and* $\mathscr{H}$ *satisfy the hypothesis of the above Theorem, then there exists* $N$ *a positive integer such for the general* $[\alpha_1], \dots, [\alpha_N] \in \widehat{X}$ *the map*

$$\bigoplus_{k=1}^{N} H^0(X,\, \mathscr{F}\otimes\alpha_k)\otimes H^0(X,\, \mathscr{H}\otimes\alpha_k^\vee) \xrightarrow{\quad m_{\alpha_k}\quad} H^0(X,\, \mathscr{F}\otimes\mathscr{H})$$

*is surjective.*

We conclude this paragraph by presenting two results on multiplication maps of sections we will be needing afterward.

**Proposition 2.1.12.** *Let $\mathscr{A}$ be an ample line bundle on an abelian variety $X$. The map*

$$m : H^0(X, \mathscr{A}^{\otimes 2}) \otimes H^0(X, \mathscr{A}^{\otimes 2} \otimes \alpha) \longrightarrow H^0(X, \mathscr{A}^{\otimes 4} \otimes \alpha). \tag{2.1.7}$$

*is surjective for the general $[\alpha] \in \widehat{X}$. If furthermore $\mathscr{A}$ does not have a base divisor, then the locus $Z \subseteq \mathrm{Pic}^0(X)$ in which it fails to be surjective has codimension at least 2.*

*Proof.* The first part of the statement is classical, for a reference see, for example [5, Proposition 7.2.2]. For what it concerns the second part, it was proved by Pareschi–Popa ([62]) as a consequence of their $M$-regularity techniques. □

## 2.2 Multiplication Maps on Abelian Varieties

A well known result by Khaled states that

**Proposition 2.2.1** ([45])**.** *Let $\mathscr{A}$ be an ample symmetric vector bundle on an abelian variety $X$. Take $k = 2n$ an even positive integer. Thus $\mathscr{A}^{\otimes k}$ is totally symmetric and for every $n \geq 1$ and every $h \in \mathbb{Z}$, $h \geq 3$ the following map is surjective for every $[\alpha] \in \mathrm{Pic}^0(X)$*

$$m_\alpha^+ : H^0(X, \mathscr{A}^{\otimes k})^+ \otimes H^0(X, \mathscr{A}^{\otimes h} \otimes \alpha) \longrightarrow H^0(X, \mathscr{A}^{\otimes k+h} \otimes \alpha).$$

The main goal of this section is to prove that the same is true for every $h = 2$ and for *general* $\alpha \in \widehat{X}$. If furthermore we assume that $\mathscr{A}$ does not have a base divisor, then we will show that the locus of $[\alpha] \in \widehat{X}$ where

$$m_\alpha^+ : H^0(X, \mathscr{A}^{\otimes 2n})^+ \otimes H^0(X, \mathscr{A}^{\otimes 2} \otimes \alpha) \longrightarrow H^0(X, \mathscr{A}^{\otimes 2(n+1)} \otimes \alpha) \tag{2.2.1}$$

fails to be surjective has codimension at least 2. We will do this by slightly modifying the methods adopted by Khaled in [45].

To this end, consider the isogeny

$$\xi : X \times X \longrightarrow X \times X$$

given by $\xi = (p_1 + p_2,\ p_1 - p_2)$

**Lemma 2.2.2.** *For any $[\alpha] \in \widehat{X}$ we have an isomorphism*

$$\xi^*(p_1^*(\mathscr{A} \otimes \beta)) \otimes p_2^*(\mathscr{A} \otimes \alpha) \xrightarrow{\Phi_\alpha^\beta} p_1^*(\mathscr{A}^{\otimes 2} \otimes \beta \otimes \alpha) \otimes p_2^*(\mathscr{A}^{\otimes 2} \otimes \beta \otimes \alpha^\vee)$$

*Proof.* We will use the "See-saw" Principle: for any $y \in X$ we have

$$\xi^*(p_1^*(\mathscr{A} \otimes \beta) \otimes p_2^*(\mathscr{A} \otimes \alpha))_{|X \times \{y\}} \simeq t_y^*\mathscr{A} \otimes t_{-y}^*\mathscr{A} \otimes t_y^*\alpha^\vee \otimes t_{-y}^*\alpha \simeq$$
$$\simeq \mathscr{A}^{\otimes 2} \otimes \beta \otimes \alpha.$$

Now we look at the restirction of $\xi^*(p_1^*(\mathscr{A} \otimes \beta) \otimes p_2^*(\mathscr{A} \otimes \alpha))$ to $\{0\} \times X$ and we get

$$\xi^*(p_1^*(\mathscr{A} \otimes \beta) \otimes p_2^*(\mathscr{A} \otimes \alpha))_{|\{0\} \times X} \simeq \mathscr{A} \otimes \beta \otimes i^*\mathscr{A} \otimes \alpha) \simeq \mathscr{A}^{\otimes 2} \otimes \beta \otimes \alpha^\vee;$$

hence the statement is proved. □

*Notation* 2.2.1. When $[\alpha] = [\beta^\vee]$ we will denote $\Phi_\alpha^\beta$ simply by $\Phi_\beta$

Composing with the Künneth isomorphism we have a map

$$\xi^* : H^0(X,\, \mathscr{A} \otimes \beta) \otimes H^0(X,\, \mathscr{A} \otimes \alpha) \longrightarrow H^0(X,\, \mathscr{A}^{\otimes 2} \otimes \beta \otimes \alpha) \otimes H^0(X,\, \mathscr{A}^{\otimes 2} \otimes \beta \otimes \alpha^\vee).$$

Taking $[\alpha] = [\beta^\vee]$ we want to characterize the image of $H^0(X, \mathscr{A} \otimes \beta) \otimes H^0(X, \mathscr{A} \otimes \beta^\vee)$ in $H^0(X,\, \mathscr{A}^{\otimes 2})^+ \otimes H^0(X,\, \mathscr{A}^{\otimes 2} \otimes \beta^{\otimes 2})$ through $\xi^*$.

In order to achieve this goal, we consider the following two automorphisms of $X \times X$

- ◇ $i_{\mathrm{L}} = (-p_1, p_2)$. It induces the automorphism $[-1]_{\mathscr{A}^{\otimes 2}} \otimes \mathrm{id}$ on

  $$H^0(X,\, \mathscr{A}^{\otimes 2})^+ \otimes H^0(X,\, \mathscr{A}^{\otimes 2} \otimes \beta^{\otimes 2});$$

  for every section $s \in H^0(X,\, \mathscr{A}^{\otimes 2}) \otimes H^0(X,\, \mathscr{A}^{\otimes 2} \otimes \beta^{\otimes 2})$, we have that

  $$s \in H^0(X,\, \mathscr{A}^{\otimes 2})^+ \otimes H^0(X,\, \mathscr{A}^{\otimes 2} \otimes \beta^{\otimes 2}) \quad \Leftrightarrow \quad [-1]_{\mathscr{A}^{\otimes 2}} \otimes \mathrm{id}(s) = s.$$

- ◇ $\widehat{\tau}$, the automorphism of $X \times X$ defined by $(-p_2, -p_1)$.

Now, for every $[\beta] \in \widehat{X}$ let

$$\psi_\beta : \mathscr{A} \otimes \beta \longrightarrow i^*\mathscr{A} \otimes \beta \tag{2.2.2}$$

be the isomorphism given by tensoring the normalized isomorhims of $\mathscr{A}$, $\psi_{\mathscr{A}}$, with the identity of $\beta$. We denote by

$$\widehat{\nu} : p_2^* i^*(\mathscr{A} \otimes \beta) \otimes p_1^* i^*(\mathscr{A} \otimes \beta^\vee) \longrightarrow: p_1^*(\mathscr{A} \otimes \beta) \otimes p_2^*(\mathscr{A} \otimes \beta^\vee)$$

the isomorphism of sheaves defined by $\widehat{\nu}(p_2^* i^* t \otimes p_1^* i^* s) = p_1^* i^*(\psi_\beta) i^* s \otimes p_2^* i^*(\psi_{\beta^\vee}) i^* t$. We obtain the following diagram

$$\begin{array}{ccc}
H^0(X \times X,\, p_1^*(\mathscr{A} \otimes \beta) \otimes p_2^*(\mathscr{A} \otimes \beta^\vee)) & \xrightarrow{\text{Künneth}} & H^0(X,\, \mathscr{A} \otimes \beta) \otimes H^0(X,\, \mathscr{A} \otimes \beta^\vee) \\
\downarrow{\scriptstyle \widehat{\tau}} & & \\
H^0(X \times X,\, p_2^* i^*(\mathscr{A} \otimes \beta) \otimes p_1^* i^*(\mathscr{A} \otimes \beta^\vee)) & & \downarrow{\scriptstyle \widehat{T}^{\mathscr{A} \otimes \beta}} \\
\downarrow{\scriptstyle \widehat{\nu}} & & \\
H^0(X \times X,\, p_1^*(\mathscr{A} \otimes \beta) \otimes p_2^*(\mathscr{A} \otimes \beta^\vee)) & \xrightarrow{\text{Künneth}} & H^0(X,\, \mathscr{A} \otimes \beta) \otimes H^0(X,\, \mathscr{A} \otimes \beta^\vee)
\end{array}$$

Where $\widehat{T}^{\mathscr{A}\otimes\beta}$ is the involution defined by $\widehat{T}^{\mathscr{A}\otimes\beta}(s\otimes t) = i^*(\psi_{\beta^\vee})i^*t \otimes i^*(\psi_\beta)i^*s$. Let us denote by

$$[H^0(X,\, \mathscr{A}\otimes\beta)\otimes H^0(X,\, \mathscr{A}\otimes\beta^\vee)]^\pm$$

the eigenspaces of

$$H^0(X,\, \mathscr{A}\otimes\beta)\otimes H^0(X,\, \mathscr{A}\otimes\beta^\vee)$$

under the action of $\widehat{T}^{\mathscr{A}\otimes\beta}$.

**Proposition 2.2.3.** *For every* $[\alpha]\in\widehat{X}$ *we have*

$$\xi^*[H^0(X,\, \mathscr{A}\otimes\beta)\otimes H^0(X,\, \mathscr{A}\otimes\beta^\vee)]^\pm \subseteq H^0(X,\, \mathscr{A}^{\otimes 2})^\pm \otimes H^0(X,\, \mathscr{A}^{\otimes 2}\otimes\beta^{\otimes 2}).$$

*Proof.* We shall write $H^0(\mathscr{A})$ instead of $H^0(X,\mathscr{A})$. This proof is just a slight modification of Khaled's proof of [46, Proposition 2.2]. First of all observe that the following diagram commutes:

$$\begin{array}{ccc}
p_1^*(\mathscr{A}^{\otimes 2})\otimes p_2^*(\mathscr{A}^{\otimes 2}\otimes\beta^{\otimes 2}) & \xleftarrow{\;p_1^*(\psi_{\mathscr{A}^{\otimes 2}})\otimes\mathrm{id}\;} & p_1^*i^*(\mathscr{A}^{\otimes 2})\otimes p_2^*(\mathscr{A}^{\otimes 2}\otimes\beta^{\otimes 2}) \\
\uparrow{\scriptstyle \Phi_\beta} & & \uparrow{\scriptstyle i_L^*(\Phi_\beta)} \\
\xi^*(p_1^*(\mathscr{A}\otimes\beta)\otimes p_2^*(\mathscr{A}\otimes\beta^\vee)) & \xleftarrow{\;\xi^*(\widehat{\nu})\;} & \xi^*(\widehat{\tau}^*(p_1^*(\mathscr{A}^\otimes\beta)\otimes p_2^*(\mathscr{A}\otimes\beta^\vee)))
\end{array}$$

In fact it certainly commutes up to a costant. By looking at the diagram in the origin of $X\times X$ this costant can be shown to be 1. Thus we have the commutative diagram

$$\begin{array}{ccccc}
H^0(\mathscr{A}^{\otimes 2}\boxtimes(\mathscr{A}^{\otimes 2}\otimes\beta^{\otimes 2})) & \xleftarrow{p_1^*(\psi_{\mathscr{A}^{\otimes 2}})\otimes\mathrm{id}} & H^0(i^*\mathscr{A}^{\otimes 2}\boxtimes(\mathscr{A}^{\otimes 2}\otimes\beta^{\otimes 2})) & \xleftarrow{i_L^*} & H^0(\mathscr{A}^{\otimes 2}\boxtimes(\mathscr{A}^{\otimes 2}\otimes\beta^{\otimes 2})) \\
\uparrow{\scriptstyle\Phi_\beta} & & \uparrow{\scriptstyle i_L^*\Phi_\beta} & & \uparrow{\scriptstyle\Phi_\beta} \\
H^0(\xi^*((\mathscr{A}\otimes\beta)\boxtimes(\mathscr{A}\otimes\beta^\vee)) & \xleftarrow{\xi^*\widehat{\nu}} & H^0(\xi^*\widehat{\tau}^*((\mathscr{A}\otimes\beta)\boxtimes(\mathscr{A}\otimes\beta^\vee)) & \xleftarrow{i_L^*} & H^0(\xi^*((\mathscr{A}\otimes\beta)\boxtimes(\mathscr{A}\otimes\beta^\vee)) \\
\uparrow{\scriptstyle\xi^*} & & \uparrow{\scriptstyle\xi^*} & & \uparrow{\scriptstyle\xi^*} \\
H^0((\mathscr{A}\otimes\beta)\boxtimes(\mathscr{A}\otimes\beta^\vee)) & \xleftarrow[\widehat{\nu}]{} & H^0(\widehat{\tau}^*((\mathscr{A}\otimes\beta)\boxtimes(\mathscr{A}\otimes\beta)) & \xleftarrow[\widehat{\tau}]{} & H^0((\mathscr{A}\otimes\beta)\boxtimes(\mathscr{A}\otimes\beta^\vee))
\end{array}$$

Hence, composing with Künneth isomorphism we obtain another commutative diagram:

$$\begin{array}{ccc}
H^0(\mathscr{A}^{\otimes 2})\otimes H^0(\mathscr{A}^{\otimes 2}\otimes\beta^{\otimes 2}) & \xleftarrow{\;\xi^*\;} & H^0(\mathscr{A}\otimes\beta)\otimes H^0(\mathscr{A}\otimes\beta^\vee) \\
{\scriptstyle p_1^*(\psi_{\mathscr{A}^{\otimes 2}})\otimes\mathrm{id}}\uparrow & & \\
H^0(\mathscr{A}^{\otimes 2})\otimes H^0(\mathscr{A}^{\otimes 2}\otimes\beta^{\otimes 2}) & & \uparrow{\scriptstyle\widehat{T}^{\mathscr{A}\otimes\beta}} \\
{\scriptstyle i_L^*}\uparrow & & \\
H^0(\mathscr{A}^{\otimes 2})\otimes H^0(\mathscr{A}^{\otimes 2}\otimes\beta^{\otimes 2}) & \xleftarrow{\;\xi^*\;} & H^0(\mathscr{A}\otimes\beta)\otimes H^0(\mathscr{A}\otimes\beta^\vee)
\end{array}$$

Therefore we have

$$\xi^*\circ\widehat{T}^{\mathscr{A}\otimes\beta} = [-1]_{\mathscr{A}^{\otimes 2}}\otimes\mathrm{id}\circ\xi^*.$$

The statement follows directly. $\square$

**Theorem 2.2.4.** *Let $\mathscr{A}$ be an ample symmetric line bundle on $X$ and take $[\alpha] \in \widehat{X}$. Then the multiplication map*

$$m : H^0(X,\, \mathscr{A}^{\otimes 2}) \otimes H^0(X,\, \mathscr{A}^{\otimes 2} \otimes \alpha) \longrightarrow H^0(X,\, \mathscr{A}^{\otimes 4} \otimes \alpha)$$

*is surjective if and only if the following multiplication map is surjective*

$$m^+ : H^0(X,\, \mathscr{A}^{\otimes 2})^+ \otimes H^0(X,\, \mathscr{A}^{\otimes 2} \otimes \alpha) \longrightarrow H^0(X,\, \mathscr{A}^{\otimes 4} \otimes \alpha)$$

*Proof.* The "if" part is strightforward. The proof of the "only if" part is a mix of Obuchi's proof of his Theorem about normal generation on abelian varieties (cfr. [56, Theorem]) and Khaled's proof of his result about normal generation on Kummer varieties (cfr. [46, Theorem 2.3]).
The multiplication map

$$m : H^0(X,\, \mathscr{A}^{\otimes 2}) \otimes H^0(X,\, \mathscr{A}^{\otimes 2} \otimes \alpha) \longrightarrow H^0(X,\, \mathscr{A}^{\otimes 4} \otimes \alpha)$$

is (as observed by Khaled) the composition of $\xi^*$ with the map $\mathrm{id} \otimes \mathrm{e}_{\mathscr{A}^{\otimes 4} \otimes \alpha^\vee}$, where $\mathrm{e}_{\mathscr{A}^{\otimes 4} \otimes \alpha^\vee}$ is the evaluation in 0 of the sections of $\mathscr{A}^{\otimes 4} \otimes \alpha^\vee$. Infact if we denote by $\Delta : X \to X \times X$ the diagonal immersion, then $m$ is just $\Delta^*$ composed with the Künneth isomorphism. Now we can write $\Delta = \xi \circ f$ where $f : X \to X \times X$ is the morphism defined by $x \mapsto (x, 0_X)$. Now observe that, modulo Künneth isomorhism,

$$f^* : H^0(X, \mathscr{A}^{\otimes 4}) \otimes H^0(X, \mathscr{A}^{\otimes 4} \otimes \alpha^\vee) \to H^0(X, \mathscr{A}^{\otimes 4})$$

is exactly $\mathrm{id} \otimes \mathrm{e}_{\mathscr{A}^{\otimes 4} \otimes \alpha^\vee}$. Hence $m = \Delta^* = \mathrm{id} \otimes \mathrm{e}_{\mathscr{A}^{\otimes 4} \otimes \alpha^\vee} \circ \xi^*$

Thus we can consider the following commutative diagram

$$\begin{array}{ccc}
\bigoplus_{[\beta^{\otimes 2}]=[\alpha]} H^0(\mathscr{A} \otimes \beta) \otimes H^0(\mathscr{A} \otimes \beta^\vee) & \xrightarrow{\quad \xi^* \quad} & H^0(\mathscr{A}^{\otimes 2}) \otimes H^0(\mathscr{A}^{\otimes 2} \otimes \alpha) \\
 & {\scriptstyle 2_X^* \otimes \mathrm{e}_{\mathscr{A} \otimes \beta^\vee}} \searrow \qquad \swarrow {\scriptstyle m} & \\
 & H^0(\mathscr{A}^{\otimes 4} \otimes \alpha) &
\end{array}$$

The upper arrow is an isomorphism by projection formula. In fact we have that

$$\xi_* \mathscr{O}_{X \times X} \simeq \bigoplus_{[\gamma] \in \widehat{X}_2} p_1^* \gamma \otimes p_2^* \gamma.$$

Then we can write

$$\begin{aligned}
H^0(\mathscr{A}^{\otimes 2}) \otimes H^0(\mathscr{A}^{\otimes} 2 \otimes \alpha) &\simeq H^0(p_1^*(\mathscr{A}^{\otimes 2}) \otimes p_2^*(\mathscr{A}^{\otimes 2} \otimes \alpha)) \simeq \\
&\simeq H^0(\xi^*(p_1^*(\mathscr{A} \otimes \beta) \otimes p_2^*(\mathscr{A} \otimes \beta^\vee))) \simeq \\
&\simeq H^0(\xi_* \xi^*(p_1^*(\mathscr{A} \otimes \beta) \otimes p_2^*(\mathscr{A} \otimes \beta^\vee))) \simeq \\
&\simeq H^0(p_1^*(\mathscr{A} \otimes \beta) \otimes p_2^*(\mathscr{A} \otimes \beta^\vee) \otimes \xi_* \mathscr{O}_{X \times X}) \simeq \\
&\simeq \bigoplus_{[\beta^{\otimes 2}]=[\alpha]} H^0(\mathscr{A} \otimes \beta) \otimes H^0(\mathscr{A} \otimes \beta^\vee).
\end{aligned}$$

Therefore the surjectivity of $m$ is equivalent to the following

$$(\dagger) \qquad 0 \notin \mathrm{Bs}(\mathscr{A} \otimes \beta^{\vee}) \text{ for every } \beta \in \widehat{X} \text{ such that } \beta^2 \simeq \alpha,$$

where, as usual $\mathrm{Bs}(\mathscr{A} \otimes \beta^{\vee})$ stands for the base locus of $\mathscr{A} \otimes \beta^{\vee}$.

Now we shall prove that if $(\dagger)$ holds, then $m^+$ is surjective. Thank to the Proposition 2.2.3 and the isomorphism

$$\bigoplus_{[\beta^2]\simeq[\alpha]} H^0(\mathscr{A} \otimes \beta) \xrightarrow{2_X^*} H^0(\mathscr{A}^{\otimes 4} \otimes \alpha),$$

yielded by $2_X^*$ and projection formula, it is enough to check that for every $[\beta]$ satisfying $[\beta^{\otimes 2}] = [\alpha]$

$$2_X^* \otimes \mathrm{e}_{\mathscr{A}\otimes\beta^{\vee}}(H^0(\mathscr{A} \otimes \beta) \otimes H^0(\mathscr{A} \otimes \beta^{\vee}))^+ = 2_X^* H^0(\mathscr{A} \otimes \beta).$$

To this goal take $s \in H^0(\mathscr{A} \otimes \beta)$, we want to provide an element

$$\sigma \in (H^0(\mathscr{A} \otimes \beta) \otimes H^0(\mathscr{A} \otimes \beta^{\vee}))^+$$

such that

$$2_X^* \otimes \mathrm{e}_{\mathscr{A}\otimes\beta^{\vee}}(\sigma) = 2_X^*(s)$$

Denote by $\lambda$ the constant $\mathrm{e}_{\mathscr{A}\otimes\beta^{\vee}}([i^*\psi_\beta \circ i^*](s))$. If $\lambda \neq 0$, take

$$\sigma := \frac{1}{\lambda} \cdot (s \otimes [i^*\psi_\beta \circ i^*](s)).$$

Suppose, otherwise, that $\lambda = 0$; since $(\dagger)$ holds, it exist a $t \in H^0(\mathscr{A} \otimes \beta)$ such that $\mathrm{e}_{\mathscr{A}\otimes\beta^{\vee}}([i^*\psi_\beta \circ i^*](t)) = 1$. Then, take $\sigma$ to be the section

$$(s+t) \otimes [i^*\psi_\beta \circ i^*]\,(s+t) - t \otimes [i^*\psi_\beta \circ i^*](t) \quad \in \quad (H^0(\mathscr{A} \otimes \beta) \otimes H^0(\mathscr{A} \otimes \beta^{\vee}))^+.$$

Applying $2_X^* \otimes \mathrm{e}_{\mathscr{A}\otimes\beta^{\vee}}$ to $\sigma$ we get

$$\begin{aligned} 2_X^* \otimes \mathrm{e}_{\mathscr{A}\otimes\beta^{\vee}}(\sigma) &= 2_X^* \otimes \mathrm{e}_{\mathscr{A}\otimes\beta^{\vee}}\left((s+t) \otimes [i^*\psi_\beta \circ i^*]\,(s+t) - t \otimes [i^*\psi_\beta \circ i^*](t)\right) = \\ &= 2_X^*(s+t) \cdot 1 - 2_X^*(t) \cdot 1 = \\ &= 2_X^*(s). \end{aligned}$$

Therefore the statement is true. □

Using the above Theorem and a well known facts about multiplication maps on abelian varieties one is able to prove the following:

**Corollary 2.2.5.** *1. For every $\mathscr{A}$ ample symmetric invertible sheaf on $X$ the multiplication map*

$$m^+ : H^0(X, \mathscr{A}^{\otimes 2})^+ \otimes H^0(X, \mathscr{A}^{\otimes 2} \otimes \alpha) \longrightarrow H^0(X, \mathscr{A}^{\otimes 4} \otimes \alpha) \tag{2.2.3}$$

*is surjective for the generic $[\alpha] \in \mathrm{Pic}^0(X)$.*

*2. If furthermore $\mathscr{A}$ does not have a base divisor, then the locus where* (2.2.3) *is not surjective has codimension at least 2.*

*Proof.* The statement follows directly from the Theorem 2.2.4 and the corresponding statements (cfr. Proposition 2.1.12) about the map

$$m : H^0(X, \mathscr{A}^{\otimes 2}) \otimes H^0(X, \mathscr{A}^{\otimes 2} \otimes \alpha) \longrightarrow H^0(X, \mathscr{A}^{\otimes 4} \otimes \alpha).$$

□

Now we are ready to challenge the result that is the main point of this paragraph.

**Theorem 2.2.6.** *Let $\mathscr{A}$ be an ample symmetric invertible sheaf on $X$, then*

*1. there exist a non-empty open subset $U \subseteq \mathrm{Pic}^0(X)$ such that for every $h,\ n \in \mathbb{Z}$ with $n \geq 1$ and and every $\alpha \in U$ the following map is surjective*

$$m_\alpha^+ : H^0(X, \mathscr{A}^{\otimes 2n})^+ \otimes H^0(X, \mathscr{A}^{\otimes 2} \otimes \alpha) \longrightarrow H^0(X, \mathscr{A}^{\otimes 2n+2}. \otimes \alpha). \tag{2.2.4}$$

*2. If furthermore $\mathscr{A}$ does not have a base divisor, then the locus $Z$ in $\mathrm{Pic}^0(X)$ where* (2.2.4) *fails to be surjective has codimension at least 2.*

*Proof.* We will apply pretty much the same techniques used by Khaled in [45] in order to prove Proposition 2.2.1. We will procede by induction on $n$, with base given by Corollary 2.2.5.

*Case $n > 1$.* Observe the following commutative diagram:

$$\begin{array}{ccc} H^0(X, \mathscr{A}^{\otimes 2})^+ \otimes H^0(X, \mathscr{A}^{\otimes 2(n-1)})^+ \otimes H^0(X, \mathscr{A}^{\otimes 2} \otimes \alpha) & \xrightarrow{\varphi_\alpha} & H^0(X, \mathscr{A}^{\otimes 2})^+ \otimes H^0(X, \mathscr{A}^{\otimes 2(n-1)+2} \otimes \alpha) \\ \downarrow & & \downarrow \psi_\alpha \\ H^0(X, \mathscr{A}^{\otimes 2n})^+ \otimes H^0(X, \mathscr{A}^{\otimes 2} \otimes \alpha) & \xrightarrow{m_\alpha} & H^0(X, \mathscr{A}^{\otimes 2n+2} \otimes \alpha) \end{array}$$

The locus of points in $\mathrm{Pic}^0(X)$ where $m_\alpha$ is not surjective is contained in the following union

$$\{[\alpha] \in \mathrm{Pic}^0(X) \mid \varphi_\alpha \text{ is not surjective}\} \cup \{[\alpha] \in \mathrm{Pic}^0(X) \mid \psi_\alpha \text{ is not surjective}\}.$$

Hence we have

$$\operatorname{codim}\{[\alpha] \mid m_\alpha \text{ is not surjective}\} \geq$$
$$\min\{\operatorname{codim}\{[\alpha] \mid \varphi_\alpha \text{ is not surjective}\}, \operatorname{codim}\{[\alpha] \mid \psi_\alpha \text{ is not surjective}\}\}$$

Since $2 + 2(n-1) \geq 3$, by Proposition 2.2.1, $\psi_\alpha$ is surjective for every $[\alpha] \in \widehat{X}$ and the latter of the two sets above is empty. By inductive hypothesis,

$$\operatorname{codim}\{[\alpha] \in \operatorname{Pic}^0(X) \mid \varphi_\alpha \text{ is not surjective}\}$$

is greater than or equal to one, in the general case, to two, when $\mathscr{A}$ has not a base divisor, therefore the Theorem holds.

□

## 2.3 Equations and Syzygies of Kummer Varieties

Putting together the results of the previous paragraphs, in this last Section we will prove Theorem A and Theorem C.

First of all, observe that the case $p = 0$ of the Theorem A follows directly as a corollary of Khaled's work (cfr. Proposition 2.2.1). Thus, for what it follows we will always consider $p \geq 1$.

Our strategy in proving both Thoerem A and C will be using the part (b) of Lemma 2.1.3 and reduce the problem to checking the surjectvity of

$$H^0(\mathcal{K}_X,\, A^{\otimes n}) \otimes H^0(\mathcal{K}_X,\, M_{A^{\otimes n}}^{\otimes p} \otimes A^{\otimes nh}) \longrightarrow H^0(\mathcal{K}_X,\, M_{A^{\otimes n}}^{\otimes p} \otimes A^{\otimes n(h+1)}) \tag{2.3.1}$$

for every $h \geq r+1$. Denoting with $\mathscr{A}$ an ample symmetric line bundle on $X$ such that $\mathscr{A}^{\otimes 2} \simeq \pi^* A$, we split the proof in several steps.

**Step 1:** *Reduction to the surjectivity of*

$$H^0(X,\, \mathscr{A}^{\otimes 2n})^+ \otimes H^0(X,\, \pi^*(M_{A^{\otimes n}})^{\otimes p} \otimes \mathscr{A}^{\otimes 2nh}) \longrightarrow H^0(X,\, \pi^*(M_{A^{\otimes n}})^{\otimes p} \otimes \mathscr{A}^{\otimes 2n(h+1)}) \tag{2.3.2}$$

Our first aim will be to show that if the map (2.3.2) is surjective, then (2.3.1) is too. We prove this in the next lemma:

**Lemma 2.3.1.** *Let $A \in \operatorname{Pic}(\mathcal{K}_X)$ and $E$ be a vector bundle on $\mathcal{K}_X$. Denote by $\mathscr{E}$ the pullback $\pi_X^* E$. If the following multiplication map is surjective,*

$$H^0(X, \mathscr{A}^{\otimes 2n})^+ \otimes H^0(X, \mathscr{E}) \longrightarrow H^0(X, \mathscr{A}^{\otimes 2n} \otimes \mathscr{E}) \tag{2.3.3}$$

*then also*

$$H^0(\mathcal{K}_X, A^{\otimes n}) \otimes H^0(\mathcal{K}_X, E) \longrightarrow H^0(\mathcal{K}_X, A^{\otimes n} \otimes E) \tag{2.3.4}$$

*is surjective.*

*Proof.* The proof is very straightforward. The idea is that the vector space $V := H^0(X, \mathscr{E})$ splits as a direct sum $V^+ \oplus V^-$, where $V^+$ (respectively $V^-$) is the subspace of invariant (respectively anti-invariant) sections with respect to the $\mathbb{Z}/2\mathbb{Z}$ action defined by

$$s \mapsto s \circ i.$$

Therefore, if we denote by $f$ the map in (2.3.3), also $f$ splits as a direct sum $f^+ \oplus f^-$ with

$$f^{\pm} : H^0(X, \mathscr{A}^{\otimes 2n})^+ \otimes H^0(X, \mathscr{E})^{\pm} \longrightarrow H^0(X, \mathscr{A}^{\otimes 2n} \otimes \mathscr{E}).$$

If $f$ is surjective, then $f^+$ is surjective. The statement follows from the trivial identification of $f^+$ with (2.3.4). □

The next step in the proof will find us with reducing our problem to the surjectivity of a map of the type (2.1.6) and hence to an $M$-regularity problem. Before going any further we need some remarks.

*Remark* 2.3.1. Suppose that $A$ is an ample line bundle on $\mathcal{K}_X$ and let $\mathscr{A}$ an invertible sheaf on $X$ such that $\pi_X^* A = \mathscr{A}^{\otimes 2}$. Take $n$ an integer such that $A^{\otimes n}$ is globally generated and consider the following exact sequence of vector bundles:

$$0 \to M_{A^{\otimes n}} \longrightarrow H^0(\mathcal{K}_X,\, A^{\otimes n}) \otimes \mathscr{O}_{\mathcal{K}_X} \longrightarrow A^{\otimes n} \to 0.$$

By pulling back via the canonical surjection $\pi_X$ we get:

$$0 \to \pi_X^*(M_{A^{\otimes n}}) \longrightarrow H^0(X,\, \mathscr{A}^{\otimes 2n})^+ \otimes \mathscr{O}_X \longrightarrow \mathscr{A}^{\otimes 2n} \to 0.$$

Hence, after defining

$$W_n := H^0(X,\, \mathscr{A}^{2n})^+$$

we have that $\pi_X^* M_{A^{\otimes n}} \simeq M_{W_n}$.

From now on given a sheaf $\mathscr{F}$ on $X$, we will often write $H^0(\mathscr{F})$ instead of $H^0(X,\, \mathscr{F})$.

**Step 2:** *Reduction to an $M$-regularity problem*

Consider the vanishing locus $V^1(M_{W_n} \otimes \mathscr{A}^{\otimes 2})$, we claim that it coincides with the locus of $\alpha \in \widehat{X}$ such that the multiplication map

$$m_\alpha^+ : H^0(\mathscr{A}^{\otimes 2n})^+ \otimes H^0(\mathscr{A}^{\otimes 2} \otimes \alpha) \longrightarrow H^0(\mathscr{A}^{\otimes 2(n+1)} \otimes \alpha)$$

is not surjective. Infact consider the short exact sequence

$$0 \to M_{W_n} \longrightarrow W_n \otimes \mathscr{A}^{\otimes 2} \otimes \alpha \longrightarrow \mathscr{A}^{\otimes 2(n+1)} \otimes \alpha \to 0.$$

Taking cohomology we get

$$H^0(M_{W_n} \otimes \alpha) \to H^0(\mathscr{A}^{\otimes 2n})^+ \otimes H^0(\mathscr{A}^{\otimes 2} \otimes \alpha) \xrightarrow{m_\alpha^+} H^0(\mathscr{A}^{\otimes 2(n+1)} \otimes \alpha) \to \cdots$$
$$\cdots \to H^1(M_{W_n} \otimes \alpha) \to H^0(\mathscr{A}^{\otimes 2n})^+ \otimes H^1(\mathscr{A}^{\otimes 2} \otimes \alpha) \to \cdots$$

Since for every topologically trivial line bundle $\alpha \in \mathrm{Pic}^0(X)$ $H^0(\mathscr{A}^{\otimes 2n})^+ \otimes H^1(\mathscr{A}^{\otimes 2} \otimes \alpha) = 0$, it follows that the surjectivity of $m_\alpha^+$ is equivalent to the vanishing of $H^1(M_{W_n} \otimes \alpha)$.

Thank of this characterization of the locus on $\mathrm{Pic}^0(X)$ where $m_\alpha^+$ fails to be surjective and using a technique introduced by Kempf and widely employed we were able to prove the following

**Lemma 2.3.2.** *Let $\mathscr{A}$ and $\mathscr{E}$ be an ample symmetric sheaf on an abelian variety $X$ and a coherent sheaf on $X$, respectively. If $\mathscr{E} \otimes \mathscr{A}^{\otimes -2}$ is $M$-regular, then the multiplication map*

$$H^0(X, \mathscr{A}^{\otimes 2n})^+ \otimes H^0(X, \mathscr{E}) \longrightarrow H^0(X, \mathscr{A}^{\otimes 2n} \otimes \mathscr{E}) \tag{2.3.5}$$

*is surjective for every $n \geq 1$.*

Before proceeding with the proof we will state an immediate corollary of this Lemma, that reduce our problem to an $M$-regularity problem

**Corollary 2.3.3.** *If $M_{W_n}^{\otimes p} \otimes \mathscr{A}^{\otimes 2(nh-1)}$ is $M$-regular, then* (2.3.2) *is surjective.*

*Proof of the Lemma.* By Proposition 2.2.6(1) we know that $V^1(M_{W_n} \otimes \mathscr{A}^{\otimes 2})$ is contained in a Zariski closed subset of $\widehat{X}$. Therefore it exists an open set $\widehat{U}_0 \subseteq \mathrm{Pic}^0(X)$ such that $\widehat{U}_0 \cap V^1(M_{W_n} \otimes \mathscr{A}^{\otimes 2}) = \emptyset$. Now observe the following commutative diagram.

$$\begin{array}{ccc}
\bigoplus_{\alpha \in \widehat{U}_0} H^0(\mathscr{A}^{\otimes 2n})^+ \otimes H^0(\mathscr{A}^{\otimes 2} \otimes \alpha) \otimes H^0(\mathscr{E} \otimes \mathscr{A}^{\otimes -2} \otimes \alpha^\vee) & \longrightarrow & H^0(\mathscr{A}^{\otimes 2n})^+ \otimes H^0(\mathscr{E}) \\
\downarrow f & & \downarrow g \\
\bigoplus_{\alpha \in \widehat{U}_0} H^0(\mathscr{A}^{\otimes 2n+2} \otimes \alpha) \otimes H^0(\mathscr{E} \otimes \mathscr{A}^{\otimes -2} \otimes \alpha^\vee) & \xrightarrow{h} & H^0(\mathscr{E} \otimes \mathscr{A}^{\otimes 2n})
\end{array}$$

The map $f = \oplus m_\alpha^+$ is surjective by our choice of the set $\widehat{U}_0$, the map $h$ is surjective by $M$-regularity hypothesis together with Theorem 2.1.10. Thus $g$ is necessarily surjective.

□

**Step 3:** *Solution of the $M$-regularity problem.*

Theorem 2.1.10 and the result in the last paragraph allow us to reduce the problem of the surjectivity of the map (2.3.5) to an $M$-regularity problem. In particular we have that, if $M_{W_n}^{\otimes p} \otimes \mathscr{A}^{\otimes 2(nh-1)}$ is $M$-regular for every $h \geq r$, then $A^{\otimes n}$ satisfies $N_p^r$. The solution of the $M$-regularity problem is presented in the next two statment.

**Proposition 2.3.4.** *Let $p$ a positive integer. Then $M_{W_n}^{\otimes p} \otimes \mathscr{A}^{\otimes m}$ satisfies I.T. with index 0 (and hence it is $M$-regular) for every $m \geq 2p+1$*

Note that Theorem A follows at once from this Proposition taking $m = 2nr - 2$.

*Proof.* We will procede by induction on $p$.

*Case* $p = 1$. Let us consider the following exact sequence:

$$0 \to M_{W_n} \longrightarrow W_n \otimes \mathscr{O}_X \longrightarrow \mathscr{A}^{\otimes 2n} \to 0.$$

twistin with $\mathscr{A}^{\otimes m} \otimes \alpha$ with $[\alpha]$ any element $\widehat{X}$ we obtain

$$0 \to M_{W_n} \otimes \mathscr{A}^{\otimes m} \otimes \alpha \longrightarrow W_n \otimes \mathscr{A}^{\otimes m} \otimes \alpha \longrightarrow \mathscr{A}^{\otimes 2n+m} \otimes \alpha \to 0. \tag{2.3.6}$$

Hence one can easily see that the vanishing of the higher cohomology of $M_{W_n} \otimes \mathscr{A}^{\otimes m} \otimes \alpha$ depends upon:

(i) the vanishing of the higher cohomology of $\mathscr{A}^{\otimes m} \otimes \alpha$ and

(ii) the surjectivity of the following multiplication map:

$$H^0(X,\, \mathscr{A}^{\otimes 2n})^+ \otimes H^0(X,\, \mathscr{A}^{\otimes m} \otimes \alpha) \longrightarrow H^0(X,\, \mathscr{A}^{\otimes 2n+m} \otimes \alpha).$$

Condition (i) holds for every $\alpha$ as long as $m \geq 1$, while, thank to Khaled result (cfr Proposition 2.2.1), we know condition (ii) holds for every $\alpha$ as long as $m \geq 3$.

*Case* $p > 1$. Suppose now that $p > 1$ and take any $[\alpha] \in \widehat{X}$. By twisting (2.3.6) by $M_{W_n}^{\otimes p-1}$ we can observe that the vanishing of higher cohomology of $M_{W_n}^{\otimes p} \otimes \mathscr{A}^{\otimes m} \otimes \alpha$ is implied by

(i) the vanishing of the higher cohomology of $M_{W_n}^{\otimes p-1} \otimes \mathscr{A}^{\otimes m} \otimes \alpha$ and

(ii) the surjectivity of the following multiplication map:

$$H^0(X, \mathscr{A}^{\otimes 2n})^+ \otimes H^0(X, M_{W_n}^{\otimes p-1} \otimes \mathscr{A}^{\otimes m} \otimes \alpha) \longrightarrow H^0(X, M_{W_n}^{\otimes p-1} \otimes \mathscr{A}^{\otimes 2n+m} \otimes \alpha) \tag{2.3.7}$$

By induction (i) holds as long as $m \geq 2p-1$. Thank to Lemma 2.3.3 and Lemma 2.1.10 we know that if $M_{W_n}^{\otimes p-1} \otimes \mathscr{A}^{\otimes m-2} \otimes \alpha$ satisfies I.T. with index 0, than (2.3.7) is surjective holds. But we use induction again and we get that $M_{W_n}^{\otimes p-1} \otimes \mathscr{A}^{\otimes m-2} \otimes \alpha$ is I.T. with index 0 whenever $m-2 \geq 2p-1$, that is whenever $m \geq 2p+1$ and hence the statement is proved. □

**Proposition 2.3.5.** *In the notations above, take $p \geq 1$ an integer. If $\mathscr{A}$ does not have a base divisor, then $M_{W_n}^{\otimes p} \otimes \mathscr{A}^{\otimes m}$ is $M$-regular for every $m \geq 2p$.*

Again Theorem C, follows at once after taking $m = 2(nr-1)$.

*Proof.* For $m \geq 2p+1$ the statement is a direct consequence of the Proposition above, hence we can limit ourselves to the case $m = 2p$. We will procede by induction on $p$.

*Case $p = 1$*: We want to prove that $\operatorname{codim} V^i(M_{W_n} \otimes \mathscr{A}^{\otimes 2}) > i$ for every $i \geq 1$. From the vanishing of the higher cohomology of $\mathscr{A}^{\otimes 2} \otimes \alpha$ for every $\alpha \in \operatorname{Pic}^0(X)$ we know that the loci

$$V^i(M_{W_n} \otimes \mathscr{A}^{\otimes 2}) = \emptyset \quad \text{for every } i \geq 2.$$

Recall that the locus $V^1(M_{W_n} \otimes \mathscr{A}^{\otimes 2})$ is the locus of points $\alpha \in \widehat{X}$ such that the multiplication

$$W_n \otimes H^0(\mathscr{A}^{\otimes 2} \otimes \alpha) \longrightarrow H^0(A^{\otimes 2n+2} \otimes \alpha)$$

is not surjective. We know from Proposition 2.2.6 that if $\mathscr{A}$ has not a base divisor then this locus has at least codimension 2 and hence the statement is proved.

*Case $p > 1$* Take $[\alpha] \in \widehat{X}$. Consider the following exact sequence

$$0 \to M_{W_n}^{\otimes p} \otimes \mathscr{A}^{\otimes 2p} \otimes \alpha \longrightarrow W_n \otimes M_{W_n}^{\otimes p-1} \otimes \mathscr{A}^{\otimes 2p} \otimes \alpha \longrightarrow M_{W_n}^{\otimes p-1} \otimes \mathscr{A}^{\otimes 2p+n} \otimes \alpha \to 0 \tag{2.3.8}$$

From Proposition 2.3.4(a) we know that for every $i \geq 1$ both $H^i(M_{W_n}^{\otimes p-1} \otimes \mathscr{A}^{\otimes 2p} \otimes \alpha)$ and $H^i(M_{W_n}^{\otimes p-1} \otimes \mathscr{A}^{\otimes 2p+n} \otimes \alpha)$ vanish. Thus the loci $V^i(M_{W_n}^{\otimes p} \otimes \mathscr{A}^{\otimes 2p})$ are empty for every $i \geq 2$. It remains to show that that

$$\operatorname{codim} V^1(M_{W_n}^{\otimes p} \otimes \mathscr{A}^{\otimes 2p}) \geq 2.$$

As before one may observe that this locus is exactly the locus in $\widehat{X}$ where the following multiplication map fails to be surjective:

$$W_n \otimes H^0(M_{W_n}^{\otimes p-1} \otimes \mathscr{A}^{\otimes 2p} \otimes \alpha) \longrightarrow H^0(M_{W_n}^{\otimes p-1} \otimes \mathscr{A}^{\otimes 2p+2n} \otimes \alpha).$$

Infact, taking cohomology in (2.3.8) and observing that for every $[\alpha] \in \widehat{X}$, $h^1(M_{W_n}^{\otimes p-1} \otimes \mathscr{A}^{\otimes 2p} \otimes \alpha) = 0$, due to Proposition 2.3.4, we have the conlcusion following the same argument that in the case $p = 0$.

Now take $[\alpha] \in V^1(M_{W_n}^{\otimes p} \otimes \mathscr{A}^{\otimes 2p})$. By inductive hypothesis the sheaf $M_{W_n}^{\otimes p-1} \otimes \mathscr{A}^{\otimes 2(p-1)}$ is $M$-regular. Corollary 2.1.11 implies that there exists a positive integer $N$ and $[\beta_1], \dots, [\beta_N] \in \widehat{X}$ such that the following is surjective.

$$\bigoplus_{k=1}^{N} H^0(\mathscr{A}^{\otimes 2n+2} \otimes \beta_k \otimes \alpha) \otimes H^0(M_{W_n}^{\otimes p-1} \otimes \mathscr{A}^{\otimes 2(p-1)} \otimes \beta_k^{\vee}) \xrightarrow{m_{\beta_k}} H^0(M_{W_n}^{\otimes p-1} \otimes \mathscr{A}^{\otimes 2p+2n} \otimes \alpha)$$

Consider the commutative square

$$\begin{array}{ccc}
\bigoplus_{k=1}^{N} H^0(\mathscr{A}^{\otimes 2n})^+ \otimes H^0(\mathscr{A}^{\otimes 2} \otimes \alpha \otimes \beta_k) \otimes H^0(M_{W_n}^{\otimes p-1} \otimes \mathscr{A}^{\otimes 2p-2} \otimes \beta_k^{\vee}) & \longrightarrow & H^0(\mathscr{A}^{\otimes 2n})^+ \otimes H^0(M_{W_n}^{\otimes p-1} \otimes \mathscr{A}^{\otimes 2p} \otimes \alpha) \\
\downarrow & & \downarrow \\
\bigoplus_{k=1}^{N} H^0(\mathscr{A}^{\otimes 2n+2} \otimes \alpha \otimes \beta_k) \otimes H^0(M_{W_n}^{\otimes p-1} \otimes \mathscr{A}^{\otimes 2(p-1)} \otimes \beta_k^{\vee}) & \longrightarrow & H^0(M_{W_n}^{\otimes p-1} \otimes \mathscr{A}^{\otimes 2n+2p})
\end{array}$$

The right arrow is not surjective by our choice of $\alpha$. The bottom arrow is surjective, hence the left arrow could not be surjective. Therefore

$$\alpha \in \bigcup_{k=1}^{N} Z_k,$$

where $Z_k$ stands for the locus of $[\beta] \in \widehat{X}$ such that the multiplication map

$$H^0(\mathscr{A}^{\otimes 2n})^+ \otimes H^0(\mathscr{A}^{\otimes 2} \otimes \beta \otimes \beta_k) \to H^0(\mathscr{A}^{\otimes 2n+2} \otimes \beta \otimes \beta_k) \tag{2.3.9}$$

fails to be surjective. Thus one has that

$$V^1(M_{W_n}^{\otimes p} \otimes \mathscr{A}^{\otimes 2p}) \subseteq \bigcup_{k=1}^{N} Z_k.$$

By Theorem 2.2.6(2) the loci $Z_k$ have codimension at least 2, therefore

$$\operatorname{codim} V^1(M_{W_n}^{\otimes p} \otimes \mathscr{A}^{\otimes 2p}) \geq \operatorname{codim} \bigcup_{k=1}^{N} Z_k \geq 2$$

and the statement is proved. □

## 2.4 Final Remarks

In this section we present some possible further developments of the achievements presented in this Chapter.

First of all we remark that we do not expect our results to be optimal for $p \geq 2$. In fact we guess, following the analogy with abelian varieties, that for large $p$, $A^{\otimes n}$ should satisfies property $N_p$ for $n \sim \frac{p}{2}$ rather than for $n \sim p$. However any result of this kind is impossible to achieve with the methods employed here. Our hope is that using Ploog's *equivariant Fourier-Mukai transform* ([67]) instead of classical integral transforms is possible to develope an $M$-regularity theory on Kummer varieties and to use it to investigates their syzygies.

Another possible development is to find an explicit basis for the quartics that cut a Kummer variety $\mathcal{K}_X$ associated to a pricipally polarized abelian variety $X$ and embedded by a divisor $2\Theta$ with $\Theta$ a principal polarization in $X$. In fact Khaled in [46] found such a base when $2\Theta$ yields a projectively normal embedding (indeed he was able to prove that Kummer varieties are defined by quartics equation just in this setting). Unfortunately the open sets of the moduli space of Kummer varieties of dimension $g$ constitued by projectively normal varieties is not much interesting from a geometric point of view. For example, for what it concerns Jacobians $J(C)$, the projective normality of the embedding $2\Theta$ is equivalent to the existence of an even theta-characteristic on $C$ (cfr. [46]); hence there are examples of Jacobians with non-projectively normal associated Kummer, as well ones whose associated Kummer variety satisfies $N_0$. Therefore, finding the equation defining Kummer varieties in a non projectively normal setting could have many applications and could lead to a better comprehension of the geometry of the moduli spaces of Kummer varieties and consequently to a better understanding of the Shottky problem.

### Changing Group

One may wonder if the techniques described in this Chapter are adapt to study the syzygies of other (singular) quotients of abelian varieties. More precisely let $G$ be a finite group of authomorphisms of an abelian variety $X$, we can consider the (G.I.T.) quotient $\pi_X^G : X \to \mathcal{K}_X^G := X/G$ and investigate the syzygies of an embedding of his.

The first problem that we encounter in this study is that we do not have the nice characterization of ample line bundles on $\mathcal{K}_X^G$ in terms of line bundles on $X$ we had for classical Kummer varieties. Namely we do not know if the pullback of any ample line bundle on $\mathcal{K}_X^G$ can be described as a (fixed) power of an ample line bundle on $X$. What we can say is that, for every $L$ $G$-invariant ample invertible sheaf on $X$ (e.g. a sheaf such that for every $g \in G$ there is an isomorphism $\psi_g : L \longrightarrow g^*L$), the sheaf $L^{\otimes|G|}$ is of the form $\pi_X^{G,\,*}\mathscr{L}$ with $\mathscr{L} \in \mathrm{Pic}(\mathcal{K}_X^G)$. Hence the first step toward the study of syzygies of these quotient varieties would be restricting onself to the analysis of the embeddings given by those line bundles $\mathscr{L}$ such that $p_X^{G,*}\mathscr{L} \simeq L^{\otimes|G|}$ for some $L \in \mathrm{Pic}(X)$.

Next we need to check if the proof in the previous Sections are "stable by changing group". After substituting the groups $H^i(X, L^{\otimes n})^+$ with the equivariant cohomology groups $H^i(X, L^{\otimes n})^G$, there is no problem whatsoever for what it concerns the results of Section 2.3. The issues arise when we are looking at the "base" of our induction, namely at section 2.2. in fact the results we find there are proved using techniques tailor made for the $\mathbb{Z}/2\mathbb{Z}$ actions. Thus, in order to study the syzygies of generalized Kummer varieties it is necessary to give answer to the following problem:

**Problem.** *Let $G$ be a finite group acting on an abelian variety $X$. Take $L$ a $G$-invariant ample line bundle on $X$. Consider $V_G(n, k)$ the locus in $\widehat{X}$ where the following is NOT surjective:*

$$m_\alpha^G : H^0(X, L^{\otimes n})^G \otimes H^0(X, L^{\otimes k} \otimes \alpha) \longrightarrow H^0(X, L^{\otimes n+k} \otimes \alpha).$$

*For which $k$ and $n$ nonnegative integers*

1. *is $V_G(n, k)$ empty?*
2. *has $V_G(n, k)$ codimension greater than 1?*
3. *has $V_G(n, k)$ positive codimension?*

A possible way to circumvent the problem could be, again, adopting the equivariant Fourier-Mukai transform and see if it can be used to develope an equivariant $M$-regularity theory.

# Part II

# Pluricanonical Maps of Varieties of Maximal Albanese Dimension

# CHAPTER 3

# TETRACANONICAL MAPS OF VARIETIES OF MAXIMAL ALBANESE DIMENSION

An interesting issue in birational geometry is studying the structure of the pluricanonical maps of smooth varieties. In particular, given $Z$ a complex smooth projective variety, one wants to find explicitly (when it exists) an integer $n_0$ such that for every $n \geq n_0$ the pluricanonical linear system $|nK_Z| \simeq |\omega_Z^{\otimes n}|$ yields a map birational equivalent to the Iitaka fibration of $X$.

For example, in the case of curves and surfaces of general type, the answer to this question had been long known: if $C$ is a curve of genus $g \geq 2$, then an easy application of Riemann-Roch Theorem tells us that the tricanonical map $\varphi_{|\omega_Z^{\otimes 3}|}$ is birational. The case of surfaces was succesfully challenged by Bombieri in [6], who proved that given $S$ is a surface of general type, then the pentacanonical, $\varphi_{|\omega_Z^{\otimes 5}|}$, is always birational.

Under the further assumption that $Z$ is of *maximal Albanese dimension* (i.e. the Albanese map $Z \longrightarrow \mathrm{Alb}(Z)$ is generically finite), the pluricanonical maps are surprisling easier to understand. In [10] Chen and Hacon proved that if $Z$ is a smoooth complex projective variety of maximal Albanese dimension, then the image of $\varphi_{|\omega_Z^{\otimes n}|}$ has dimension $\kappa(X)$ for any $n \geq 6$; if furthermore $Z$ is of general type then this map its birational onto its image for any $n \geq 6$. Moreover if the Albanese image of $Z$ is not ruled by tori) than the complete linear system $|mK_Z|$ induces a birational map for every $m \geq 3$. These statements were later recoverd by Pareschi–Popa ([60]) as an application of their Fourier-Mukai based techniques. More recently it was shown by Jiang [36], applying ideas from [60], that, if $Z$ is a smoooth projective variety of maximal Albanese dimension, then $\varphi_{|\omega_Z^{\otimes 5}|}$ is a model of the Iitaka fibration. The main result of this Chapter, whose proof uses $M$-regularity techniques introduced by Pareschi and Popa in [60], is an improvement of Jiang's theorem:

**Theorem 3.A.** *If $Z$ is a complex projective smooth variety of maximal Albanese dimension and general type, then the tetracanonical map $\varphi_{|\omega_Z^{\otimes 4}|}$ is birational onto its image.*

The argument consists in showing that reducible divisors $D_\alpha + D_{\alpha^\vee}$, with $[\alpha] \in \mathrm{Pic}^0(Z)$, $D_\alpha \in |\omega_Z^{\otimes 2} \otimes \alpha|$ and $D_{\alpha^\vee} \in |\omega_Z^{\otimes 2} \otimes \alpha^\vee|$ separate points in a suitable open set of $Z$. The crucial point is that, for all $[\alpha] \in \mathrm{Pic}^0(Z)$, the sections of $\omega_Z^{\otimes 2} \otimes \alpha$ passing trough a general point of $Z$ have good generation properties. In order to achieve this we use Pareschi-Popa theory of $M$-regularity and continuous global generation ([61, 60]), joint with a theorem of Chen-Hacon [10] on the fact that the variety $V^0(\omega_Z)$ spans $\mathrm{Pic}^0(Z)$.

It is worth mentioning that this achievement is not sharp. In fact in the next Chapter we will present a proof, obtained in collaboration with Z. Jiang and M. Lahoz, of the birationality of the tricanonical map for varieties of general type and maximal Albanese dimension. However tha aforementioned proof, since it relays on induction on dimension, is less explicit than the one exposed here.

In what follows $Z$ will always be a smooth complex variety of general type and maximal Albanese dimension while $\omega_Z$ shall denote its dualizing sheaf. By $\mathrm{Alb}(Z)$ we will mean the Albanese variety of $Z$.

## 3.1 Background Material

### 3.1.1 Asymptotic Multiplier Ideals and Related Vanishing Properties

In this paragraph we briefly recall the basic properties of multiplier ideals and asymptotic multiplier ideals. For a complete treatment of this matter we recommend Chapters 9 to 11 of [50].

Assume that $Y$ is a smooth variety of dimension $n$ and let $D$ be an integral divisor on $Y$. Let $V \subseteq H^0(Y, \mathscr{O}_Y(D))$ be a non-zero finite dimensional linear subspece. We recall that a *log resolution of the linear series* $|V|$ is a projective birational mapping

$$\mu : Y' \to Y$$

with $Y'$ non-singular and such that the $\mu^*|V| = |W| + F$ where $F + \mathrm{exc}(\mu)$ is a divisor with simple normal crossing support and $W \subseteq H^=(Y', \mathscr{O}_{Y'}(\mu^* D - F))$ is a base point free linear series. (Here $\mathrm{exc}(\mu)$ stands for the sum of exceptional divisors of $\mu$).

Now, following [50, Definition 9.2.10 ], given $D$, $V$, $\mu$ as above and $c > 0$ a rational number one can define the *multiplier ideal* $\mathscr{J}(c \cdot |V|)$ corresponding to $c$ and $|V|$ as

$$\mathscr{J}(c \cdot |V|) := \mu_* \mathscr{O}_{Y'}(K_{Y'/Y} - \lfloor c \cdot F \rfloor).$$

It can be checked that this object does not depend on the choice of the log resolution $\mu$. If we assume that the linear series $|D|$ has non-negative Iitaka dimension and we take $p$ a positive integer, we can form the multiplier ideal sheaf associated to the complete linear series $|pD|$

$$\mathscr{J}(\frac{c}{p}|pD|) \subseteq \mathscr{O}_Y.$$

We have the following Lemma

**Lemma 3.1.1.** *For every integer $k \geq 1$ there is an inclusion*

$$\mathscr{J}(\frac{c}{p}|pD|) \subseteq \mathscr{J}(\frac{c}{pk}|pkD|)$$

This result together with the ascending chain condition on ideals tells us that the family

$$\{\mathscr{J}(\frac{c}{p}|pD|)\}_{p\geq 0}$$

admits a unique maximal element $\mathscr{J}(c||D||)$, called the *asymptotic multiplier ideal associated to $c$ and $|D|$.* The main reason for the introduction of multiplier ideals is the possibility to extend classical results, such as Kodaira vanishing for big and nef line bundles. In fact we have the following important achievement that alone is the cause of our need for multiplier ideals.

**Theorem 3.1.2** (Nadel vanishing for asymptotic multiplier ideals [50, Thm. 11.2.12])**.** *In the notation above, if $D$ has non-negative Iitaka dimension, we have the following vanishing:*

*(i) If $A$ is big and nef integral divisor on $Y$, then*

$$H^i(Y, \mathscr{O}_Y(K_Y + mD + A) \otimes \mathscr{J}(||mD||)) = 0 \quad \textit{for } i > 0$$

*(ii) If $D$ is big then the same statement is true assuming just the nefness of $A$. In particular if $D$ is big then*

$$H^i(Y, \mathscr{O}_Y(K_Y + mD) \otimes \mathscr{J}(||mD||)) = 0 \quad \textit{for } i > 0$$

Since multiplier ideals are, in particular, ideal sheaves, it is very natural to ask ourselves questions about the geometry of their zero locus. We have the following results:

**Lemma 3.1.3.** *Given a $L$ line bundle on $Y$ of non-negative Iitaka dimension*

*(i) Then*

$$H^0(Y, L \otimes \mathscr{J}(||L||)) = H^0(Y, L),$$

*i.e. the zero locus of $\mathscr{J}(||L||)$ is contained in the base locus of $L$.*

*(ii) For every non negative integer k,*

$$\mathscr{J}(||L^{\otimes(k+1)}||) \subseteq \mathscr{J}(||L^{\otimes k}||)$$

*Proof.* Both results are exposed in [50], part (i) is Proposition 11.2.10, while part (ii) is Theorem 11.1.8. □

As an easy corollary we have the following Lemma due to Pareschi and Popa that we will be needing afterward.

**Lemma 3.1.4** ([60, Lemma 6.3(b)]). *Let $Y$ be a smooth complex variety of general type and denote by $K_Y$ its canonical divisor. For any $m > 1$, the zero locus of $\mathscr{J}(||(m-1))K_Y||)$ is contained in the base locus of $\mathscr{O}_Y(mK_Y) \otimes \alpha$ for every $[\alpha] \in \mathrm{Pic}^0(Y)$.*

For the reader's benefit we include the short proof.

*Proof.* Since bigness is a numerical property, all line bundles $\omega_Y \otimes \alpha$ are big. By Nadel's vanishing (Theorem 3.1.2) for any $m > 1$ we have

$$h^i(Y, \omega_Y^{\otimes m} \otimes \beta \otimes \mathscr{J}(||(\omega_Y \otimes \alpha)^{\otimes(m-1)}||)) = 0$$

for every $i > 0$ and every $[\alpha]$, $[\beta] \in \mathrm{Pic}^0(Y)$. Since the Euler characteristic is invariant under smooth deformations, $h^0(Y, \omega_Y^{\otimes m} \otimes \beta \otimes \mathscr{J}(||(\omega_Y \otimes \alpha)^{\otimes(m-1)}||))$ does not depend on $[\beta]$. Denote this quantity by $\lambda(\alpha)$. Now, Lemma 3.1.3 implies that for any $[\beta] \in \mathrm{Pic}^0(Y)$

$$h^0(Y, \omega_Y^{\otimes m} \otimes \beta \otimes \mathscr{J}(||(\omega_Y \otimes \alpha)^{\otimes(m-1)}||)) \leq h^0(Y, \omega_Y^{\otimes m} \otimes \beta)$$

and equality holds for $[\beta] = [\alpha^{\otimes m}]$. As a consequence we have that, by semicontinuity, for any $[\beta]$ in a neighbourhood of $[\alpha^{\otimes m}]$, $h^0(Y, \omega_Y^{\otimes m} \otimes \beta) = \lambda(\alpha)$. Since this holds for every $[\alpha]$ it follows that $h^0(Y, \omega_Y^{\otimes m} \otimes \beta)$ is constant (and equal to $\lambda(\mathscr{O}_{\mathrm{Pic}^0(Y)})$) for any $[\beta]$ and $m > 1$. The statement follows directly. □

In what follows we will adopt the following notation. Given $\mathscr{L}$ an invertible sheaf on a smooth projective variety $Y$ such that $H^0(Y, \mathscr{L}) \neq 0$; we will denote by $\mathscr{J}(||\mathscr{L}||)$ the asymptotic multiplier ideal sheaf associated to the complete linear series $|D|$ with $D$ an effective integral divisor with $\mathscr{O}_Y(D) \simeq \mathscr{L}$.

### 3.1.2 Iitaka Fibration

We conclude this section of preliminaries by recalling the definition and the basic properties of the Iitaka fibration associated to a line bundle. We will discuss Iitaka fibrations associated to the canonical sheaf of a maximal Albanese dimension variety more throughly in Chapter 4, but since we spoke about it in the introduction to this Chapter, we deemed it best to expose this material in here. We follow [50, Section 2.1.C].

Let $Y$ be a smooth projective variety and consider $L$ a line bundle such that its Iitaka dimension $\kappa(Y, L)$ is non negative. Let

$$\mathbf{N}(Y, L) := \{k \in \mathbb{N} \mid h^0(Y, L^{\otimes k}) \neq 0\},$$

and denote by $\varphi_k : Y \dashrightarrow Z_k \subseteq \mathbb{P}$ the rational map associated to the complete linear system $|L^{\otimes k}|$. As $k$ grows these maps turn out to have a "limit" in the following sense:

**Theorem 3.1.5** ([50, Theorem 2.1.19]). *In the notation above if $k \in \mathbf{N}(Y, L)$ is big enough then the map $\varphi_k$ is birational equivalent to a given fiber space*

$$\varphi_\infty : Y_\infty \to Z_\infty \tag{3.1.1}$$

*of normal varieties and the restriction of $L$ to a general fiber of $\varphi_\infty$ has Iitaka dimension 0. In particular it exists a commutative diagram*

$$\begin{array}{ccc} Y_\infty & \xrightarrow{u_\infty} & Y \\ {\scriptstyle \varphi_\infty}\downarrow & & \downarrow{\scriptstyle \varphi_k} \\ Z_\infty & \underset{u_k}{\dashrightarrow} & Z_k \end{array}$$

*where the horizontal maps are birational and $u_\infty$ is a morphism. Moreover setting $L_\infty := u_\infty^* L$ and taking $F$ to be a very general fiber of $\varphi_\infty$ then*

$$\kappa(F, L_{\infty|F}) = 0.$$

Observe that $\dim Z_\infty \simeq \kappa(Y, L)$.

The algebraic fiber space (3.1.1) is called the *Iitaka fibration associated to $L$* and it is unique up to birational equivalence. When $L = \omega_Y$ the canonical line bundle on $Y$, then we speak about the *Iitaka fibration of $Y$*.

Since the issues we will be addressing have a birational nature, we will usually assume that $Y_\infty = Y$.

It has been proved by Chen and Hacon [11] that $Y$ is of maximal Albanese dimension if and only if $Z_\infty$ is. In this case the general fiber of $\varphi_\infty$ is again of maximal Albanese dimension. Since it has Kodaira dimension 0, a result of Kawamata [40] implies that it is birational to an abelian variety $K$.

### 3.1.3 Generation Properties of M-regular Sheaves

In the first Chapter we recalled the notion of $M$-regular sheaf on an abelian variety and we listed some of their main properties. However for what we want to prove here these notions are not sufficient: since we will need to work with the generation properties of $M$-regular sheaves. Below we recall some basic terminology about generation of sheaves and we present results of Pareschi and Popa that investigate the generation properties for Mukai-regular sheaves.

Given $\mathscr{F}$ a coherent sheaf on variety $Y$ over an algebraically, closed field $k$, we say that it is *globally generated* (in brief $GG$) if the evaluation map

$$H^0(Y, \mathscr{F}) \otimes \mathscr{O}_Y \to \mathscr{F}$$

is surjective. More generally, given $T \subseteq Y$ a proper subvariety, we say that $\mathscr{F}$ is *globally generated away from* $T$ if for every $y \notin T$ the map

$$H^0(Y, \mathscr{F}) \otimes k(y) \to \mathscr{F} \otimes k(y)$$

is surjective. The following notion has been introduced to further weaken the concept of global generation.

**Definition 3.1.6** ([60, Definition 2.10])**.** Let $Y$ be a variety

1. A sheaf $\mathscr{F}$ on $Y$ is *continuosly globally generated* (in brief $CGG$) if the sum of the evaluation maps
$$\mathcal{E}v_U : \bigoplus_{\alpha \in U} H^0(\mathscr{F} \otimes \alpha) \otimes \alpha^{-1} \longrightarrow \mathscr{F}$$
is surjective for every $U$ non empty open set of $\mathrm{Pic}^0(Y)$.

2. More generally, given $T$ a proper subvariety of $Y$, $\mathscr{F}$ is said to be *continuosly globally generated away from* $T$ if $\mathrm{Supp}(\mathrm{Coker}(\mathcal{E}v_U))$ is contained in $T$ for every $U$ non empty open set of $\mathrm{Pic}^0(Y)$.

The following proposition relates the behavior of the $CGG$ property with respect to the twisting by a line bundle.

**Proposition 3.1.7** ([61, Proposition 2.12])**.** *In the settings above, assume that $\mathscr{F}$ is CGG away from a subvariety $T$ and let $\mathscr{A}$ a line bundle that is CGG away from a subvariety $W$ of $Y$, then for every $\alpha \in \mathrm{Pic}^0(Y)$ the sheaf $\mathscr{F} \otimes \mathscr{A} \otimes \alpha$ is globally generated away from $T \cup W$.*

We conclude this paragraph by recalling a result of Pareschi and Popa that will enable us to study generation properties of certain sheaves on an abelian variety and derive from them the generation properties of $\omega_Z^{\otimes 4}$ on $Z$.

**Proposition 3.1.8** ([60])**.** *An $M$-regular sheaf on an abelian variety is continuosly globally generated.*

### 3.1.4 On the Geometry of Generic Vanishing Loci

In this subsection we expose some results of Chen–Hacon [10, Theorem 1], Green–Lazarsfeld [27, Theorem 0.1] and Simpson [69, Section 4,6,7] about the geometry of the generic vanishing loci introduced in the first chapter.

**Theorem 3.1.9** (Subtorus Theorem [29, 69])**.** *Let $Y$ be a compact Kähler manifold, and $W$ an irreducible component of $V^i(\omega_Y)$ for some $i$. Then*

*(a) There exist a torsion point $[\beta]$ and a subtorus $T$ of $\mathrm{Pic}^0(Y)$ such that $W = [\beta] + T$.*

*(b) There exist a normal analytic variety $X$ of dimension $\leq d-i$, such that (any smooth model of) $X$ has maximal Albanese dimension and a morphism with connected fibres $f : Y \longrightarrow X$ such that $T$ is contained in $f^*(\mathrm{Pic}^0(X))$.*

We recall that the morphism in (b) arises as the Stein factorization of the composition $\pi \circ \mathrm{alb}_Y$ where $\pi : \mathrm{Alb}(Y) \to \widehat{T}$ is the map dual to the inclusion $T \hookrightarrow \mathrm{Pic}^0(Y)$.

**Theorem 3.1.10** ([11])**.** *Let $Y$ be a complex smooth projective variety of maximal Albanese dimension and let $f : Y \to Z$ a model of the Iitaka fibration of $Y$. Then the translates through the origin $\mathscr{O}_Y$ of the components of $V^0(Y, \omega_Y)$ generate $\mathrm{Pic}^0(Z)$ as a subabelian variety of $\mathrm{Pic}^0(Y)$. In particular of $Y$ is of general type then they generate the whole $\mathrm{Pic}^0(Y)$.*

## 3.2 Proof of Theorem 3.A

We will prove this slightly more general statement:

**Theorem 3.2.1.** *Let $Z$ be a smooth complex projective variety of maximal Albanese dimension and of general type, then, for every $[\alpha] \in \mathrm{Pic}^0(Z)$ the linear system $|\omega_Z^{\otimes 4} \otimes \alpha|$ is birational.*

Note that also in the setting of Jiang and Chen-Hacon what is proved is the birationality of $\omega_Z^{\otimes 5} \otimes \alpha$ for all $[\alpha] \in \mathrm{Pic}^0(Z)$.
We will need the following easy lemma.

**Lemma 3.2.2.** *Let $E$ an effective divisor on an abelian variety $X$. Take $T_1, \ldots, T_k$ subtori of $X$ such that they generate $X$ as an abstract group and let $\gamma_i$, $i = 1, \ldots, k$ some points of $X$. Then $E \cap (\gamma_i + T_i) \neq \emptyset$ for at least one $i$.*

*Proof.* If $E$ is ample, then the statement is easily seen to be true.
Thus we can suppose, without a loss of generality, that $E$ is not ample. Then it exists an abelian variety $A \neq 0$, a surjective homomorphism $f : X \longrightarrow A$ and an ample divisor $H$ on $A$ such that $E = f^*H$. Suppose now that $E$ does not intersect $\gamma_i + T_i$ for every $i$, it follows that $H$ does not intersect $f(\gamma_i + T_i)$ for every $i$. Hence, since $H$ is ample, the image $\gamma_i + T_i$ through $f$ is a point, and therefore $T_i$ is contained in the kernel of $f$ for every $i$. But this is impossible since the $T_i$'s generate $X$.

□

*Proof of Theorem 3.2.1.*
Let $a_Z : Z \longrightarrow \mathrm{Alb}(Z)$ Albanese map of $Z$. Let us state the following claim.

*Claim 1:* For the generic $z \in Z$ the sheaf $a_{Z\,*}(\mathscr{I}_z \otimes \omega_Z^{\otimes 2} \otimes \mathscr{J}(||\omega_Z||))$ is $M$-regular.

*Remark* 3.2.1. In [60] Pareschi and Popa proved this claim to be true in the case of varieties of Albanese general type, and this was the key point of their proof of the birationality of tricanonical maps (cf. [60, Theorem 6.1]). However, in their argument, the assumption that the Albanese image of $Z$ is not ruled by subtori enters crucially. The proof of this claim is, therefore, the point in which our work diverges from [60].

Before proceeding with the proof of the Claim 1 let us see how it implies Theorem A, this argument follows the one Pareschi and Popa in [60]. First of all we simplify a bit the notation by letting $\mathscr{F} := \omega_Z^{\otimes 2} \otimes \mathscr{J}(||\omega_Z||)$. Now observe that, since the Albanese morphism $a_Z$ is generically finite, a well known extension of Grauert-Riemenschneider vanishing theorem (see, for example [60, proof of Proposition 5.4]) yields that, for any $\alpha \in \mathrm{Pic}^0(Z)$, the higher direct images $R^i a_{Z\,*}(\omega_Z^{\otimes 2} \otimes \alpha \otimes \mathscr{J}(||\omega_Z \otimes \alpha||))$ vanish. Therefore, for every $i \geq 0$ we get the equality

$$V^i(Z, \omega_Z^{\otimes 2} \otimes \alpha \otimes \mathscr{J}(||\omega_Z \otimes \alpha||)) = V^i(\mathrm{Alb}(Z), a_{Z\,*}(\omega_Z^{\otimes 2} \otimes \alpha \otimes \mathscr{J}(||\omega_Z \otimes \alpha||)))$$

By Nadel vanishing for multiplier ideals (Theorem 3.1.2(ii)) the loci at the left hand side are empty when $i > 0$. In particular $a_{Z\,*}(\omega_Z^{\otimes 2} \otimes \alpha \otimes \mathscr{J}(||\omega_Z \otimes \alpha||))$ is $M$-regular on $\mathrm{Alb}(Z)$ and hence cgg everywhere (Proposition 3.1.8). Since $a_Z$ is generically finite, we get that $\omega_Z^{\otimes 2} \otimes \alpha \otimes \mathscr{J}(||\omega_Z \otimes \alpha||)$ is *CGG* away from $W$, the exceptional locus of $a_Z$. This means by definition that the map

$$\bigoplus_{[\beta] \in \Omega} H^0(\omega_Z^{\otimes 2} \otimes \alpha \otimes \mathscr{J}(||\omega_Z \otimes \alpha||) \otimes \beta) \otimes \beta^{-1} \xrightarrow{ev_z} \omega_Z^{\otimes 2} \otimes \alpha \otimes \mathscr{J}(||\omega_Z \otimes \alpha||) \otimes \mathbb{C}(z) \tag{3.2.1}$$

is surjective for every $z \notin W$ and for every $\Omega \subseteq \mathrm{Pic}^0(Z)$ non-empty open set. Now we recall that the base locus of $\mathscr{J}(||\omega_Z \otimes \alpha||)$ is contained in the base locus of $\omega_Z^{\otimes 2} \otimes \alpha^{\otimes 2}$ for every $[\alpha] \in \mathrm{Pic}^0(Z)$ (cfr Lemma 3.1.4); hence we can consider the closed set

$$W' = W \cup \{\text{zero locus of } \mathscr{J}(||\omega_Z \otimes \alpha||)\}.$$

It follows easily from (3.2.1) that, if $z \notin W'$ the map

$$\bigoplus_{[\beta]\in\Omega} H^0(\omega_Z^{\otimes 2} \otimes \alpha \otimes \beta) \otimes \beta^{-1} \xrightarrow{ev_z} \omega_Z^{\otimes 2} \otimes \alpha \otimes \mathbb{C}(z)$$

is surjective for every $\Omega$ open dense set in $\mathrm{Pic}^0(Z)$ and thus $\omega_Z^{\otimes 2} \otimes \alpha$ is CGG away from $W'$.

Afterwards, since we are supposing that the claim is true, we can take $U$ a non empty open set of $Z$ such that for every $z \in U$ the sheaf $a_{Z\,*}(\mathscr{I}_z \otimes \mathscr{F})$ is $M$-regular and hence $CGG$. Again, it follows $\mathscr{I}_z \otimes \mathscr{F}$ is $CGG$ away from $W$. Finally by [61, Proposition 2.12] $\mathscr{I}_z \otimes \mathscr{F} \otimes \omega_Z^{\otimes 2} \otimes \alpha$ is globally generated away from $W'$ and therefore $\mathscr{I}_z \otimes \otimes \omega_Z^{\otimes 4} \otimes \alpha$ is globally generated away from $T = W' \cup \{\text{zero locus of } \mathscr{J}(||\omega_Z||)\}$ and we can conclude that $\omega_Z^{\otimes 4} \otimes \alpha$ is very ample outside a proper subvariety of $Z$ and hence $\varphi_{|\omega_Y^{\otimes 4} \otimes \alpha|}$ is birational for every $[\alpha] \in \mathrm{Pic}^0(Z)$. This prove the Theorem.

Now we proceed with the proof Claim 1. Recall that, by the Subtorus Theorem (Theorem 3.1.9), we can write

$$V^0(\omega_Z) = \bigcup_{i=1}^{k} ([\beta_k] + T_k)$$

where the $[\beta_j]$ are torsion points of $\mathrm{Pic}^0(Z)$ and the $T_j \subseteq \mathrm{Pic}^0(Z)$ subtori. Hence, for all $i$ and for all $[\alpha] \in T_i$ the line bundle $\omega_Z \otimes \alpha \otimes \beta_i$ has sections, and therefore its base locus $Bs(\omega_Z \otimes \alpha \otimes \beta_i)$ is a proper subvariety of Z). We take an non-empty Zariski open set U contained in the complement of

$$W \cup \bigcup_{i=1}^{k} \bigcap_{\alpha\in T_i} \mathrm{Bs}(\omega_Z \otimes \alpha \otimes \beta_i).$$

Given $z \in U$ we want to prove that for every $i \geq 1$

$$\mathrm{Codim}\, V^i(a_{Z\,*}(\mathscr{I}_z \otimes \mathscr{F})) \geq i + 1.$$

Consider the following short exact sequence:

$$(3.2.2) \qquad 0 \to \mathscr{I}_z \otimes \mathscr{F} \otimes \gamma \longrightarrow \mathscr{F} \otimes \gamma \longrightarrow \mathscr{F} \otimes \gamma \otimes \mathbb{C}(z) \to 0.$$

Since $z$ does not belong to the exceptional locus of the Albanese map of $Z$, by pushing forward for $a_Z$ we still get a short exact sequence:

$$(3.2.3) \qquad 0 \to a_{Z\,*}(\mathscr{I}_z \otimes \mathscr{F} \otimes \gamma) \longrightarrow a_{Z\,*}(\mathscr{F} \otimes \gamma) \longrightarrow a_{Z\,*}(\mathscr{F} \otimes \gamma \otimes \mathbb{C}(z)) \to 0.$$

In particular $R^1 a_{Z\,*}(\mathscr{I}_z \otimes \mathscr{F} \otimes \gamma) = 0$. Since, as already observed, an extension of Grauert-Riemenschneider vanishing, yields the vanishing the $R^j a_{Z\,*}(\omega_Z^{\otimes 2} \otimes \mathscr{J}(||\omega_Z||) \otimes \gamma)$ for every $j > 0$ and every $[\gamma] \in \mathrm{Pic}^0(Z)$, we have that $R^j a_{Z*}(\mathscr{I}_z \otimes \mathscr{F} \otimes \gamma) = 0$ for every $j > 0$ and every $\gamma$. Therefore, for every $j \geq 0$

$$V^i(a_{Z*}(\mathscr{I}_z \otimes \mathscr{F})) = V^i(\mathscr{I}_z \otimes \mathscr{F}) \tag{3.2.4}$$

Now we look at the long cohomology sequence of the short exact sequence (3.2.2). By Nadel vanishing (Theorem 3.1.2) we have that, for every $i \geq 1$ and $[\gamma] \in \mathrm{Pic}^0(Z)$,

$$H^i(Z, \omega_Z^{\otimes 2} \otimes \mathscr{J}(||\omega_Z||) \otimes \alpha) = 0$$

and therefore, by (3.2.4) the loci $V^i(a_{Z*}(\mathscr{I}_z \otimes \mathscr{F}))$ are empty for every $i \geq 2$.
It remains to see that $\mathrm{Codim} V^1(a_{Z*}(\mathscr{I}_z \otimes \mathscr{F})) > 1$. Before proceeding further we state the following:

*Claim 2:* For every $j = 1, \ldots, k$, $V^1(\mathscr{I}_z \otimes \mathscr{F}) \cap (T_j + [\beta_j^{\otimes 2}]) = \emptyset$.

The reader may observe that by Theorem 3.1.10 the $T_j$'s generate $\mathrm{Pic}^0(Z)$, and this, together with Lemma 3.2.2, implies directly Claim 1. In fact the locus $V^1(\mathscr{I}_z \otimes \mathscr{F}))$ is a proper subvariety of $\mathrm{Pic}^0(Z)$ that cannot possibly contain a divisor, hence its codimension shall be greater than 1. By (3.2.4) also $\mathrm{Codim} V^1(a_{Z*}(\mathscr{I}_z \otimes \mathscr{F})) > 1$.

To prove Claim 2 we reason in the following way. First we remark that there is an alternative description of the locus $V^1(\mathscr{I}_z \otimes \mathscr{F})$: in fact by (3.2.2) and the vanishing of $H^1(Z, \mathscr{F} \otimes \gamma)$ for every $[\gamma]$ (again due to Nadel vanishing) we have that $H^1(Z, \mathscr{I}_z \otimes \mathscr{F} \otimes \gamma) \neq 0$ if and only if the map $H^0(Z, \mathscr{F} \otimes \gamma) \longrightarrow \mathscr{F} \otimes \gamma \otimes \mathbb{C}(z)$ is not surjective. Therefore we can write

$$V^1(\mathscr{I}_z \otimes \mathscr{F}) = \{[\gamma] \in \mathrm{Pic}^0(Z) \text{ such that } z \notin \mathrm{Bs}(\mathscr{F} \otimes \gamma)\}.$$

Since we chose $z \notin \bigcap_{\alpha \in T_j} \mathrm{Bs}(\omega_Z \otimes \alpha \otimes \beta_i)$ it exists $V_j$ a dense open set of $T_j$ such that, for every $[\delta] \in V_j$, $z$ is not in $\mathrm{Bs}(\omega_Z \otimes \beta_i \otimes \delta)$. Take $[\gamma] \in T_j$; the intersection $V_j \cap ([\gamma] - V_j)$ is still a non-empty open set of $T_j$. For $[\delta] \in V_j \cap ([\gamma] - V_j)$ we write

$$\omega_Z^{\otimes 2} \otimes \beta_j^{\otimes 2} \otimes \gamma \simeq (\omega_Z \otimes \beta_j \otimes \delta) \otimes (\omega_Z \otimes \beta_j \otimes \gamma \otimes \delta^{-1});$$

since, thanks to our choice of $\delta$, both $\omega_Z \otimes \beta_j \otimes \delta$ and $\omega_Z \otimes \beta_j \otimes \gamma \otimes \delta^{-1}$ are generated in $z$, also the left hand side is. Therefore $z \notin \mathrm{Bs}(\omega_Z^{\otimes 2} \otimes \beta_j^{\otimes 2} \otimes \gamma)$. Now we use again Lemma 3.1.4: the zero locus of $\mathscr{J}(||\omega_Z||)$ is contained in the base locus of $\omega_Z^{\otimes 2} \otimes \alpha$ for every $[\alpha] \in \mathrm{Pic}^0(Z)$; thus we have the equality $\mathrm{Bs}(\omega_Z^{\otimes 2} \otimes \beta_j^{\otimes 2} \otimes \delta) = \mathrm{Bs}(\mathscr{F} \otimes \beta_j^{\otimes 2} \otimes \delta)$ and, consequently, $\delta \otimes \beta_j^{\otimes 2} \notin V^1(\mathscr{I}_z \otimes \mathscr{F})$ and the claim is proved. □

# CHAPTER 4

# EFFECTIVE IITAKA FIBRATIONS

The results presented in this chapter are joint work with Z. Jiang and M. Lahoz ([39]) and are an improvement of the material exposed in the previous Chapter and published in [72]. In particular, our main achievement is the following:

**Theorem 4.A.** *Let $X$ be a smooth projective variety of maximal Albanese dimension. Then,*

1. *the linear system $|4K_X|$ induces the Iitaka fibration of $X$;*
2. *if $X$ is of general type, then the linear system $|3K_X|$ induces a birational map.*

This is clearly the optimal bound for varieties of maximal Albanese dimension. On one hand, the varieties of general type whose bicanonical map is not birational have been studied in [2, 48]. On the other hand, we produce varieties of dimension at least 4, whose tricanonical map does not induce the Iitaka fibration (see Example 4).

We observe that when $\chi(X, \omega_X) > 0$, the birationality of the tricanonical map was proved by Chen and Hacon [9, Thm. 5.4]. Hence, we restrict ourselves to the case $\chi(X, \omega_X) = 0$. In that situation we have a special fibration where the $m$-th pluricanonical linear system restricts surjectively to the general fiber for $m \geq 3$ (see Lemma 4.3.5). On the base of this fibration, we construct two positive line bundles (see Lemmas 4.3.3 and 4.3.4). One of them induces a birational map on the base, and the other one is used to prove the appropiate effectiveness. We use them to apply Lemma 4.2.1, which allows us to apply induction on the dimension of $X$. The lack of effectivity is what forces us to consider the tetracanonical map for non-general type varieties (as we note in Remark 4.5.2).

### Notation

In the sequel, $X$ will always be a smooth complex variety of maximal Albanese dimension. We denote by $a_X : X \to A_X$ the Albanese morphism of a smooth projective variety.

We will not distinguish in this Chapter between line bundles and divisors on $X$. Thus given two line bundles $L$ and $M$ we will denote the twist of $M$ and $L$ either by $L \otimes M$ or by $L + M$. The dual of an invertible sheaf will be denoted either by $L^\vee$ or by $-L$. The reason for doing this is that sometimes the notation with tensor products is cumbersome.

Given a smooth variety $X$, by $P_m(X)$ we will mean the $m$-th plurigenus of $X$, i.e. $h^0(X, \omega_X^{\otimes m})$.

## 4.1 Background Material

In this section we present some results we will be needing afterwards.

### Vanishing Theorems

We begin with two vanishing results. The first one tells us when certain higher direct image are zero.

**Lemma 4.1.1** ([50, Lemma 4.3.10]). *Let $h : X \to Z$ be a morphism of irreducible projective varieties and let $H$ be a very ample divisor on $Z$. Suppose that $\mathscr{G}$ is a coherent sheaf on $X$ with the property that*

$$H^i(X, \mathscr{G} \otimes h^*(H^{\otimes m})) = 0 \quad \textit{for every } i \geq 1 \textit{ and every } m >> 0.$$

*Then $R^j h_*(\mathscr{G}) = 0$ for every $j > 0$.*

The second is a very important vanishing result of Kawamata-Viehweg that extends Kodaira vanishing to big and nef $\mathbb{Q}$-divisors.

**Theorem 4.1.2** (Kawamata–Viehweg vanishing Theorem[50, §15]). *Let $X$ be a smooth projective variety of dimension $n$, and let $N$ be an integral divisor on $X$. Assume that $N$ is numerically equivalent to a sum $B + \Delta$ with $B$ a nef and big divisor and $\Delta = \sum_i a_i \Delta_i$ a $\mathbb{Q}$-divisor with simple normal crossing support and fractional coefficents*

$$0 \leq a_i < 1 \quad \textit{for every } i.$$

*Then for every $i > 0$*

$$h^i(X, \mathscr{O}_X(K_X + N)) = 0$$

### A result of Barja–Lahoz–Naranjo–Pareschi

In [2] the authors introduced the following:

**Hypothesis 4.1.3** ([2, Hypothesis 4.7])**.** Let $X$ be a variety of dimension $n$ and maximal Albanese dimension equipped with a generically finite morphism to an abelian variety $a : X \to A$ such that the map $a^* : \widehat{A} \to \mathrm{Pic}^0(X)$ is an embedding and

$$\mathrm{codim} V^i_a(X, \omega_X) \geq i + 1 \quad \text{for every } 0 < i < n$$

The reason for recalling this setting is the following result of [2]:

**Proposition 4.1.4** ([2, Proposition 4.10])**.** *Let $X$ be a smooth projective variety satisfying Hypothesis 4.1.3. Then $\chi(\omega_X) = 0$ if and only if $X$ is not of general type. Furthermore, if this is the case, $X$ is birational to an abelian variety.*

### Some Technical Lemmas of Jiang

We conclude this section by listing some results from [37] that we will be using.

Given a surjective morphism of smooth projective varieties, $f : X \to Y$, and $D$ a $\mathbb{Q}$-divisor on $X$ we will say that the *Iitaka model of* $(X, D)$ *dominates* $Y$ if it exists a positive integer $N$, and an ample divisor $H$ on $Y$ such that $ND - f^*H$ is an effective divisor.

**Lemma 4.1.5** ([37, Lemma 2.1])**.** *Suppose that $F : X \to Y$ is as above. Let $L$ be a $\mathbb{Q}$ divisor on $X$ such that the Iitaka model of $(X, L)$ dominates $Y$, and let $D$ a nef $\mathbb{Q}$-divisor on $Y$ such that $L + f^*D$ is an integral divisor on $X$. Then for every $i \geq 1$, every $j \geq 0$ and every $[\alpha] \in \mathrm{Pic}^0(X)$*

$$h^i(Y, R^j f_*(\mathscr{O}_X(K_X + L + f^*D)) \otimes \mathscr{J}(||L||) \otimes \alpha)) = 0$$

**Lemma 4.1.6** ([37, Lemma 2.3])**.** *Suppose that $f : X \to Y$ is an algebraic fiber space between smooth projective varieties and suppose that the $m$-th plurigenus of $X$, $P_m(X)$, in non zero for some $m \geq 2$. Assume that $H$ is a big $\mathbb{Q}$-divisor on $Y$. Then, the Iitaka model of $(X, (m-1)K_{X/Y} + f^*H)$ dominates $Y$, for any $m \geq 2$. Furthermore denoting by $\mathscr{F}$ the sheaf*

$$f_*(\mathscr{O}_X(K_X + (m-1)K_{X/Y}) \otimes \mathscr{J}(||(m-1)K_{X/Y} + f^*H||)$$

*we get that $\mathscr{F}$ is a non zero sheaf of rank $P_m(X_y)$ with $X_y$ the general fiber of $f$.*

Observe that if $Y$ is of general type, then we can take $H = K_Y$ and thus we get

1. the Iitaka model of $(X, K_X + (m-2)K_{X/Y})$ dominates $Y$;

2. $f_*(\mathscr{O}_X(2K_X+(m-2)K_{X/Y}\otimes \mathscr{J}(||K_X+(m-2)K_{X/Y}||)$ is a non zero sheaf of rank $P_m(X_y)$.

Now let $g : X \to Z$ be an algebraic fiber space between smooth projective varieties such that is not zero for some $m \geq 2$. Fix $H$ a big, base point free divisor on $Z$. Consider the following ideal sheaves

$$\mathscr{J}_{m-1,n} := \mathscr{J}(||(m-1)K_{X/Z} + \frac{1}{n}g^*H||)$$

**Lemma 4.1.7** ([37, Lemma 2.4]). *The following inclusions hold*

$$\mathscr{J}_{m-1,n} \supseteq \mathscr{J}_{m-1,n+1},$$

*and, in addition, it exists an integer $N$ such that for every $n \geq N$ there is an equality*

$$\mathscr{J}_{m-1,n} = \mathscr{J}_{m-1,N}.$$

In the above assumptions introduce the following sheaves:

$$\mathscr{F}_{m-1,H} := g_*(\mathscr{O}_X(K_X + (m-1)K_{X/Z}) \otimes \mathscr{J}_{m-1,N}.$$

**Lemma 4.1.8** ([37, Lemma 2.5]). *In the above assumptions and notation, suppose furthermore that $Z$ is of maximal Albanese dimension and that $H$ is a big, base point free divisor pulled back from* Alb$(Z)$*. Then $\mathscr{F}_{m-1,H}$ is a non zero GV sheaf.*

## 4.2 Preliminaries

We begin with some easy lemmas.

**Lemma 4.2.1.** *Let $f : X \xrightarrow{g} Z \xrightarrow{h} Y$ be fibrations between smooth projective varieties. Let $L$ be an line bundle on $X$. If the following two conditions hold:*

*1) The image of $H^0(X,L) \to H^0(X_y, L|_{X_y})$ induces a map birationally equivalent to $g_{|X_y} : X_y \to Z_y$ for a general fiber $X_y$ of $f$;*

*2) There are line bundles $H_i$, $1 \leq i \leq M$, on $Y$ such that $L - f^*H_i$ is effective and the multiple evaluation map*

$$\varphi_Y : Y \to \mathbf{P}(H^0(Y,H_1)^*) \times \cdots \times \mathbf{P}(H^0(Y,H_M)^*)$$

*is birational.*

*Then, the linear system $|L|$ induces a map birationally equivalent to $g : X \to Z$.*

*Proof.* Since the $L - f^*H_i$ are effective we have a projection $\pi$ that induces the following diagram

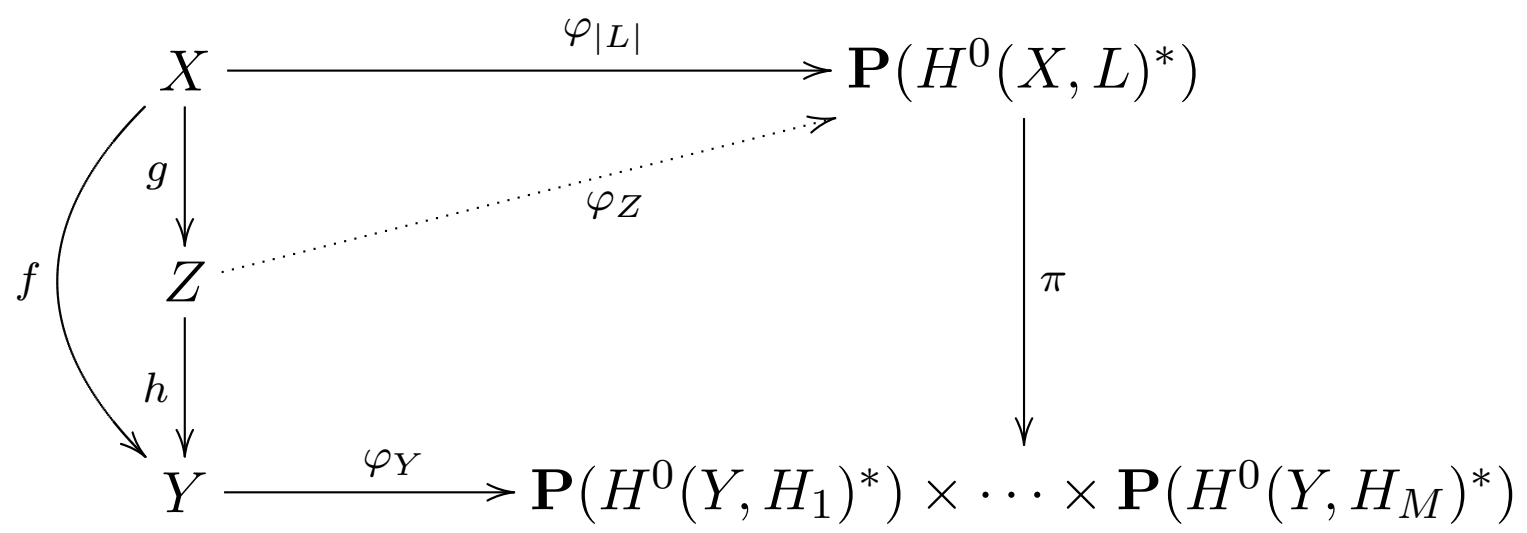

Condition 2) guarantees that $\varphi_{|L|}$ separates generic fibres of $f$ and condition 1) shows that the map $\varphi_{|L|}$ factorizes to $\varphi_Z \circ g$ and a general fiber of $h$ is mapped birationally via $\varphi_Z$. □

We will need the following lemma to ensure the birationality of $\varphi_Y$ in the previous lemma.

**Lemma 4.2.2.** *Let $\widehat{\pi} : \widehat{Y} \to Y$ be an abelian cover with Galois group $G$ of smooth projective varieties of maximal Albanese dimensions. We denote $b_{\widehat{Y}} = a_Y \circ \widehat{\pi} : \widehat{Y} \to A_Y$. Assume that $V^0_{b_{\widehat{Y}}}(\widehat{Y}, \omega_{\widehat{Y}}) = \mathrm{Pic}^0(Y)$ and*

$$\widehat{\pi}_*\omega_{\widehat{Y}}^{\otimes 2} = \bigoplus_{\chi \in G^*} \mathcal{H}_\chi,$$

*where $\mathcal{H}_\chi$ is the line bundle corresponding to the character $\chi \in G^*$.*

*Then, there exists $\mathcal{H}_{\chi_0}$ such that the multiple evaluation map*

$$\varphi_{[\alpha_1]\cdots[\alpha_M]} : Y \to \mathbf{P}(H^0(Y, \mathcal{H}_{\chi_0} \otimes [\alpha_1])^*) \times \cdots \times \mathbf{P}(H^0(Y, \mathcal{H}_{\chi_0} \otimes [\alpha_M])^*)$$

*is birational for some $[\alpha_i] \in \mathrm{Pic}^0(Y)$, $1 \leq i \leq M$.*

If $\widehat{\pi}$ is an isomorphism, then Lemma 4.2.2 is essentially a part of the proof of [9, Thm. 4.4].

*Proof.* We first write

$$\widehat{\pi}_*\omega_{\widehat{Y}} = \bigoplus_{\chi \in G^*} \mathcal{L}_\chi.$$

Since $V^0(\omega_{\widehat{Y}}, b_{\widehat{Y}}) = \mathrm{Pic}^0(Y)$, we conclude that there exists $\chi_1$ such that $V^0_{a_Y}(Y, \mathcal{L}_{\chi_1}) = \mathrm{Pic}^0(Y)$. Then, for any $[\alpha] \in \mathrm{Pic}^0(Y)$, the line bundle $\mathcal{L}_{\chi_1}^{\otimes 2} \otimes \alpha$ is globally generated on the open dense subset

$$Y \setminus \cap_{[\alpha] \in \mathrm{Pic}^0(Y)} \mathrm{Bs}(|\mathcal{L}_{\chi_1} \otimes \alpha|).$$

Since there is the natural $G$-map of vector bundles on $Y$,

$$(\widehat{\pi}_*\omega_{\widehat{Y}})^{\otimes 2} \to \widehat{\pi}_*\omega_{\widehat{Y}}^{\otimes 2},$$

if we take $\chi_0 = \chi_1^2$, then we have an inclusion $\mathcal{L}_{\chi_1}^{\otimes 2} \hookrightarrow \mathcal{H}_{\chi_0}$. Hence, there is an open dense subset $U$ of $Y$ such that for any $[\alpha] \in \mathrm{Pic}^0(Y)$, the line bundle $\mathcal{H}_{\chi_0} \otimes \alpha$ is globally generated on $U$, i.e. $\varphi_{\alpha|U}$ is a morphism.

On the other hand, we consider $\widehat{\pi}_*(\omega_{\widehat{Y}}^{\otimes 2} \otimes \mathscr{J}(||\omega_{\widehat{Y}}||))$. Since $\mathscr{J}(||\omega_{\widehat{Y}}||) \hookrightarrow \mathscr{O}_{\widehat{Y}}$ is $G$-invariant, we can write

$$\widehat{\pi}_*(\omega_{\widehat{Y}}^{\otimes 2} \otimes \mathscr{J}(||\omega_{\widehat{Y}}||)) = \bigoplus_{\chi \in G^*} \mathcal{H}_\chi \otimes \mathscr{I}_\chi,$$

where $\mathscr{I}_\chi$ is an ideal sheaf on $Y$. Moreover, we have

$$H^0(Y, \mathcal{H}_\chi \otimes \mathscr{I}_\chi) \simeq H^0(Y, \mathcal{H}_\chi),$$

and

$$H^i(Y, \mathcal{H}_\chi \otimes \mathscr{I}_\chi \otimes \alpha) = 0,$$

for any $[\alpha] \in \mathrm{Pic}^0(Y)$ and $i \geq 1$.

Therefore, $\mathcal{H}_{\chi_0} \otimes \mathscr{I}_{\chi_0} \otimes \alpha$ is globally generated on the open subset $U$, so $\mathrm{Supp}(\mathscr{O}_Y/\mathscr{I}_{\chi_0})$ is contained in $Y \backslash U$. Now, let $V = U \backslash \mathrm{Exc}(a_Y)$ and for any point $y \in V$, from the exact sequence

$$0 \to \mathscr{I}_y \otimes \mathcal{H}_{\chi_0} \otimes \mathscr{I}_{\chi_0} \to \mathcal{H}_{\chi_0} \otimes \mathscr{I}_{\chi_0} \to \mathbf{C}_y \to 0,$$

we see that $a_{Y*}(\mathscr{I}_y \otimes \mathcal{H}_{\chi_0} \otimes \mathscr{I}_{\chi_0})$ is a $M$-regular sheaf, so it is *cgg* (see 3.1.7). Hence for any $z \in V$ different from $y$, there exists $[\alpha] \in \mathrm{Pic}^0(Y)$ such that $\mathscr{I}_y \otimes \mathcal{H}_{\chi_0} \otimes \mathscr{I}_{\chi_0} \otimes \alpha$ is globally generated on $z$.

This shows that for any two different points $y, z \in V$ there exists $[\alpha] \in \mathrm{Pic}^0(Y)$ and a divisor $D_\alpha$ in $|\mathcal{H}_{\chi_0} \otimes \alpha|$ such that $y \in D_\alpha$ but $z \notin D_\alpha$. Therefore $\varphi_\alpha(y) \neq \varphi_\alpha(z)$.

We take $\alpha_1, \ldots, \alpha_M$ such that $\varphi_{\alpha_1 \cdots \alpha_M}$ becomes stable, namely $\varphi_{\alpha_1 \cdots \alpha_M}$ is birational equivalent to $\varphi_{\alpha_1 \cdots \alpha_M P}$ for any $P \in \mathrm{Pic}^0(Y)$. Then, $\varphi_{\alpha_1 \cdots \alpha_M}$ is birational. □

The following lemma should be compared to [8, Lem. 3.1].

**Lemma 4.2.3.** *Let $f : X \to Y$ be a surjective morphism between smooth projective varieties. Assume that $X$ is of maximal Albanese dimension. Then, $K_{X/Y}$ is effective.*

*Proof.* We have the natural inclusion $f^*\Omega_Y^1 \xrightarrow{i} \Omega_X^1$. Denote by $\mathscr{F}$ the saturation of $i(f^*\Omega_Y^1)$. Then, $\det(\mathscr{F}) - f^*K_Y$ is an effective divisor on $X$. We then consider the exact sequence

$$0 \to \mathscr{F} \to \Omega_X^1 \to \mathcal{Q} \to 0.$$

Since $X$ is of maximal Albanese dimension, $\Omega_X^1$ is generically globally generated and so is $\mathcal{Q}$.

Hence $\det(\mathcal{Q})$ is also an effective divisor. Hence $K_{X/Y} = \det(\mathscr{F}) - f^*K_Y + \det(\mathcal{Q})$ is effective. □

## 4.3 Positive Bundles on the Base and Surjectivity of the Restriction Map to a Fiber

We will use the following definition that it is strongly related to Hypotheses 4.1.3.

**Definition 4.3.1.** Let $\mathscr{F}$ be a coherent sheaf on an abelian variety $A$. We say that $\mathscr{F}$ is almost $M$-regular if $V^0(A,\mathscr{F}) = \widehat{A}$, $\operatorname{codim}_{\widehat{A}} V^i(A,\mathscr{F}) \geq i+1$, for $1 \leq i \leq \dim A - 1$, and $\dim V^{\dim A}(A,\mathscr{F}) = 0$.

Let $X$ be a smooth projective variety of dimension $n$ and maximal Albanese dimension. We know that the pushforward of the canonical bundle $a_{X*}\omega_X$ is a $GV$-sheaf but it often fails to be $M$-regular, which makes the tricanonical map difficult to study. Hence we consider the set

$$\mathcal{S}_X := \{0 < j < n \mid V^j_{a_X}(X,\omega_X) \text{ has a component of codimension } j\},$$

which measures how far $a_{X*}\omega_X$ is away from being almost $M$-regular.

**Setting 4.3.2.** Assume that $\mathcal{S}_X$ is not empty. We denote by $k$ the maximal number of $\mathcal{S}_X$ and $[\beta] + \widehat{B} \subseteq V^k_{a_X}(X,\omega_X)$ a codimension-$k$ component, where $[\beta]$ is zero or a torsion element of $\mathrm{Pic}^0(X)\backslash\widehat{B}$. Let the following commutative diagram

$$\begin{array}{ccc} X & \xrightarrow{a_X} & A_X \\ {\scriptstyle f}\downarrow & & \downarrow{\scriptstyle \mathrm{pr}} \\ Y & \xrightarrow{a_Y} & B \end{array}$$

be a suitable birational modification of the Stein factorization of the composition $\mathrm{pr} \circ a_X$, such that $Y$ is smooth.

**Lemma 4.3.3.** *Assume $\mathcal{S}_X$ is not empty, so we are in Setting 4.3.2. Then, in some birational model of $f : X \to Y$, there exists a line bundle $\mathcal{L}$ on $Y$ such that $a_{Y*}\mathcal{L}$ is almost $M$-regular, $V^0_{a_Y}(Y,\mathcal{L}) = \widehat{B}$, and $\mathscr{O}_X(K_X) \otimes \beta^{\otimes j} \otimes f^*\mathcal{L}^\vee$ has a non-trivial section for some $j \in \mathbf{Z}$. Moreover,*

*1) if $\beta$ is trivial, we can take $\mathcal{L}$ to be $\omega_Y$ and $j = 0$;*

*2) if $[\beta] \in \mathrm{Pic}^0(X)\backslash\widehat{B}$, then we can take $\mathcal{L}$ such that $a_{Y*}\mathcal{L}$ is $M$-regular.*

*Proof.* We know by [29, Theorem 0.1] that the dimension of a general fiber of $f$ is $k$.

Assume first that $\beta$ is trivial. For $[\alpha_B] \in \widehat{B} - \cup_j V^1(R^j f_*\omega_X)$

$$\begin{aligned} 0 < h^k(X,\omega_X \otimes a_X^*\alpha_B) &= h^k(X,\omega_X \otimes f^*\alpha_B) = h^0(Y, R^k f_*\omega_X \otimes \alpha_B) \\ &= h^0(Y,\omega_Y \otimes \alpha_B). \end{aligned}$$

So $\chi(Y, \omega_Y) > 0$ by semicontinuity. Hence $V^0_{a_Y}(\omega_Y) = \mathrm{Pic}^0(Y)$.

Moreover, the pull-back by pr of any codimension-$j$ component of $V^j_{a_Y}(Y, \omega_Y)$ is a codimension-$(j+k)$ component of $V^{j+k}_{a_X}(X, \omega_X)$. Hence by the maximality of $k$, we know that $\mathrm{codim}_{\widehat{B}} V^i(\omega_Y, a_Y) \geq i+1$ for all $0 < i < \dim Y$ and $a_{Y*}\omega_Y$ is almost $M$-regular. By Lemma 4.2.3, $\mathscr{O}_X(K_X - f^*K_Y)$ is effective.

Now, assume that $[\beta] \in \mathrm{Pic}^0(X) \backslash \widehat{B}$. We may choose $[\beta]$ such that $\langle[\beta]\rangle \cap \widehat{B} = 0$, where $G := \langle[\beta]\rangle$ is the subgroup generated by $[\beta]$. We then consider the étale cover $\pi : \widetilde{X} \to X$ induced by $G$ and, after modifications, we have the following diagram

(4.3.1)
$$\begin{array}{ccccc} \widetilde{X} & \xrightarrow{\pi} & X & \xrightarrow{a_X} & A_X \\ \downarrow{\scriptstyle \widetilde{f}} & & \downarrow{\scriptstyle f} & & \downarrow{\scriptstyle \mathrm{pr}} \\ \widehat{Y} & \xrightarrow{\widehat{\pi}} & Y & \xrightarrow{a_Y} & B \end{array} \qquad b_{\widetilde{X}} : \widetilde{X} \to A_X, \quad b_{\widehat{Y}} : \widehat{Y} \to B$$

where all varieties are smooth and the vertical morphisms are fibrations.

The same arguments in the previous case show that $V^0_{b_{\widehat{Y}}}(\widehat{Y}, \omega_{\widehat{Y}}) = \widehat{B}$ and $\mathrm{codim}_{\widehat{B}} V^i_{b_{\widehat{Y}}}(\widehat{Y}, \omega_{\widehat{Y}}) \geq i+1$ for all $0 < i < \dim \widehat{Y}$.

Since $[\beta] + \widehat{B} \subset V^k_{a_X}(X, \omega_X)$, $R^k f_*(\omega_X \otimes Q) \neq 0$. If we denote by $X_y$ a general fiber of $f$, then we know that $\beta_{|X_y}$ is $\mathscr{O}_{X_y}$. Hence

$$\widehat{\pi}_*\omega_{\widehat{Y}} = R^k(f \circ \pi)_*\omega_{\widetilde{X}} = \oplus_i (R^k f_*(\omega_X \otimes \beta^{\otimes i})) = \oplus_i \mathcal{L}_i,$$

where $\mathcal{L}_0 = \omega_Y$ and after modifications, we may assume that all $\mathcal{L}_i$ are line bundles on $Y$. Moreover, we have

$$\begin{aligned} \ker(\widehat{B} \to \mathrm{Pic}^0(\widehat{Y})) &= \widehat{B} \cap \ker(\widehat{A}_X \to \widehat{A}_{\widetilde{X}}) \\ &= \widehat{B} \cap \langle[\beta]\rangle = 0. \end{aligned}$$

Hence for $i \neq 0$, $V^{\dim Y}_{a_Y}(Y, \mathcal{L}_i) = \emptyset$ and $a_{Y*}\mathcal{L}_i$ is $M$-regular.

We know $K_{\widetilde{X}/\widehat{Y}}$ is effective by Lemma 4.2.3 and we have a nonzero map $\widehat{\pi}^*\mathcal{L}_i \to \omega_{\widehat{Y}}$ for every $i$. Hence,

$$\begin{aligned} 0 &< h^0(\widetilde{X}, \omega_{\widetilde{X}} \otimes \widetilde{f}^*\omega^{\vee}_{\widehat{Y}}) \leq h^0(\widetilde{X}, \omega_{\widetilde{X}} \otimes \widetilde{f}^*\widehat{\pi}^*\mathcal{L}^{\vee}_i) \\ &= \sum_j h^0(X, \omega_X \otimes \beta^{\otimes j} \otimes f^*\mathcal{L}^{\vee}_i). \end{aligned}$$

So there exists $j$ such that $\omega_X \otimes \beta^{\otimes j} \otimes f^*\mathcal{L}^{\vee}_i$ is effective. Thus in 2), we can take $\mathcal{L}$ to be any $\mathcal{L}_i$ for $i \neq 0$. $\square$

**Lemma 4.3.4.** *Assume $\mathcal{S}_X$ is not empty, so we are in Setting 4.3.2. Let $\mathcal{L}$ be the line bundle obtained in Lemma 4.3.3. Then, after modifying $Y$, there exists a line bundle $\mathcal{H}$ on $Y$ and $i \in \mathbf{Z}$ such that $\mathcal{H} \otimes \mathcal{L}^{\otimes -2}$ and $\mathscr{O}_X(2K_X) \otimes \beta^{\otimes i} \otimes f^*\mathcal{H}^\vee$ have non-trivial sections. Moreover we can take $[\alpha_i] \in \mathrm{Pic}^0(Y)$, $1 \leq i \leq M$ such that the multiple evaluation map*

$$\varphi_Y : Y \to \mathbf{P}(H^0(Y, \mathcal{H} \otimes \alpha_1)^*) \times \cdots \times \mathbf{P}(H^0(Y, \mathcal{H} \otimes \alpha_M)^*)$$

*is birational.*

*Proof.* If $\beta$ is trivial, we just take $\mathcal{H}$ to be $\omega_Y^{\otimes 2}$. By Lemma 4.2.2, we conclude.

If $[\beta] \in \mathrm{Pic}^0(X) \backslash \widehat{B}$, then we use the same notations as in the proof of Lemma 4.3.3. Notice that in diagram (4.3.1) we have $\beta_{|X_y}$ is the trivial bundle. Hence, a general fiber of $\widetilde{f}$ is birational via $\pi$ to a general fiber of $f$. So $G$ is also the Galois group of the field extension $k(\widehat{Y})/k(Y)$. After modifications of $\widetilde{f} : \widetilde{X} \to \widehat{Y}$, we may assume that $G$ acts also on $\widehat{Y}$, $\widetilde{f}$ is a $G$-equivariant morphism, and $\widehat{\pi} : \widehat{Y} \to Y$ is a $G$-cover of smooth projective varieties.

Now, we apply Lemma 4.2.2 to $\widehat{\pi}$. Take $\mathcal{H}$ to be the direct summand of $\widehat{\pi}_*\omega_{\widehat{Y}}^{\otimes 2}$ which contains $\mathcal{L}^2$.

As in the proof of Lemma 4.3.3, since $\widehat{\pi}^*\mathcal{H} \hookrightarrow \omega_{\widehat{Y}}^{\otimes 2}$ and $K_{\widetilde{X}/\widehat{Y}}$ is effective, we conclude that there exists $i \in \mathbf{Z}$ such that $\mathscr{O}_X(2K_X) \otimes \beta^{\otimes i} \otimes f^*\mathcal{H}^\vee$ is effective. □

**Lemma 4.3.5.** *Assume $\mathcal{S}_X$ is not empty, so we are in Setting 4.3.2. Let $\mathcal{L}$ be the line bundle obtained in Lemma 4.3.3. Then, for $y$ a general point of $Y$, the restriction map*

$$H^0(X, \mathscr{O}_X(mK_X - (m-3)f^*\mathcal{L}) \otimes \alpha) \to H^0(X_y, \mathscr{O}_{X_y}(mK_{X_y}) \otimes \alpha)$$

*is surjective, for any $m \geq 2$ and $[\alpha] \in V^0_{a_X}(X, \omega_X^{\otimes m})$.*

*Proof.* We just prove the statement for $\alpha = \mathscr{O}_X$, the same arguments works for the general case.

There are two distinguished cases, whether $Q$ is trivial or not, which we deal with slightly different techniques.

**Case A**. Assume that $\beta$ is trivial.

We have seen in the proof of Lemma 4.3.3 that $\chi(Y, \omega_Y) > 0$, so $Y$ is of general type. By Lemma 4.1.6, the Iitaka model of $(X, K_X + (m-2)K_{X/Y})$ dominates $Y$ and by Lemma 4.1.5 there exists an asymptotic multiplier ideal sheaf $\mathscr{I} := \mathscr{J}(||K_X + (m-2)K_{X/Y}||)$ on $X$ such that $a_{Y*}f_*(\mathscr{O}_X(2K_X + (m-2)K_{X/Y}) \otimes \mathscr{I})$ is a sheaf satisfying I.T. with index 0. In particular it is $M$-rgular and hence, by 3.1.7, *cgg*. It follows that

$$\mathscr{F} := f_*(\mathscr{O}_X(2K_X + (m-2)K_{X/Y}) \otimes \mathscr{I})$$

is *cgg* outside the exceptional locus of $a_Y$.

We conclude, similarly to [2, Prop. 4.4. and Cor. 4.11] that, $\mathscr{F} \otimes \omega_Y$ is globally generated in an open dense subset of $Y$. Indeed, we first notice that $|\omega_Y \otimes \alpha|$ is not empty for all $[\alpha] \in \widehat{B}$. Take $\alpha_i$, $1 \leq i \leq N$ such that the evaluation map

$$\oplus_{i=1}^{N} H^0(Y, \mathscr{F} \otimes \alpha_i) \otimes \alpha_i^{\vee} \to \mathscr{F}$$

is surjective. Then

$$\oplus_{i=1}^{N} H^0(Y, \mathscr{F} \otimes \alpha_i) \otimes H^0(Y, \omega_Y \otimes \alpha_i^{\vee}) \otimes \mathscr{O}_Y \to \mathscr{F} \otimes \omega_Y$$

is surjective over $Y \setminus \cup_{1 \leq i \leq N} \mathrm{Bs}(|\omega_Y \otimes \alpha_i|)$. Finally this evaluation map factors through $H^0(Y, \mathscr{F} \otimes \omega_Y) \otimes \mathscr{O}_Y \to \mathscr{F} \otimes \omega_Y$.

Moreover, by the second part of Lemma 4.1.6, $\mathscr{F}$ is a non-zero sheaf on $Y$ of rank $\alpha_m(X_y)$. Hence, over a general point $y \in Y$, $\mathscr{F}|_y$ is isomorphic to $H^0(X_y, \mathscr{O}_{X_y}(mK_{X_y}))$. Since

$$\mathscr{F} \otimes \omega_Y \subset f_*(\mathscr{O}_X(mK_X - (m-3)f^*K_Y))$$

and they have the same rank $\alpha_m(X_y)$, we conclude the proof of the lemma when $\beta$ is trivial.

**Case B**. If $\beta$ is non-trivial, we use the same commutative diagram as in the proof of Lemma 4.3.3:

$$\begin{array}{ccccc} \widetilde{X} & \xrightarrow{\pi} & X & \xrightarrow{a_X} & A_X \\ \downarrow{\scriptstyle \widetilde{f}} & & \downarrow{\scriptstyle f} & & \downarrow{\scriptstyle \mathrm{pr}} \\ \widehat{Y} & \xrightarrow{\widehat{\pi}} & Y & \xrightarrow{a_Y} & B \end{array} \qquad b_{\widetilde{X}} : \widetilde{X} \to A_X, \quad b_{\widehat{Y}} : \widehat{Y} \to B$$

We claim that the Iitaka model of $(X, (m-1)K_X - (m-2)f^*\mathcal{L}_i)$ dominates $Y$. Indeed, by Lemma 4.1.6, the Iitaka model of $(\widetilde{X}, (m-1)K_{\widetilde{X}/\widehat{Y}} + \widetilde{f}^*K_{\widehat{Y}})$ dominates $\widehat{Y}$. Since $(m-1)K_{\widetilde{X}/\widehat{Y}} + \widetilde{f}^*K_{\widehat{Y}} \preceq \pi^*((m-1)K_X - (m-2)f^*\mathcal{L}_i)$ and $\pi$ is birational to an étale cover, the claim is clear.

Then, by Lemma 4.1.5 there is an ideal $\mathscr{I}_i$ of $X$ such that

$$a_{Y*}(f_*\mathscr{O}_X(mK_X) \otimes \mathscr{I}_i \otimes \mathcal{L}_i^{\otimes -(m-2)})$$

is a sheaf satisfying I.T. with index 0, and $f_*\mathscr{O}_X(mK_X) \otimes \mathscr{I}_i \otimes \mathcal{L}_i^{\otimes -(m-2)}$ has rank $\alpha_m(X_y)$ (see for example [37, Proof of Lem. 3.9]). We conclude as before, that for $i \neq 0$, $f_*(\mathscr{O}_X(mK_X - (m-3)f^*\mathcal{L}_i))$ is globally generated over an open dense subset of $Y$. $\square$

## 4.4 General Type Case

In this section $X$ will be a variety of general type.

**Theorem 4.4.1.** *Let $X$ be a smooth projective variety, of maximal Albanese dimension and general type. Then, the linear system $|\omega_X^{\otimes 3} \otimes \alpha|$ induces a birational map, for any $[\alpha] \in \mathrm{Pic}^0(X)$.*

*Proof.* We reason by induction on the dimension of $X$, that we will denote by $n$. Note that for $n = 1$ the result is well known. So we assume that for any $[\alpha_Y] \in \mathrm{Pic}^0(Y)$, $|3K_Y + \alpha_Y|$ induces birational map for any smooth projective variety $Y$ of maximal Albanese dimension, general type, and $\dim Y \leq n - 1$.

Observe that $\mathcal{S}_X$ is empty if and only if $a_{X*}\omega_X$ is almost $M$-regular. Since $X$ is of general type, if $\mathcal{S}_X$ is empty, then $\chi(X, \omega_X) > 0$ by Proposition 4.1.4. Since Chen–Hacon have proved that $|3K_X|$ induces a birational map in this situation from now on, we will assume $\mathcal{S}_X$ is not empty.

As in the last section, we are in Setting 4.3.2. So, let $k$ be the maximal number of $\mathcal{S}_X$ and let $[\beta] + \widehat{B}$ be an irreducible component of $V^k_{a_X}(X, \omega_X)$. Note that we can choose $[\beta]$ such that $\langle[\beta]\rangle \cap \widehat{B} = 0$. Taking the étale cover induced by $G := \langle[\beta]\rangle$, we get the same commutative diagram as in the proof of Lemma 4.3.5 (see (4.3.1))

$$
\begin{array}{ccccc}
 & & \overset{b_{\widetilde{X}}}{\frown\!\!\!\longrightarrow} & & \\
\widetilde{X} & \xrightarrow{\pi} & X & \xrightarrow{a_X} & A_X \\
\big\downarrow{\scriptstyle \widetilde{f}} & & \big\downarrow{\scriptstyle f} & & \big\downarrow{\scriptstyle \mathrm{pr}} \\
\widehat{Y} & \xrightarrow{\widehat{\pi}} & Y & \xrightarrow{a_Y} & B, \\
 & & \underset{b_{\widehat{Y}}}{\smile\!\!\!\longrightarrow} & &
\end{array}
\tag{4.4.1}
$$

where $\widetilde{f}$ and $\widetilde{g}$ are obtained by taking respectively the Stein factorization of $\mathrm{pr} \circ b_{\widetilde{X}}$ and $\mathrm{pr}_Z \circ b_{\widetilde{X}}$ and after modifications, we may assume $\widehat{\pi}$ is a $G$-cover of smooth varieties.

Let $y \in Y$ be a general point and denote by $X_y$ a general fiber of $f$.

By Lemma 4.3.5, the restriction map

$$
H^0(X, \mathscr{O}_X(3K_X) \otimes \alpha) \to H^0(X_y, \mathscr{O}_{X_y}(3K_{X_y}) \otimes \alpha) \tag{4.4.2}
$$

is surjective for any $[\alpha] \in \mathrm{Pic}^0(X)$ and, by induction hypothesis,

$$
|3K_{X_y} + \alpha_{|X_y}|
$$

induces a birational map.

We have also produced interesting line bundles on $Y$ in Lemma 4.3.3 and Lemma 4.3.4. Let $\mathcal{H}$ be the line bundle on $Y$ constructed in Lemma 4.3.4. According to Lemma 4.2.1, in order to conclude the proof of the theorem, we just need to prove the following claim.

**Claim ‡.** For every $[\alpha] \in \widehat{A}_X$ and every $[\alpha'] \in \widehat{B}$, the line bundle

$$3K_X + \alpha - f^*(\mathcal{H} + \alpha')$$

has a non-trivial section.

Let be $\mathscr{J} := \mathscr{J}(||2K_X - f^*\mathcal{H} + \frac{1}{N}f^*H||)$, where $N$ is an integer large enough and $H$ is an ample divisor on $Y$. For any $[\alpha] \in \widehat{A}_X$, we define

$$\mathscr{F}_\alpha := f_*\big(\mathscr{O}_X(3K_X - f^*\mathcal{H}) \otimes \mathscr{J} \otimes \alpha)\big).$$

Observe that to conclude the proof of the claim it is enough to see that $V^0_{a_Y}(Y, \mathscr{F}_\alpha) = \widehat{B}$.

For any ample divisor $H'$ on $Y$, we have that

$$(\dagger)\ \ H^i(X, \mathscr{O}_X(3K_X - f^*\mathcal{H}) \otimes \mathscr{J} \otimes \alpha \otimes f^*\mathscr{O}_Y(H')) = 0,$$

for any $i > 0$. We postpone the proof of (†) to the end of the proof of this theorem. From (†) we deduce that

$$R^i f_*\big(\mathscr{O}_X(3K_X - f^*\mathcal{H}) \otimes \mathscr{J} \otimes \alpha)\big) = 0,$$

for any $i > 0$ (see Lemma 4.1.1). Therefore,

$$\chi(Y, \mathscr{F}_\alpha) = \chi(X, \mathscr{O}_X(3K_X - f^*\mathcal{H}) \otimes \mathscr{J} \otimes \alpha)$$

is constant for $P \in \widehat{A}_X$.

By Lemma 4.3.3 and Lemma 4.3.4, there exist integers $i$ and $j$ and effective divisors $D_1 \in |K_X + i\beta - f^*\mathcal{L}|$ and $D_2 \in |2K_X + j\beta - f^*\mathcal{H}|$. Let $m = i + j$ and write $D = D_1 + D_2 \in |3K_X + m\beta - f^*\mathcal{H} - f^*\mathcal{L}|$, i.e.

$$H^0(X, \mathscr{O}_X(3K_X + m\beta - f^*\mathcal{H} - f^*\mathcal{L})) = H^0(X, \mathscr{O}_X(D)) \neq 0.$$

Since

$$\begin{aligned}
\mathscr{J} &= \mathscr{J}(||2K_X - f^*\mathcal{H} + \frac{1}{N}f^*H||) & \\
&\supset \mathscr{J}(||2K_X - f^*\mathcal{H})||) & H \text{ is ample on } Y \\
&= \mathscr{J}(||2K_X + j\beta - f^*\mathcal{H}||) & [\beta] \text{ is torsion} \\
&\supset \mathscr{O}_X(-D_2) & \text{by Lemma 3.1.3 ,}
\end{aligned}$$

we have

$$\begin{aligned}
H^0(Y, \mathscr{F}_{m\beta} \otimes \mathcal{L}^{-1}) &= H^0(X, \mathscr{O}_X(3K_X + m\beta - f^*\mathcal{H} - f^*\mathcal{L}) \otimes \mathscr{J}) \\
&\supset H^0(X, \mathscr{O}_X(D_1)) \neq 0. \qquad (4.4.3)
\end{aligned}$$

Therefore, since $V^0_{a_Y}(Y, \mathcal{L}) = \widehat{B}$, we have $h^0(Y, \mathscr{F}_{m\beta} \otimes \alpha') > 0$ for all $[\alpha'] \in \widehat{B}$.

On the other hand, we see by Lemma 4.1.8 that $\mathscr{F}_\alpha$ is a $GV$-sheaf for any $[\alpha] \in \widehat{A}_X$. Therefore, for $[\alpha'] \in \widehat{B}$ general and any $[\alpha] \in \widehat{A}_X$,

$$h^0(Y, \mathscr{F}_\alpha \otimes \alpha') = \chi(B, \mathscr{F}_\alpha) = \chi(B, \mathscr{F}_{m\beta}) = h^0(B, \mathscr{F}_{m\beta} \otimes \alpha') > 0.$$

Hence, by semicontinuity, for any $\alpha \in \widehat{A}_X$, $V^0_{a_Y}(\mathscr{F}_\alpha) = \widehat{B}$.

*Proof of* (†) Notice that $2K_{\widetilde{X}/\widehat{Y}} \preceq \pi^*(2K_X - f^*\mathcal{H})$ and $K_{\widetilde{X}/\widehat{Y}} + \frac{1}{N}\widetilde{f}^*\widehat{\pi}^*H$ is a big $\mathbb{Q}$-divisor on $\widetilde{X}$. Hence, $2K_X - f^*\mathcal{H} + \frac{1}{N}f^*H$ is a big $\mathbb{Q}$-divisor on $X$. So (†) is a consequence of Kawamata-Viehweg vanishing theorem (Theorem 4.1.2). □

## 4.5 Iitaka Fibrations

In this section $X$ will not necessarily be a variety of general type.

**Theorem 4.5.1.** *Let $X$ be a smooth projective variety, of maximal Albanese dimension. Then, the linear system $|4K_X + \alpha|$ induces a model of the Iitaka fibration of $X$, for any $[\alpha] \in V^0_{a_X}(X, \omega_X^{\otimes 2})$.*

Before starting the proof the Theorem 4.5.1, which is parallel to the proof of Theorem 4.4.1, let us fix the notation.

**Setting 4.5.2.** Consider the following diagram:

$$\begin{array}{ccc} X & \xrightarrow{a_X} & A_X \\ \downarrow{\scriptstyle g} & & \downarrow{\scriptstyle \mathrm{pr}_Z} \\ Z & \xrightarrow{a_Z} & A_Z \end{array}$$

where $g : X \to Z$ which is a model of the Iitaka fibration of $X$ and $Z$ is smooth. Let $K$ be the kernel of $\mathrm{pr}_Z$. We denote by $X_z$ a general fiber of $g$, which is birational to its Albanese variety $\widetilde{K}$, and the natural map $\widetilde{K} \to K$ is an isogeny. We know that $\mathrm{pr}_Z^* \widehat{A}_Z$ is an irreducible component of

$$\mathscr{K} := \ker(\widehat{A}_X \to \mathrm{Pic}^0(X_z))$$

and denote by $\mathcal{Q} := \mathscr{K} / \mathrm{pr}_Z^* \widehat{A}_Z$. Observe that $\mathscr{Q}$ can be also identified with $\ker(\widehat{K} \to \widehat{\widetilde{K}})$.

*Remark* 4.5.1. The group $\mathcal{Q}$ is often non-trivial and this is exactly the reason why the tricanonical map can not always induce the Iitaka fibration. In some specific cases, given information about $\mathcal{Q}$, we could prove that the tricanonical map or some twisted tricanonical map (the maps induced by $|3K_X + \alpha|$ for some $[\alpha] \in \mathrm{Pic}^0(X)$) would induce the Iitaka fibration (see Remark 4.5.2).

Nevertheless we will construct a variety of maximal Albanese dimension (see Example 4), where NONE of the twisted tricanonical maps is birationally equivalent to the Iitaka fibration.

Before proving the theorem, we start with an easy observation.

**Lemma 4.5.3.** *The kernel $\mathscr{K}$ defined in Setting 4.5.2, satisfies*

$$\mathscr{K} = V^0_{a_X}(X, \omega_X^{\otimes m}) \qquad \textit{for all } m \geq 2.$$

*Proof.* It is clear that $V^0_{a_X}(X, \omega_X^{\otimes m}) \subseteq \mathscr{K}$. If $[\alpha] \in \mathscr{K}$, then $g_*(\omega_X^{\otimes m} \otimes \alpha)$ is a nontrivial torsion-free sheaf. By Lemma 4.1.5, $g_*(\omega_X^{\otimes m} \otimes \mathscr{J}(\|\omega_X^{\otimes(m-1)}\|) \otimes \alpha)$ is a sheaf satisfying I.T. with index 0 for any $m \geq 2$. Hence, we conclude since $0 < h^0(Z, g_*(\omega_X^{\otimes m} \otimes \mathscr{J}(\|\omega_X^{\otimes(m-1)}\|) \otimes \alpha)) \leq h^0(X, \omega_X^{\otimes m} \otimes \alpha)$. □

*Proof of Theorem 4.5.1.* We will prove the theorem by induction on the dimension of $X$. We suppose the statement is true in dimension $\leq n-1$ and assume $\dim X = n$. If $X$ is of general type, then we are back to Theorem 4.4.1. Hence we can assume $\kappa(X) = \dim Z = n - l$, for some number $l > 0$. In particular, $\mathcal{S}_X$ is not empty.

Hence, we are in Setting 4.3.2. Let $k$ be the maximal number of $\mathcal{S}_X$ and let $[\beta] + \widehat{B}$ be an irreducible component of $V^k_{a_X}(X, \omega_X)$. Since $R^k f_*(\omega_X \otimes \beta \otimes \alpha) \neq 0$, for all $[\alpha] \in \widehat{B}$, we observe that $\widehat{B} \hookrightarrow \mathrm{pr}_Z^* \widehat{A}_Z$. Hence, Setting 4.3.2 and 4.5.2 combine in the following commutative diagram

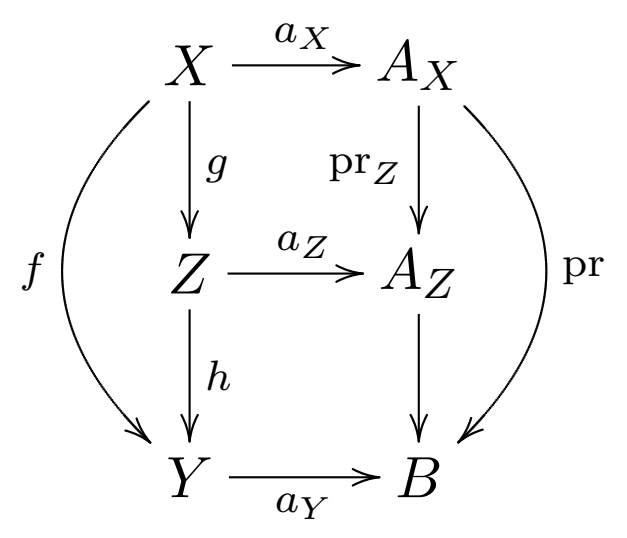

Let $y \in Y$ be a general point and denote by $X_y$ and $Z_y$ general fibers of $f$ and $h$. By Easy Addition Formula (e.g. [35, Thm. 10.4]), $\dim Y + \kappa(X_y) \geq \kappa(X) = \dim Z$. Hence $\kappa(X_y) \geq \dim Z_y$ and then $g_{|X_y} : X_y \to Z_y$ is the Iitaka fibration of $X_y$.

By Lemma 4.3.5 and 4.5.3, the restriction map

$$H^0(X, \mathscr{O}_X(4K_X + \alpha)) \to H^0(X_y, \mathscr{O}_{X_y}(4K_{X_y} + \alpha_{|X_y}))$$

is surjective, for any $[\alpha] \in \mathscr{K}$. Notice that $[\alpha_{|X_y}] \in V^0_{a_{X_y}}(X_y, \omega_{X_y}^{\otimes 2})$, so by induction hypothesis,

$$|4K_{X_y} + \alpha_{|X_y}|$$

induces the Iitaka fibration $g_{|X_y} : X_y \to Z_y$.

Let $\mathcal{H}$ be the line bundle on $Y$ constructed in Lemma 4.3.4. Then, by Lemma 4.2.1, we just need to prove the following claim to finish the proof of the theorem.

**Claim.** For every $[\alpha] \in \mathscr{K}$ and every $[\alpha'] \in \widehat{B}$,

$$4K_X + \alpha - f^*(\mathcal{H} + \alpha')$$

has a non-trivial section.

Let be $\mathscr{J} := \mathscr{J}(||3K_X - f^*\mathcal{H} + \frac{1}{N} f^*H||)$, where $N$ is an integer large enough and $H$ is an ample divisor on $Y$. For any $\alpha \in \mathscr{K}$, we define

$$\mathscr{G}_\alpha := g_*\big(\mathscr{O}_X(4K_X + \alpha - f^*\mathcal{H}) \otimes \mathscr{J}\big).$$

Observe that to conclude the proof of the claim it is enough to see that $V^0_{a_Y}(Y, h_*\mathscr{G}_\alpha) = \widehat{B}$.

By Lemma 4.1.5, we have

$$H^i(Z, \mathscr{G}_\alpha \otimes \beta'' \otimes h^*H') = 0,$$

for any $i \geq 1$, any ample divisor $H'$ on $Y$, and any $[\beta''] \in \widehat{A}_Z$. Hence, by Lemma 4.1.1

$$R^i h_*(\mathscr{G}_\alpha \otimes \beta'') = 0,$$

for any $i > 0$. Therefore,

$$\chi(Y, h_*(\mathscr{G}_\alpha \otimes \beta'')) = \chi(Z, \mathscr{G}_\alpha \otimes \beta'') \tag{4.5.1}$$

is constant for $\beta'' \in \widehat{A}_Z$.

By Lemma 4.3.3 and Lemma 4.3.4, we know there exists $m \in \mathbf{Z}$ such that

$$H^0(X, \mathscr{O}_X(3K_X + m\beta - f^*\mathcal{H} - f^*\mathcal{L})) \neq 0.$$

Observe that $[\alpha - m\beta]$ is not necessarily in $\mathrm{pr}_Z^* \widehat{A}_Z$. But, since $[\alpha - m\beta] \in \mathscr{K}$, we have that $a_{Z*}g_*\mathscr{O}_X(K_X + \alpha - m\beta)$ is a non-trivial $GV$-sheaf. In particular, $V^0_{a_Z}(Z, g_*\mathscr{O}_X(K_X + \alpha - m\beta)) \neq \emptyset$. Hence there exists $[\beta_0] \in \mathrm{pr}_Z^* \widehat{A}_Z$ such that $[\alpha - m\beta + \beta_0] \in V^0_{a_X}(X, \omega_X)$.

Therefore

$$\begin{aligned} & 4K_X + \alpha + \beta_0 - f^*\mathcal{H} - f^*\mathcal{L} \\ = \; & (K_X + \alpha - m\beta + \beta_0) + (3K_X + m\beta - f^*\mathcal{H} - f^*\mathcal{L}) \end{aligned}$$

is the sum of two effective divisors. By the same argument as in (4.4.3),

$$H^0(Z, \mathscr{G}_\alpha \otimes h^*\mathcal{L}^{-1} \otimes \beta_0) \neq 0.$$

We know that $V^0_{a_Y}(Y,\mathcal{L})=\widehat{B}$ (see Lemma 4.3.3) and $h_*(\mathscr{G}_\alpha\otimes\beta')$ is a $GV$-sheaf for any $\alpha\in\mathscr{K}$ and any $\beta'\in\widehat{A}_Z$ (see Lemma 4.1.8). Hence, for $\alpha'\in\widehat{B}$ general,

$$\begin{aligned} h^0(Y,h_*\mathscr{G}_\alpha\otimes\alpha') &= \chi(Y,h_*\mathscr{G}_\alpha)=\chi(Y,h_*(\mathscr{G}_\alpha\otimes\beta_0)) && \text{by (4.5.1)}\\ &= h^0(Y,h_*(\mathscr{G}_\alpha\otimes\beta_0)\otimes\alpha')>0. \end{aligned}$$

By semicontinuity, $V^0_{a_Y}(h_*\mathscr{G}_\alpha)=\widehat{B}$ for any $[\alpha]\in\mathscr{K}$. □

Now, we can make more precise Remark 4.5.1.

*Remark* 4.5.2. In the previous proof, observe that if $[\alpha-m\beta]$ lies in $\mathrm{pr}_Z^*\widehat{A}_Z$, then $3K_X+\alpha+\beta_0-f^*\mathcal{H}-f^*\mathcal{L}$ is effective for some $\beta_0\in\mathrm{pr}^*\widehat{A}_Z$. So, we could have improved the result to the tricanonical map (assuming the induction hypothesis). In particular, if $\mathscr{Q}:=\mathscr{K}/\mathrm{pr}_Z^*\widehat{A}_Z$ is trivial for $X$ and the successive fibers of the induction process, then the tricanonical map twisted by an element in $\mathscr{K}$ induces the Iitaka fibration.

Moreover, if for some $[\alpha]\in\mathscr{K}$, $[\alpha]+\mathrm{pr}_Z^*\widehat{A}_Z$ is an irreducible component of $V^0_{a_X}(X,\omega_X)$, then we can again prove that the tricanonical map twisted by an element in $\mathscr{K}$ induces the Iitaka fibration. This shows that varieties of maximal Albanese dimension, where none of the twisted tricanonical map is birational equivalent to the Iitaka fibration, are closely related to varieties of maximal Albanese dimension, of general type with vanishing holomorphic Euler characteristic.

We finish with an example of maximal Albanese dimension, whose tricanonical map does not induce the Iitaka fibration. This example is based on the famous Ein-Lazarsfeld threefold, which is constructed in [16, Ex. 1.13] and further investigated in [7].

**Example 4.** We take three bielliptic curves $C_i$ of genus 2, $i=1,2,3$. Let $\rho_i:C_i\to E_i$ be the double cover over an elliptic curve $E_i$ and denote by $\tau_i$ the involution of fibers of $\rho_i$. We write

$$\rho_*\mathscr{O}_{C_i}=\mathscr{O}_{E_i}\oplus\mathcal{L}_i^{-1},$$

where $\mathcal{L}_i$ is a line bundle on $E_i$ of degree 1.

Let $Y$ be the threefold $(C_1\times C_2\times C_3)/(\tau_1,\tau_2,\tau_3)$, which has only rational singularities. We know that $a_Y:Y\to E_1\times E_2\times E_3$ is a $(\mathbf{Z}/2\mathbf{Z}\times\mathbf{Z}/2\mathbf{Z})$-cover.

We then take an abelian variety $A$ and a $(\mathbf{Z}/2\mathbf{Z}\times\mathbf{Z}/2\mathbf{Z})$-étale cover $\widetilde{A}\to A$. Set $\{[\mathscr{O}_A,][\alpha_1],[\alpha_2],[\alpha_3]\}$ to be the kernel $\widehat{A}\to\widehat{\widetilde{A}}$.

Denote $H=(\mathbf{Z}/2\mathbf{Z}\times\mathbf{Z}/2\mathbf{Z})$ and let $X'$ be the variety $(Y\times\widetilde{A})/H$, where $H$ acts diagonally on $Y\times\widetilde{A}$. Notice that $X'$ has only rational singularities and let $X$ be a resolution of singularities of $X'$. The Albanese morphism

$$a_X:X\to E_1\times E_2\times E_3\times A$$

is birationally a $(\mathbf{Z}/2\mathbf{Z}\times\mathbf{Z}/2\mathbf{Z})$-cover.

After permutation of $\{[\alpha_i], i = 1, 2, 3\}$, we have

$$\begin{aligned} a_{X*}\omega_X^{\otimes 3} &\simeq \left(\mathcal{L}_1^{\otimes 2} \boxtimes \mathcal{L}_2^{\otimes 2} \boxtimes \mathcal{L}_3^{\otimes 2} \boxtimes \mathscr{O}_A\right) \oplus \left(\mathcal{L}_1^{\otimes 3} \boxtimes \mathcal{L}_2^{\otimes 3} \boxtimes \mathcal{L}_3^{\otimes 2} \boxtimes \alpha_1\right) \\ &\oplus \left(\mathcal{L}_1^{\otimes 3} \boxtimes \mathcal{L}_2^{\otimes 2} \boxtimes \mathcal{L}_3^{\otimes 3} \boxtimes \alpha_2\right) \oplus \left(\mathcal{L}_1^{\otimes 2} \boxtimes \mathcal{L}_2^{\otimes 3} \boxtimes \mathcal{L}_3^{\otimes 3} \boxtimes \alpha_3\right) \end{aligned}$$

It is easy to check that for any $[\alpha] \in \mathrm{Pic}^0(X)$, the linear series $|3K_X + \alpha|$ can not induce the Iitaka fibration $X \to E_1 \times E_2 \times E_3$.

Using, the notation of Setting 4.5.2, observe that

$$\mathscr{K} = V^0_{a_X}(X, \omega_X^{\otimes 2}) = \bigcup_{[\beta] \in \{[\mathscr{O}_A][, \alpha_1], [\alpha_2], [\alpha_3]\}} E_1 \times E_2 \times E_3 \times \{[\beta]\}$$

and $\mathcal{Q} = \mathbb{Z}/2\mathbb{Z} \times \mathbb{Z}/2\mathbb{Z}$. Indeed, $\mathcal{Q}$ can be identified with $\{[\mathscr{O}_A][, \alpha_1], [\alpha_2], [\alpha_3]\}$.

# Part III

# Theta Divisors in the Classification of Varieties

# CHAPTER 5

# BIRATIONALITY OF THE A LBANESE MAP AND OTHER RESULTS

Smooth models of theta divisors are, after abelian varieties, the simplest and most important among irregular varieties. An interesting research started by Ein–Lazarsfeld ([16]) is to characterize irreducible theta divisors by their birational invariants, in a similar way to what Kawamata ([40]) did for abelian varieties. In this setting we recall the following result of Hacon–Pardini ([32]).

**Theorem.** *Let $X$ be a smooth complex variety of dimension $n$, let $A$ be an abelian variety of dimension $n+1$, and let $a : X \to A$ a generically finite morphism. Assume that*

*a)* $f^* : H^0(X, \Omega^n_X) \to H^0(A, \Omega^n_A)$ *is an isomorphism;*

*b)* $h^1(X, \omega_X \otimes a^*\alpha) = 1$ *for every* $[\alpha] \in \widehat{A} \backslash \{[\mathscr{O}_A]\}$ *and all* $i > 0$. *Then $f$ is birational into its image $f(X)$, $A$ is principally polarized and $f(X)$ is a theta divisor in $A$.*

Not long ago Barja–Lahoz–Naranjo–Pareschi, and independently Lazarsfeld–Popa, provided this new characterization.

**Theorem** ([2],[52])**.** *Let $X$ be a smooth complex variety that admits $a : X \to A$ a generically finite morphism into an abelian variety such that*

*a)* $\chi(\omega_X) = 1$;

*b)* $\dim A > \dim X$;

*c) for every* $0 < i \leq \dim X$, $V^i_a(X, \omega_X) = \{\mathscr{O}_A\}$.

*Then $A$ is principally polarized,* $\dim A = n+1$, *$a$ has generic degree 1 and $X$ is a smooth model of a principal polarization in $A$.*

This was more recently refined in [59] (see Theorem 5.3.1).

It is tempting to manage to give similar characterizations for *products* of smooth irreducible theta divisors. To this aim Pareschi conjectured the following:

**Conjecture** (Pareschi). *If $X$ is a variety of maximal Albanese dimension, whose Albanese map is not fibered in subtori, and $\chi(\omega_X) = 1$, then it is birational to a product of irreducible $\Theta$-divisors in principally polarized abelian varieties.*

It is well known that subvarieties of abelian varieties either are fibered in subtori or (their resolution of singularities) are of general type. Certainly products of theta divisors belong to the latter of these two categories. If the conjecture were proved to be true, we would have a partial converse of this fact.

The above conjecture is known to be true for surfaces. In fact Beauville proved in [3] that the irregularity of surfaces with $\chi = 1$ is at most 4 and surfaces with irregularity 4 and Euler characteristic equal to one are products of curves of geneus 2 (and therefore products of irreducible theta divisors). Surface $S$ with $\chi(\omega_S) = 1$ and irregularity $q = 3$ were studied in [33] and [66], where was proved that the moduli space of these irregular surfaces has just two connected components

(a) $\mathcal{M}_6 := \{[X] \mid K_X^2 = 6\}$; the elements of this component are symmetric products of curves of genus 3 (and therefore are birational to $\Theta$-divisors in principally polarized abelian variety).

(b) $\mathcal{M}_8 := \{[X] \mid K_X^2 = 8\}$; in this case $X$ is a quotient of a product $C_2 \times C_3$ with $C_i$ a curve of genus $i$, where $C_2$ has an ellpitic involution $\sigma_2$ and $C_3$ as a free involution $\sigma_3$, and $X$ is the quotient by the diagonal action.

Notice that the Albanese images of surfaces of type (b), contrary to those in (a), are fibered in subtori. Therefore the surfaces of maximal Albanese dimension with Euler characteristic equal to one and Albanese image not ruled by subtori are birational to a product of theta divisors.

Other examples of varieties that satisfy Pareschi's conjecture are provided by *highly irregular varieties*, where the meaning of "highly irregular" is to be taken with respect of the following statement of Hacon–Pardini which generalize in higer dimension a result of Beauville ([3]):

**Theorem** (Hacon–Pardini [34]). *The irregularity of smooth complex variety $X$ of dimension $n$, $\chi(X) = 1$ and of maximal Albanese dimension is less or equal $2n$ and equality occurs if and only if it is birational to a product of curves of genus 2.*

Hence varieties $X$ with $\chi(X) = 1$ and maximal irregularity satisfy Pareschi's conjecture.

In this Chapter we study smooth projective varieties with $\chi(\omega_X) = 1$, whose Albanese image is not fiberd in tori, in order to challenge Pareschi's conjecture. We present some partial results. First of all, in the same setting as in the conjecture above, we were able to prove the following theorem that could be regarded as a further evidence of trueness of the conjecture, since the Albanese map of a desingularisation of a product of theta divisor is obviously birational.

**Theorem 5.A.** *The Albanese map of $X$ has degree 1.*

Afterwards we used this fact to further refine the known cohomological characterizations of (smooth models of) theta divisors recalled at the beginning of this introduction. In particular we prove

**Theorem 5.B.** *Let $X$ as above and suppose that its Albanese image is normal. Consider the generic vanishing loci*

$$V^i_{\mathrm{alb}_X}(\omega_X) := \{[\alpha] \in \mathrm{Pic}^0(X) \mid h^i(\omega_X \otimes \alpha) \neq 0\}.$$

*If $[\mathscr{O}_X]$ is an isolated point of $V^i_{\mathrm{alb}_X}(\omega_X)$ for every $i > 0$ then $X$ is birational to an irreducible theta divisors in an principally polarized abelian variety.*

We remark that the hypothesis about the normality of the Albanese image of $X$ is not very restrictive thank to a theorem of Ein–Lazarsfeld ([16]) that states that irreducible theta divisors are always normal with at worst rational singularities. Furthermore we believe that this assumption could eventually be removed.

We would like to remark that Theorem 5.B constitutes a further evidence of the validity of Pareschi's conjecture. In fact if this were to be proved, then theta divisors would be the sole irregular varieties whose generic vanishing loci have an isolated point.

We hope that Theorem 5.B (or a refinements of its), thank to the better understanding of theta divisors that it allows, would help us to succed in the proof of Pareschi's conjecture.

This Chapter is organized in the following manner. In section 5.1 we present some background material and some preliminary lemmas. Section 5.2 is dedicated to the proof of Theorem 5.A. To this aim we present some results on the Fourier-Mukai transform of the structure sheaf $\mathscr{O}_X$ and we expose a birational criterion similar to [32, Theorem 3.1]. In the last section we prove Theorem 5.B.

## 5.1 Preliminaries

In this section we recall some known results which will be needed afterwards. We begin with two Theorems of Hacon–Pardini tha generalize the work of Kollar and Green–Lazarsfeld.

**Theorem 5.1.1** ([34, Theorem 2.1]). *Let $X$, $Y$ and $Z$ be projective varieties, with $X$ smooth. Consider $f : X \to Y$ and $g : Y \to Z$ two surjective morphisms and take $[\tau] \in \mathrm{Pic}^0(X)$ a torsion point. Then:*

*(i)* $R^i g_* R^j f_*(\omega_X \otimes \tau)$ *is torsion free for every* $i, j \geq 0$*;*

*(ii)* $R^i g_* R^j f_*(\omega_X \otimes \tau) = 0$ *for every* $j \geq 0$, $i > \dim ; Y - \dim Z$*;*

*(iii)* $R^k (g \circ f)_*(\omega_X \otimes \tau) \simeq \bigoplus_{i=0}^k R^i g_* R^{k-i} f_*(\omega_X \otimes \tau)$.

**Theorem 5.1.2** (Generic Vanishing [34, Theorem 2.2]). *Let $f : X \to Y$ and $g : Y \to Z$ be morphisms of smooth projective varieties and let $a : Z \to A$ a morphism to an abelian variety. Take $[\tau] \in \mathrm{Pic}^0(X)$ a torsion point and set $\mathscr{F} := R^h g_* R^j f_*(\omega_X \otimes \tau)$. Then:*

*(i) set $k = \dim(Z) - \dim a(Z)$, $\mathrm{gv}_a(\mathscr{F}) \geq -k$. In particular if $Z$ is of maximal Albanese dimension then $\mathscr{F}$ is a GV-sheaf.;*

*(ii) every irreducible component of $V_a^i(\mathscr{F})$ is a translate of a subtorus of $\widehat{A}$ by a torsion point.*

In [34] Hacon and Pardini proved (a more general version the of the) following

**Lemma 5.1.3** ([34, Lemma 3.3]). *Let $X$ be a maximal Albanese dimesnion with $a : X \to A$ a generically finite morphism into an abelian variety. Given $W$ a component of $V_a^k(\omega_X \otimes \tau)$, denote by $g : X \to \mathrm{Pic}^0(W)$ the induced map. Then*

$$\dim g(X) \leq \dim X - k.$$

We include the short proof for the reader's benefit.

*Proof.* We can write $W = [\beta] + \widehat{B}$ with $[\beta] \in \widehat{A}$ a torsion point and $\widehat{B} \subseteq \widehat{A}$ a subtorus. Let $B := \mathrm{Pic}^0(\widehat{B}) = \mathrm{Pic}^0(W)$. Then by the Theorem 5.1.2 above, for the generic $\alpha \in \widehat{B}$ and for every $j \geq 0$ and $s > 0$ one has:

$$H^s(B, R^j g_*(\omega_X \otimes a^*\beta) \otimes \alpha) = 0.$$

Therefore the Leray spectral sequence degenrates and, thank to Theorem 5.1.1, we have the isomorphism

$$H^k(X, \omega_X \otimes a^*\beta \otimes g^*\alpha) = H^0(B, R^i g_*(\omega_X \otimes a^*\beta) \otimes \alpha).$$

Now assume by contradiction that the relative dimension of $g$ is strictly smaller than $k$, by 5.1.1 $R^k g_*(\omega_X \otimes a^*\beta) = 0$ and $H^k(\omega_X \otimes a^*(\beta \otimes \alpha)) = 0$, contradicting the hypothesis that $W \subseteq V_a^k(\omega_X)$. $\square$

We can evince some more complete information about the dimension of the generic fiber of $g$:

**Lemma 5.1.4.** *In the notation of the above Lemma, suppose furthermore that the generic vanishing index of $X$ (cfr Definition 1.2.2), $gv_a(\omega_X)$, is equal to $i$ and let $W$ be a component of $V^k_a(\omega_X)$ of codimension $k+i$. Then*

$$\dim g(X) = \dim X - k.$$

*Proof.* We will show that the general fiber of $g$, $G$, has dimension $k$. The argument we present is due to Jiang ([38]).

As before, we can write $W = [\beta] + \widehat{B}$ with $[\beta] \in \widehat{A}$ a torsion point and $\widehat{B} \subseteq \widehat{A}$ a subtorus. Let $B := \mathrm{Pic}^0(\widehat{B}) = \mathrm{Pic}^0(W)$. Now consider $\pi : A \to B$ the dual map to the inclusion $\widehat{B} \hookrightarrow \widehat{A}$. Thus we are in the following situation:

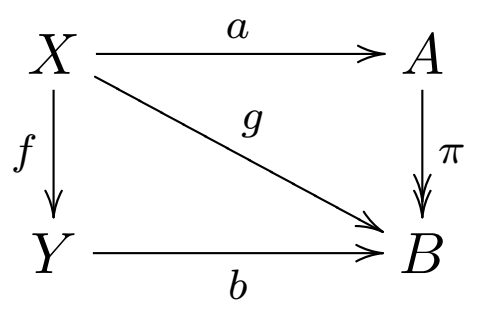

where $f$ is (a smooth model of) the Stein factorization of $g$ and $b$ is generically finite. Since $\dim Y = \dim g(X)$, by the above Lemma, the generic fiber $F$ of $f$ has dimension greater or equal $k$. We claim that equality holds, the statement would follow at once. Suppose that equality did not hold and take $j \geq 0$ and $Z$ an irreducible component of $V^j_b(\omega_Y)$. Then $\pi^*Z \subseteq V^{j+\dim F}_a(\omega_X)$ and we would have, for any $j \geq 0$

$$k + j + i < \dim F + j + i \leq \mathrm{codim}_{\widehat{A}}\pi^*Z = \mathrm{codim}_{\widehat{B}}Z + k + i.$$

Therefore, for any $j \geq 0$, $\mathrm{codim}_{\widehat{B}}V^j_b(\omega_Y) \geq j+1$. In particular we have that $V^0_b(\omega_Y)$ is a *proper* subset of $\widehat{B}$ and hence $\chi(\omega_Y) = 0$, but $gv_b(\omega_Y) \geq 1$ contradicting Proposition 1.2.10. □

## 5.2 When the Albanese Map Has Degree 1

Let $X$ a smooth projective variety over the complex numbers of dimension $n$. From now on we will work under the following hypothesis:

**Hypothesis 5.2.1.** The variety $X$ is of maximal Albanese dimension and

(i) $\chi(\omega_X) = 1$

(ii) $q(X) := h^1(X, \mathcal{O}_X) > \dim X$.

Furthermore it will be given $a : X \to A$ a generically finite map to an abelian variety of dimension $q$ such that

(iii) the pullback $a^* : \widehat{A} \to \mathrm{Pic}^0(X)$ is an embedding and

(iv) the image of $a$ is not fibered in subtori of $A$.

*Remark* 5.2.1. We recall that condition (iv) tells us that $\mathrm{gv}_a(\omega_X) \geq 1$. Infact, suppose that $\mathrm{gv}_a(\omega_X) = 0$ and take $[\beta] + \widehat{T} \subseteq V_a^k(\omega_X)$ a component of codimension $k$. Consider, as before, the following commutative diagram

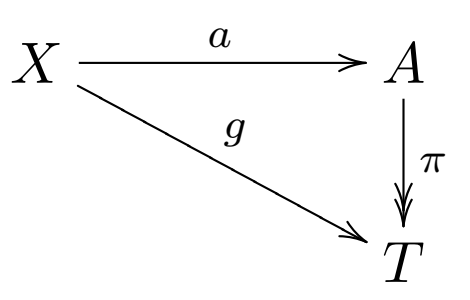

where $T := \mathrm{Pic}^0(\widehat{T})$, and $\pi : A \to T$ is the map dual to the inclusion $\widehat{T} \subseteq \widehat{A}$. By Lemma 5.1.3 the general fiber of $g$ has dimension greater or equal $k$. This means that, via $a_{|G}$, $G$ is mapped surjectively to the $k$-dimensional general fiber of $\pi$, that is a translate of the subtorus $K := \mathrm{Ker}(\pi) \subseteq A$. Hence the general fiber of $\pi_{|a(X)} : a(X) \to \widehat{T}$ is still isomorphic to $K$, implying that $a(X)$ is ruled by subtori. Therefore $gv_a(\omega_X) \geq 1$. Now we combine (i) with proposition 1.2.10 and get that $gv_a(\omega_X) = 1$.

The main purpose of this section is to prove the following statement.

**Theorem 5.2.2.** *Let $(X, a)$ be a pair that satisfies Hypothesis 5.2.1, then the map $a$ is birational onto its image.*

Observe as this would immediately imply Theorem 5.A in the introduction.

### 5.2.1 The Fourier-Mukai Transform of the Structure Sheaf

We will begin by providing information about the complex $\mathbf{R}S\mathbf{R}a_*\mathscr{O}_X$ under the assumption 5.2.1. By Remark 5.2.1 and the W.I.T crieterion (Theorem 1.2.4), we know that its cohomology is concentrated in degree $n = \dim X$, while Proposition 1.2.6 and Theorem 1.2.9 tell us that its only non vanishing cohomology sheaf $\widehat{\mathbf{R}a_*\mathscr{O}_X}$ is a torsion free sheaf of generic rank $\chi(\omega_X) = 1$, hence it is an ideal sheaf twisted by some line bundle:

$$\widehat{\mathbf{R}a_*\mathscr{O}_X} \simeq \mathscr{I}_Z \otimes L.$$

Observe that the codimension of $Z$ is greater than 1. In fact, suppose by contradiction that exists $W \subseteq Z$ a component of codimension 1, then $W$ is a component of the support of $\mathscr{E}xt^1(\mathscr{O}_Z \otimes L, \mathscr{O}_A)$. Now consider the short exact sequence

$$0 \to \mathscr{I}_Z \otimes L \to L \to \mathscr{O}_Z \otimes L \to 0 \tag{5.2.1}$$

Applying $\mathbf{R}\mathscr{H}om(-,\mathscr{O}_A)$ and taking cohomology we get

$$0 \to L^\vee \xrightarrow{\varphi} L^\vee \xrightarrow{\psi} \mathscr{E}xt^1(\mathscr{O}_Z \otimes L, \mathscr{O}_A) \to 0.$$

Now $\varphi$ is certainly an isomorphism, therefore $\psi$ is the zero map and $\mathscr{E}xt^1(\mathscr{O}_Z \otimes L, \mathscr{O}_A) = 0$ contradicting the inclusion

$$W \subseteq \operatorname{Supp} \mathscr{E}xt^1(\mathscr{O}_Z \otimes L, \mathscr{O}_A)$$

*Remark* 5.2.2. a)The subscheme $Z$ above is supported on the union of the generic vanishing loci $V^i_a(\omega_X)$ for $i \geq 1$. In fact let $W$ be a component of $Z$, for what observed above we can suppose that $\operatorname{codim}_{\widehat{A}} W = i+1$ for some $i > 1$. Then $W$ is included in the support of $\mathscr{E}xt^{i+1}(\mathscr{O}_Z \otimes L, \mathscr{O}_{\widehat{A}})$. Applying again $\mathbf{R}\mathscr{H}om(-,\mathscr{O}_{\widehat{A}})$ to (5.2.1) and taking cohomology we get that for every $j \geq 1$

$$\mathscr{E}xt^{j+1}(\mathscr{O}_Z \otimes L, \mathscr{O}_{\widehat{A}}) \simeq \mathscr{E}xt^j(\mathscr{I}_Z \otimes L, \mathscr{O}_{\widehat{A}}).$$

In particular $W \subseteq \operatorname{Supp} \mathscr{E}xt^i(\mathscr{I}_Z \otimes L, \mathscr{O}_{\widehat{A}})$. We apply Theorem 1.2.6(b) that tells us

$$\operatorname{Supp} \mathscr{E}xt^i(\mathscr{I}_Z \otimes L, \mathscr{O}_{\widehat{A}}) = \operatorname{Supp} \mathscr{E}xt^i(\widehat{\mathbf{R}a_*\mathscr{O}_X}, \mathscr{O}_{\widehat{A}}) \simeq R^i S(\mathbf{R}a_*\omega_X).$$

In particular we obtain that $W \subseteq \operatorname{Supp}(R^i\Phi_{\mathscr{P}_a}(\omega_X))$. By base change we have an inclusion $\operatorname{Supp}(R^i\Phi_{\mathscr{P}_a}(\omega_X)) \subseteq V^i_a(\omega_X)$. Therefore $W \subseteq V^i_a(\omega_X)$.

b) The sheaf $\mathscr{I}_Z \otimes L$ is a GV sheaf on $\widehat{A}$. Infact, applying Proposition 1.2.6 with $\mathscr{F} \simeq \omega_X$, we have the following isomorphisms of objects in $\mathbf{D}(A)$

$$\mathbf{R}\Delta(\mathscr{I}_Z \otimes L) \simeq (-1_{\widehat{A}})^*\mathbf{R}\Phi_{\mathscr{P}_a}(\omega_X) \simeq (-1_{\widehat{A}})^*\mathbf{R}S\mathbf{R}a_*\omega_X; \tag{5.2.2}$$

where the last equality is nothing else that equation (1.1.3). Now we apply to both sides the functor $\mathbf{R}\widehat{S}$ and we use Mukai's duality thoerem:

$$\mathbf{R}\widehat{S}\mathbf{R}\Delta(\mathscr{I}_Z \otimes L) \simeq \mathbf{R}\widehat{S}(-1_{\widehat{A}})^*\mathbf{R}S\mathbf{R}a_*\omega_X \simeq \mathbf{R}a_*\omega_X[-q].$$

Since $a$ is generically finite, Grauert-Riemenshneider implies that the latter is a sheaf in degree $q$. It follows from the W.I.T. crieterion (Theorem 1.2.4) that $\mathscr{I}_Z \otimes L$ is GV.

c) As in [59], it is easy to compute the Fourier-Mukai transform of $\mathscr{I}_Z \otimes L$ using Mukai's theorem and what already proven:

$$\mathbf{R}\widehat{S}(\mathscr{I}_Z \otimes L) \simeq \mathbf{R}\widehat{S}\mathbf{R}S(\mathbf{R}a_*\mathscr{O}_X)[n] \simeq (-1_A)^*\mathbf{R}a_*\mathscr{O}_X[n-q]. \tag{5.2.3}$$

**Proposition 5.2.3.** *The sheaf $L$ is ample.*

*Proof.* Recall that given $x \in A$ we denote by $P_x$ the topologically trivial line bundle on $A$ associated to $x$, $\mathscr{P}_{|\{x\}\times\widehat{A}}$. Let us suppose that $L$ is not ample. Hence, by some standard facts about line bundles on abelian varieties, there exists $T$ a subtorus of $A$ such that for every $x \notin T$, $h^0(\widehat{A}, L \otimes P_x) = 0$. In particular, if $x \notin T$, then also $-x \notin T$ and we have:

$$0 = H^0(\widehat{A}, \mathscr{I}_Z \otimes L \otimes P_{-x}) \simeq \mathrm{Hom}_{\mathbf{D}(\widehat{A})}(\mathscr{O}_{\widehat{A}}, \mathscr{I}_Z \otimes L \otimes P_{-x}) \simeq$$
$$\simeq \mathrm{Hom}_{\mathbf{D}(\widehat{A})}(\mathbf{R}\Delta_{\widehat{A}}(\mathscr{I}_Z \otimes L), P_{-x}).$$

Now we apply to both $P_{-x}$ and $\mathbf{R}\Delta_{\widehat{A}}(\mathscr{I}_Z \otimes L)$ the Fourier-Mukai transform. Using that it is an equivalence (and hence, in particular, it is fully faithful) we get:

$$0 = \mathrm{Hom}_{\mathbf{D}(\widehat{A})}(\mathbf{R}\widehat{S} \circ \mathbf{R}\Delta_{\widehat{A}}(\mathscr{I}_Z \otimes L), \mathbf{R}\widehat{S}(P_{-x})) \simeq$$
$$\simeq \mathrm{Hom}_{\mathbf{D}(\widehat{A})}((-1_A)^* a_*(\omega_X)[-q], \mathbb{C}(-x)[-q]) \simeq \tag{5.2.4}$$
$$\simeq \mathrm{Hom}_{\mathbf{D}(\widehat{A})}((-1_A)^* a_*(\omega_X), \mathbb{C}(-x)). \tag{5.2.5}$$

Where, as usual $\mathbb{C}(-x)$ stands for the skyscraper sheaf at the point $-x$, (5.2.4) is a consequence of Remark 5.2.2(b) and some usual computations with integral transform (combine Proposition 1.1.3 with the fact that $\mathbf{R}\widehat{S}(\mathscr{O}_{\widehat{A}} = \mathbb{C}(0)[-q])$ .

From the equality

$$\mathrm{Hom}_{\mathbf{D}(\widehat{A})}((-1_A)^* a_*(\omega_X), \mathbb{C}(-x)) = 0$$

we deduce $x$ is not in the support of $a_*\omega_X$ and hence it is not contained in $a(X)$. It follows that $a(X) \subseteq T$ contradicting the fact that it generates $A$. □

### 5.2.2 A Birationality Criterion

In this paragraph we present a birationality criterion, similar to the one exposed in [32, Theorem 3.1], which we will apply in order to prove Theorem 5.2.2.

**Proposition 5.2.4.** *Let $g : X \to Y$ be a generically finite morphism of smooth varieties of maximal Albanese dimension and let $b : Y \to A$ be a generically finite morphism into an abelian variety such that*

*(i)* $gv_b(g_*\omega_X) \geq 1$;

*(ii)* $\chi(X, \omega_X) = \chi(Y, \omega_Y)$.

*Then $g$ is birational.*

*Proof.* First of all observe that, since $g$ is generically finite, by Grauert–Riemenshneider vanishing we have $\chi(X, \omega_X) = \chi(Y, g_*\omega_X)$ and hence, using (ii),

$$\chi(Y, \omega_Y) = \chi(Y, g_*\omega_X). \tag{5.2.6}$$

Now there is a natural inclusion $\omega_Y \hookrightarrow g_*(\omega_X)$. In fact, we have a natural inclusion of differential 1-forms:

$$g^*\Omega^1_Y \hookrightarrow \Omega^1_X.$$

Taking determinants on both sides and observing that exterior powers commute with pullbacks, we get another inclusion

$$g^*\omega_Y \hookrightarrow \omega_X.$$

Pushing forward we obtain

$$g_*(g^*\omega_Y) \hookrightarrow g_*\omega_X;$$

we use projection formula in order to achieve

$$\omega_Y \otimes g_*\mathscr{O}_X \hookrightarrow g_*\omega_X.$$

Finally, by tensoring by $\omega_Y$ the natural inclusion $\mathscr{O}_Y \hookrightarrow g_*\mathscr{O}_X$ we get the seeked inclusion.

Therefore we can construct a short exact sequence

$$(5.2.7) \qquad 0 \to \omega_Y \longrightarrow g_*\omega_X \longrightarrow \mathscr{Q} \to 0.$$

By (5.2.6) and the additivity of the Euler characteristic we have

$$(5.2.8) \qquad \chi(Y, \mathscr{Q}) = 0.$$

Now we claim

$$(5.2.9) \qquad gv_b(\mathscr{Q}) \geq 1.$$

Before proceeding further let us see how (5.2.9) implies the statement. By the W.I.T. criterion (Theorem 1.2.4) we have that $\mathbf{R}\Phi_{\mathscr{P}_b}(\mathbf{R}\Delta_Y(\mathscr{Q}))$ is a sheaf $\mathbf{R}b_*\widehat{\mathbf{R}\Delta_Y}(\mathscr{Q})$ of generic rank 0 in degree $\dim Y$. But Theorem 1.2.9 tells us that $\mathbf{R}b_*\widehat{\mathbf{R}\Delta_Y}(\mathscr{Q})$ is torsion free. Therefore we necessarily have $\mathbf{R}b_*\widehat{\mathbf{R}\Delta_Y}(\mathscr{Q}) = 0$ which implies

$$0 = \mathbf{R}\Phi_{\mathscr{P}_b}(\mathbf{R}\Delta_Y(\mathscr{Q})).$$

By applying to (5.2.7) the functor $\mathbf{R}\Phi_{\mathscr{P}_b} \circ \mathbf{R}\Delta_Y$ we get an isomorphism of objects in $\mathbf{D}(\widehat{A})$:

$$\mathbf{R}\Phi_{\mathscr{P}_b}(\mathscr{O}_Y) \simeq \mathbf{R}\Phi_{\mathscr{P}_b}(\mathbf{R}\Delta_Y(g_*\omega_X)).$$

Now we use Grauert–Riemenshneider vanishing theorem together with Grothendieck–Verdier duality to obtain:

$$\mathbf{R}\Phi_{\mathscr{P}_b}\mathbf{R}\Delta(g_*\omega_X) \simeq \mathbf{R}\Phi_{\mathscr{P}_b}\mathbf{R}\Delta\mathbf{R}g_*\omega_X \simeq \mathbf{R}\Phi_{\mathscr{P}_b}\mathbf{R}g_*\mathbf{R}\Delta\omega_X \simeq \mathbf{R}\Phi_{\mathscr{P}_{b\circ g}}(\mathscr{O}_X).$$

In particular we get the following isomorphisms

(5.2.10)

$$\mathbf{R}S(\mathbf{R}b_*\mathscr{O}_Y) \simeq \mathbf{R}\Phi_{\mathscr{P}_b}(\mathscr{O}_Y) \simeq \mathbf{R}\Phi_{\mathscr{P}_b}(\mathbf{R}\Delta_Y(g_*\omega_X)) \simeq \mathbf{R}\Phi_{\mathscr{P}_{b\circ g}}(\mathscr{O}_X) \simeq \mathbf{R}S(\mathbf{R}(b\circ g)_*\mathscr{O}_X).$$

It follows from Mukai's inversion theorem and (5.2.10) above that $\mathbf{R}b_*\mathscr{O}_Y$ is isomorphic to $\mathbf{R}(b\circ g)_*\mathscr{O}_X$ as objects in $\mathbf{D}(A)$. In particular we get that $b_*\mathscr{O}_Y \simeq (b\circ g)_*\mathscr{O}_X$ as sheaves on $A$. Therefore we have

$$\deg b = \operatorname{rank} b_*\mathscr{O}_Y = \operatorname{rank}(b\circ g)_*\mathscr{O}_X = \deg b\cdot\deg g.$$

Necessarily the generic rank of $g$ must be 1, and the statement is proved.

In order to conlude the proof we just have to verify (5.2.9). To this aim, consider again the short exact sequence (5.2.7) and twist it by $b^*\alpha$ with $[\alpha]\in\widehat{A}$. Taking cohomology we get

$$\cdots\to H^i(Y, g_*\omega_X\otimes b^*\alpha)\to H^i(Y,\mathscr{Q}\otimes b^*\alpha)\to H^{i+1}(Y,\omega_Y\otimes b^*\alpha)\to\cdots$$

Therefore there is an inclusion

$$V_b^i(\mathscr{Q})\subseteq V_b^i(\omega_X)\cup V_b^{i+1}(\omega_Y).$$

Now by Theorem 1.2.3 $\operatorname{codim}V_b^{i+1}(\omega_Y)\geq i+2$. In addition, by condition (i) in the statement we have that also $\operatorname{codim}V_b^i(g_*\omega_X)\geq i+1$. It follows at once that $\operatorname{codim}V_b^i(\mathscr{Q})\geq i+1$ too and (5.2.9) holds.

□

Now we are ready to challenge Theorem 5.2.2. In fact consider $(X,a)$ a pair satisfying Hypothesis 5.2.1. Denote by $Y'$ the image of $a$ and take $b: Y\longrightarrow Y'$ to be a desingularization of its. Then there exists a variety $X'$ and a generically finite morphism $g: X'\longrightarrow Y$ such that $b\circ g$ is birationally equivalent to $a: X\longrightarrow Y'$. Since the issues we are addressing are birational in nature, we can and will assume, without loss of generality, that $X=X'$ and we have the following commutative diagram:

(5.2.11)

$$\begin{array}{ccccc} X & \xrightarrow{\quad a\quad} & Y & \hookrightarrow & A\\ & {\scriptstyle g}\searrow & \nearrow{\scriptstyle b} & & \\ & & Y' & & \end{array}$$

Given the birationality of $b$, to ensure the birationality of $a$ it is enough to prove that $g$ is birational. Let us show that the hypothesis of Proposition 5.2.4 are satisfied.

First of all, observe that, as in the proof of the birationality criterion above, we have

$$\mathbf{R}\Phi_{\mathscr{P}_b}\mathbf{R}\Delta(g_*\omega_X)\simeq\mathbf{R}\Phi_{\mathscr{P}_b}\mathbf{R}\Delta\mathbf{R}g_*\omega_X\simeq\mathbf{R}\Phi_{\mathscr{P}_b}\mathbf{R}g_*\mathbf{R}\Delta\omega_X\simeq\mathbf{R}\Phi_{\mathscr{P}_a}\mathbf{R}\Delta\omega_X.$$

Since $(X, a)$ satisfies Hypothesis 5.2.1, by Remark 5.2.1, $gv_a(\omega_X) \geq 1$. In particular (cfr. Theorems 1.2.4 and 1.2.9), $\mathbf{R}\Phi_{\mathscr{P}_a}\mathbf{R}\Delta\omega_X$ is a torsion free sheaf in degree $\dim X = \dim Y$. It follows that also $\mathbf{R}\Phi_{\mathscr{P}_b}\mathbf{R}\Delta(g_*\omega_X)$ is a torsion free sheaf in degree $\dim Y$, hence, applying again Theorems 1.2.4 and 1.2.9, we get that $gv_b(g_*\omega_X) \geq 1$.

Thus, in order to conclude, we just need to verify that $\chi(\omega_Y) = \chi(\omega_X) = 1$. Observe that, thanks to the natural inclusion $j : \omega_Y \to g_*(\omega_X)$, the generic vanishing properties of $\omega_Y$ and $g_*\omega_X$, and Grauert–Riemenshneider vanishing, if we take $[\alpha] \in \widehat{A}$ general we get

$$\chi(Y, \omega_Y) = h^0(Y, \omega_Y \otimes b^*\alpha) \leq h^0(Y, g_*\omega_X \otimes b^*\alpha) = \chi(Y, g_*\omega_X) = \chi(\omega_X) = 1.$$

Let us suppose that $\chi(Y, \omega_Y) = 0$, by Proposition 1.2.10 we would have that $gv_b(\omega_Y) = 0$. But then $b(Y) = Y' = a(Y)$ would be fibered in tori contradicting Hypothesis 5.2.1. □

## 5.3 A Characterization of Theta Divisors

In this section we further refine a cohomological characterization of $\Theta$ divisors due to Pareschi who in [59] proved the following

**Theorem 5.3.1.** *Let $X$ be a smooth projective variety of maximal Albanese dimension with $\chi(\omega_X) = 1$ and $gv(\omega_X) = 1$. Suppose furthermore that for every $0 < i < \dim X$ the codimension of $V^i(\omega_X) > i + 1$. Then $X$ is birational to a theta-divisor.*

We will prove the follwing

**Theorem 5.3.2.** *Let $(X, a)$ be a pair satisfying Hypothesis 5.2.1. Suppose furthermore that $[\mathscr{O}_A] \in \widehat{A}$ is an isolated point of $V^i_a(\omega_X)$ for every $i > 0$ and that $a(X)$ is normal. Then the abelian variety $A$ is principally polarized and $X$ is birational to a $\Theta$-divisor in $A$.*

We will begin by proving the preliminary result.

**Proposition 5.3.3.** *Let $(X, a)$ be as in the Theorem above. Then the dimension of $A$ is $n + 1$.*

*Proof.* The assumptions on $X$, namely that both the Euler characteristic of $\omega_X$ and its generic vanishing index with respect to the map $a$ are equal to 1, ensure that $\dim A > n$, $\mathbf{R}a_*\mathscr{O}_X$ is a W.I.T. object with index $q = \dim A$ and that its Fourier-Mukai transform, $\widehat{\mathbf{R}a_*\mathscr{O}_X}$ is a torsion free sheaf of generic rank 1 (crf. Proposition 1.2.6 and Theorem 1.2.9), i.e. it is an ideal sheaf $\mathscr{I}_Z$ twisted by some line bundle $L$. Moreover, thank to the fact that any codimension $i + 1$ component of $V^i_a(X, \omega_X)$ does not pass through $0_{\widehat{A}}$, we can write $Z = S \cup W$ with $S$ a scheme supported on $0_A$ and $W$ a scheme not containing $0_A$. Then we can look at the short exact sequence

$$0 \to \mathscr{I}_Z \otimes L \longrightarrow L \longrightarrow L_{|Z} \to 0$$

Taking Ext-sheaves and observing the vanishing, we have that for every $i \geq 1$

$$\mathscr{E}xt^i(\mathscr{I}_Z \otimes L, \mathscr{O}_{\widehat{A}}) \simeq \mathscr{E}xt^{i+1}(L_{|Z}, \mathscr{O}_{\widehat{A}}) \tag{5.3.1}$$

Now any component of the support of $\mathscr{E}xt^j(\mathscr{I}_Z \otimes L, \mathscr{O}_{\widehat{A}})$ is a component of the support of $\mathbf{R}^j S(\mathbf{R}a_*\omega_X)$ that is empty for $j > n$. On the other side $S$ is a component of codimension $q$ of $Z$, hence it is a component of $\mathscr{E}xt^q(\mathscr{O}_Z \otimes L, \mathscr{O}_{\widehat{A}})$. Hence $\mathscr{E}xt^q(L_{|Z}, \mathscr{O}_{\widehat{A}})$ is non-zero. Equation (5.3.1) yields that $\mathscr{E}xt^{q-1}(\mathscr{I}_Z \otimes L, \mathscr{O}_{\widehat{A}}) \neq 0$, and therefore $q - 1 \leq n$. Thus we get that $q = n + 1$. Furthermore, since by [2, Proposition 6.1]

$$k(\widehat{0}) \simeq R^n\Phi_{\mathscr{P}_a}(\omega_X) \simeq \mathscr{E}xt^n(\mathscr{I}_Z \otimes L, \mathscr{O}_{\widehat{A}}) \simeq \mathscr{E}xt^q(L_{|Z}, \mathscr{O}_{\widehat{A}})$$

we also get that $S$ is the reduced subscheme supported in 0. $\square$

We will also need some Lemmas.

**Lemma 5.3.4.** *Let $f : Y \to Z$ is a surjective morphism between varieties with $Y$ normal and $Z$ smooth, then the general fiber of $F$ is normal.*

*Proof.* A similar argument could be found in [70]. Since normality is a local property we can think $Y = \mathrm{Spec}(B)$, $Z = \mathrm{Spec}(A)$ with, $A$ and $B$ $\mathbb{C}$-algebras with $B$ a normal over $\mathbb{C}$, and $f$ induced by a morphism $\varphi : A \to B$. Since $f$ is surjective, then any prime ideal of $A$ is contracted. In particular (0) is contracted and $\varphi$ is injective. In this setting the general fiber is just $[\varphi(A\backslash\{0\})^{-1}] \cdot B$. Since normality is stable by localization, $[\varphi(A\backslash\{0\})^{-1}] \cdot B$ is normal. We get the statement. $\square$

**Lemma 5.3.5.** *Let $A$ be an abelian variety of dimension $q \geq 2$ and take $D$ an ample divisor in $A$. If $[\alpha] \in \mathrm{Pic}^0(A)$ is such that $\alpha_{|D} \simeq \mathscr{O}_D$ then $\alpha \simeq \mathscr{O}_A$.*

*Proof.* Consider the standard short exact sequence

$$0 \to \mathscr{O}_A(-D) \longrightarrow \mathscr{O}_A \longrightarrow \mathscr{O}_D \to 0$$

and tensor it by a topologically trivial line bundle $\alpha$ as in the statement. Suppose that $\alpha$ is not trivial. Taking cohomology we get

$$0 \to H^0(D, \mathscr{O}_D) \to H^1(A, \mathscr{O}_A(-D) \otimes \alpha) \to \cdots$$

Since $q \geq 2$, then $H^1(A, \mathscr{O}_A(-D) \otimes \alpha) = 0$ and $H^0(D, \mathscr{O}_D) = 0$, yielding a contradiction. Necessarily we must have $\alpha \simeq \mathscr{O}_A$. $\square$

As a corollary we get

**Corollary 5.3.6.** *In the notation of the above Lemma, suppose furthermore that $D$ is normal and take $\mu : Y \to D$ a desingularisation. If $[\alpha] \in \widehat{A}$ is such that $\mu^*\alpha$ is trivial then $\alpha \simeq \mathscr{O}_A$.*

*Proof.* By the normality of $D$ combined with projection formulas we have

$$\mathscr{O}_D \simeq \mu_* \mathscr{O}_Y \simeq \mu_*(\mu^*\alpha) \simeq \alpha \otimes \mathscr{O}_D.$$

Lemma 5.3.5 implies then that $\alpha$ is trivial. □

Now we are ready to challenge the Theorem.

*Proof of Theorem 5.3.2.* Let $a$ be the Albanese map of $X$. We can suppose that there exists $0 < k < \dim X$ and $W$ a component of $V_a^k(\omega_X)$ of codimension $k+1$. Otherwise we would be in the hypothesis of Theorem 5.3.1 and there would be nothing to prove. By the Subtorus Theorem (Theorem 3.1.9) $W = [\alpha] + \widehat{B}$ with $\widehat{B}$ a subabelian variety of $\widehat{A}$ and $[\alpha] \in \widehat{A}$ a torsion point. We will show that $[\alpha] \in \widehat{B}$, contraddicting the assumption that $[\mathscr{O}_A]$ is an isolated point of $V_a^k(\omega_X)$.

As usual, call $\pi : A \to B$ the map dual to the inclusion $\widehat{B} \hookrightarrow \widehat{A}$. Let $g := \pi \circ a$ and denote by $K$ the kernel of $\pi$. We are in the following situation:

$$\begin{array}{ccc} G & \xrightarrow{a_{|G}} & K \\ {\scriptstyle j}\downarrow & & \downarrow{\scriptstyle i} \\ X & \xrightarrow{a} & A \\ & \searrow^{g} & \downarrow{\scriptstyle \pi} \\ & & B \end{array}$$

Where $G$ is the general fiber of $g$ and $j$ is the natural inclusion.

By Lemma 5.1.4 combinded with Proposition 5.3.3, we have that

$$\dim g(X) = \dim X - k = \dim A - k - 1 = \dim B,$$

hence $g$ is surjective. In particular we have that also $\pi_{|a(X)}$ is surjective and $G$ is of dimension $k$, hence via the restirction of $a$ maps to a divisor $D$ (of a translate of) $K$. This divisor $D$ is a general fiber of a surjective map from a normal variety into a smooth one and therefore is normal by Lemma 5.3.4. We recall now that $a$ has generic degree one (Theorem 5.2.2), therefore we can think of $a_{|F} : F \to K$ as a desingularisation of $D$ and

$$a_{|F}^*(\alpha_{|K}) \simeq (a^*\alpha)_{|F} \in \mathrm{Pic}^0(F).$$

Now observe that, by Theorem 5.1.2 with $g = a = \mathrm{id}_B$ the sheaves $R^j f_*(\omega_X \otimes a^*\alpha)$ satisfy GV. It follows that for $[\beta] \in \widehat{B}$ general

$$h^0(B, R^k f_*(\omega_X \otimes a^*\alpha) \otimes \beta) \simeq h^k(X, \omega_X \otimes a^*\alpha \otimes f^*\beta) \neq 0$$

where the last inequality depends on the fact that $[\alpha] + \widehat{B} \subseteq V^k(\omega_X)$. We can conclude that $R^k f_*(\omega_X \otimes a^*\alpha) \neq 0$. Since by Theorem 5.1.1 it is torsion free, for the general $b \in B$ we have

$$0 \neq R^k f_*(\omega_X \otimes a^*\alpha) \otimes \mathbb{C}(b) \simeq H^k(F, \omega_F \otimes (a^*\alpha)_{|F}).$$

We deduce that $(a^*\alpha)_{|F} \simeq a^*_{|F}(\alpha_{|K})$ is trivial and Corollary 5.3.6 immediatiely implies that $\alpha_{|K} \simeq \mathscr{O}_K$ dualizing

$$0 \to K \longrightarrow A \longrightarrow B \to 0$$

we get

$$0 \to \widehat{B} \longrightarrow \widehat{A} \longrightarrow \widehat{K} \to 0$$

where the rightmost arrow is the restriction to $K$. Hence $[\alpha] \in \mathrm{Ker}((-)_{|K} : \widehat{A} \to \widehat{K}) = \widehat{B}$ and the statement is proved. □

# BIBLIOGRAPHY

[1] M Aprodu and G Farkas. Green's conjecture for the general cover. To appear in "Vector bundles and compact moduli-Athens, Georgia 2010", Contemporary Mathematics.

[2] M. A. Barja, M. Lahoz, J. C. Naranjo, and G. Pareschi. On the bicanonical map of irregular varieties. Preprint arXiv:0907.4363.

[3] A. Beauville. L'inegalité $p_g \geq 2q - 4$ pour les surfaces de type général. *Bull. Soc. Math. France*, 110:343–346, 1982.

[4] A. Beauville. The Coble hypersurface. *C. R. Acad. Sci. Paris*, 337(I):189–194, 2003.

[5] C. Birkenhake and H. Lange. *Complex abelian varieties*. Springer-Verlag, 1992.

[6] E. Bombieri. Canonical models for surfaces of general type. *Inst. Hautes Études Sci. Publ. Math.*, 42:127–219, 1973.

[7] J. A. Chen, O. Debarre, and Z. Jiang. Varieties with vanishing holomorphic Euler characteristic. Preprint arXiv:1105.3418, 2011.

[8] J. A. Chen and C. D. Hacon. Characterization of abelian varieties. *Invent. Math.*, 143:435–447, 2001.

[9] J. A. Chen and C. D. Hacon. Linear series of irregular varieties. Proceedings of the symposium on Algebraic Geometry in East Asia. World Scientific., 2002.

[10] J.A. Chen and C. D. Hacon. Pluricanonical maps of varieties of maximal Albanese dimension. *Mathematische Annalen*, 320(2):367–380, 2001.

[11] J.A. Chen and C. D. Hacon. On the Irregularity of the Iitaka fibration. *Comm. in Algebra*, 32(1):203–215, 2004.

[12] M Cornalba, X. Gomez-Mont, and A. Verjosvsky, editors. *Proceedings Of The College On Riemann Surfaces, International Centre For Theoretical Physics.* World Scientific, 1989.

[13] L. Di Biagio. Pluricanonical systems for 3-folds, 4-folds and $n$-folds of general type. Ph-D Thesis. Defended on July 2009 at Università di Roma "La Sapienza".

[14] L. Di Biagio. Pluricanonical systems for 3-folds and 4-folds of general type. Preprint arXiv:1001.3340v2.To appear in Mathematical Proceedings of the Cambridge Philosophical Society.

[15] L. Ein and R. Lazarsfeld. Syzigies and Koszul cohomology of smooth projective variety of any dimension. *Invent. math.*, 111:51–67, 1993.

[16] L. Ein and R. Lazarsfeld. Singularities of theta divisors and the birational geometry of irregular varieties. *Journal of the American Mathematical Society*, 10(1):243–258, 1997.

[17] L. Ein and R. Lazarsfeld. Asymptotic syzygies of algebraic varieties. March 2011.

[18] D. Eisenbud. *Commutative Algebra with a View Toward Algebraic Geometry.* Graduate Texts in Mathematics. Springer-Verlag, 2004.

[19] D. Eisenbud. *The Geometry of Syzygies: a Second Course in Commutative Algebra and Algebraic Geometry.* Birkhäuser, 2005.

[20] E. G. Evans and P. Griffith. The syzygy problem. *Ann. Math.*, 114:323–333, 1981.

[21] G Farkas. Syzygies of curves and the effective cone of $M_g$. *Duke Mathematical Journal*, 135:53–98, 2006.

[22] G Farkas. Aspects of the birational geometry of $M_g$. *Surveys in Differential Geometry*, 14:57–111, 2010.

[23] F. J. Gallego and B. P. Purnapranjna. Syzygies of projective surfaces: an overview. *J. Ramanujan Math. Soc.*, 14:65–93, 1999.

[24] M. Green. Koszul cohomology and the geometry of projective varieties I. *J. Diff. Geometry*, 19:125–171, 1984.

[25] M. Green. Koszul cohomology and the geometry of projective varieties II. *J. Diff. Geometry*, 20:279–289, 1984.

[26] M. Green. Koszul cohomology and geometry. In Cornalba et al. [12], page 666.

[27] M. Green and R. Lazarsfeld. Some results on the syzygies of finite set and algebraic curves.

[28] M. Green and R. Lazarsfeld. On the projective normality of complete linear series on an algebraic curve. *Invent. math.*, 83:73–90, 1986.

[29] M. Green and R. Lazarsfeld. Higher obsturciotns to deforming cohomology groups of line bundles. *J. Amer. Math. Soc*, 1(4):87–103, 1991.

[30] C. D. Hacon. Fourier-Mukai transform, generic vanishing theorems and polarizations in abelian varieties. *Mathematische Zeischrift*, 235:717–726, 2000. Also preprint arXiv: 9902078v1.

[31] C. D. Hacon and J. McKernan. Boundness of pluricanonical maps of varieties of general type. *Invent. Math.*, 166(1):1–25, 2006.

[32] C. D. Hacon and R. Pardini. On the birational geometry of varieties of maximal Albanese dimension. *J. Reine Angew. Math*, 546, 2002.

[33] C. D. Hacon and R. Pardini. Surfaces with $p_g = q = 3$. *Transactions of the A. M. S.*, 354:2631–2638, 2002.

[34] C. D. Hacon and R. Pardini. Birational charaterization of product of curves of genus 2. *Mathematical Research Letters*, 12:129–140, 2005.

[35] S. Iitaka. *Algebraic Geometry: An Introduction to Birational Geometry of Algebraic Varieties*. Springer, New York, 1981.

[36] Z. Jiang. On varieties of maximal Albanese dimension. September 2009. arXiv: 0909.4817.

[37] Z. Jiang. An effective version of a theorem of Kawamata on the Albanese map. *Communication in Contemporary Mathematics*, 13(3):509–532, 2011. arXiv: 0909.3973.

[38] Z. Jiang. Private e-mail communication. Unpublished, 2011.

[39] Z. Jiang, M. Lahoz, and S. Tirabassi. On the Iitaka fibration of irregular varieties. Preprint arXiv:1111.6279. Submitted, 2011.

[40] Y. Kawamata. Characterization of Abelian Varieties. *Comp. Math*, 43:253–276, 1981.

[41] G. Kempf. Linear system on abelian varieties. *Am. J. Math.*, 111:65–94, 1989.

[42] G. Kempf. Projective coordinate ring of abelian varieties. *Algebraic analysis, geometry and number theory*, pages 225–236, 1989. MR **98m**:14047.

[43] G. Kempf. *Complex Abelian Varieties and Theta Functions.* Springer-Verlag, 1990.

[44] G. Kempf. Equations of Kummer varieties. *Am. J. Math.*, 114 (n. 2):229–232, 1992.

[45] A. Khaled. Equations des variétés de Kummer. *Math. Ann.*, 295:685–701, 1993.

[46] A. Khaled. Projective normality and equations of Kummer varieties. *J. reine angew. Math.*, 465:197–217, 1995.

[47] S. Koizumi. Theta relations and projective normality of abelian varieties. *Am. J. Math*, 98:865–889, 1976.

[48] M. Lahoz. Generic vanishing index and the birationality of the bicanonical map of irregular varieties. Preprint arXiv:1007.3026, 2010.

[49] R. Lazarsfeld. A sampling of vector bundles techniques in the study of linear series. In Cornalba et al. [12], pages 500–559.

[50] R. Lazarsfeld. *Positivity in algebraic geometry I & II*, volume 4. Springer-Verlag, 2004.

[51] R. Lazarsfeld, G. Pareschi, and M. Popa. Local positivity, multiplier ideals, and syzygies of abelian varieties. Preprint arXiv:1003.4470, 2010.

[52] R. Lazarsfeld and M. Popa. BGG correspondence of cohomology of compact Kähler manifolds, and numerical invariants. Preprint arXiv:0907.0651v1, 2009.

[53] S. Mukai. Duality between $D(X)$ and $D(\widehat{X})$ with its application to picard sheaves. *Nagoya Math. J.*, 81:153–175, 1981.

[54] D. Mumford. On the equations defining abelian varieties. *Invent. Mat*, 1:287–354, 1966.

[55] D. Mumford. *Abelian varieties.* Oxford Uni. Press., 1974.

[56] A. Ohbuchi. A note on the normal generation of ample line bundles on abelian varieties. *Proceedings of the Japan Academy, Series A, Mathematical Sciences*, 64(4):119–120, 1988.

[57] G. Ottaviani and R. Paoletti. Syzygies of the Veronese embedding. *Compositio Mathematica*, 124:31–37, 2001.

[58] G. Pareschi. Syzigies of abelian varieties. *J. Amer. Math. Soc.*, 13(3):651–664, 2000.

[59] G. Pareschi. Basic results on irregular varieties via Fourirer-Mukai methods. Preprint, 2010.

[60] G. Pareschi and M. Popa. Regularity on abelian varieties III: relationship with Generic Vanishing and applications. Preprint ArXiv:0802.1021, to appear in the Proceedings of the Clay Mathematics Institute **14** (2011).

[61] G. Pareschi and M. Popa. Regularity on abelian varieties I. *J. Amer. Math. Soc.*, 16:285–302, 2003.

[62] G. Pareschi and M. Popa. Regularity on abelian varieties II: basic results on linear series and defining equations. *J. Algebraic Geometry*, 13:167–193, 2004.

[63] G. Pareschi and M. Popa. Castelnuovo theory and the geometric Shottky problem. *J. Reine Angew. Math.*, 615:25–44, 2008.

[64] G. Pareschi and M. Popa. Strong Generic Vanishing and Higer-dimensional Castelnuovo–de Franchis inequality. *Duke J. Math.*, 150(2):269–285, 2009.

[65] G. Pareschi and M. Popa. GV-sheaves, Fourier-Mukai transform and Generic Vanishing. *Amer. J. Math*, 133(1):235–271, 2011. Also preprint arXiv 0608127v4. MR:2752940.

[66] G. P. Pirola. Surfaces with $p_g = q = 3$. *Manuscripta Mathematica*, 108(2):163–170, 2002.

[67] D. Ploog. Equivariant autoequivalences for finite group actions. *Advances In Mathematics*, 216:62–74, 2007.

[68] R. Sasaki. Bounds on the degree of equations defining Kummer varieties. *J. Math. Soc. Japan*, 34 (n. 2):223–239, 1981.

[69] C. Simpson. Subspaces of moduli spaces of rank one local systems. *Ann. Sci.Éc. Norm. Sup*, 26:361–401, 1993.

[70] K. Swhede. Answer to a question. MathOverflow.

[71] S. Takayama. Pluricanonical systems on algebraic varieties of general type. *Invent. Math.*, 165(3):551–587, 2006.

[72] S. Tirabassi. On the tetracanonical map of varieties of general type and maximal Albanese dimension. *Collectanea Mathematica*, 2011. Also preprint availabe at arXiv:1103.5236.

[73] C Voisin. Green's generic syzygy conjecture for curves of even genus lying on a K3 surface . *J. European Math. Society*, 4:363–404, 2002.

[74] C Voisin. Green's canonical syzygy Conjecture for generic curves of odd genus. *Compositio Mathematica*, 141:1163–1190, 2005.

# LIST OF SYMBOLS

| | |
|---|---|
| $\mathrm{Alb}(X)$ | Albanese variety of $X$, page viii |
| $\mathrm{alb}_X$ | Albanese morphism of $X$, page viii |
| $q(X)$ | Irregularity of $X$, page x |
| $\chi(X)$, $\chi(\omega_X)$ | Euler characteristic of $X$, page ix |
| $\mathbf{Coh}(X)$ | Abelian category of coherent sheaves on $X$, page xi |
| $\mathbf{D}(X)$ | Bounded derived category of coherent sheaves on $X$, page xi |
| $H^i(X, \mathscr{F})$, $H^i(\mathscr{F})$ | i-th cohomology group of $X$, page xi |
| $H^i(X, \mathscr{F})$, $h^i(\mathscr{F})$ | Dimension of $H^i(X, \mathscr{F})$ over the base field, page xi |
| $k(x)$ | Skyscraper sheaf at $x$, page xi |
| $\mathrm{Bs}(V)$ | Base locus of $V \subseteq H^0(X, \mathscr{F})$, page xi |
| $\widehat{A}$ | $\mathrm{Pic}^0(A)$, dual abelian variety to $A$, page 1 |
| $\mathscr{P}_A$, $\mathscr{P}$ | Poincaré line bundle on $A \times \widehat{A}$, page 1 |
| $\mathbf{R}S_A$ | Classical Fourier-Mukai functor in the abelian variety $A$., page 1 |
| $-1_A$ | Multiplication by -1 in $A$, page 1 |
| $[\cdot]$ | Shift functor in the derived category, page 1 |
| $\mathbf{R}S$, $\mathbf{R}\widehat{S}$ | Classical Fourier Mukai functor, page 2 |
| $R^iS(\mathscr{F})$,$R^i\widehat{S}(\mathscr{F})$ | Cohomology groups of the complex $\mathbf{R}S(\mathscr{F})$, $\mathbf{R}\widehat{S}(\mathscr{F})$, page 2 |
| $\widehat{\mathscr{F}}$ | Fourier-Mukai transform of $\mathscr{F}$, page 2 |
| W.I.T.$(j)$ | Objects satisfying W.I.T. with index $j$, page 2 |
| I.T.$(j)$ | Sheaves satisfying I.T. with index $j$., page 2 |
| $t_p$ | Translation by $p$, page 2 |
| $P_p$ | $\mathscr{P}_{\{p\}\times\widehat{A}}$, page 2 |
| $\mathbf{R}\Phi_{\mathscr{E}}$ | Integral transform with kernel $\mathscr{E}$, page 3 |
| $\mathscr{P}_a$ | $(a\times\mathrm{id})^*\mathscr{P}$ , page 3 |
| $\mathbf{R}\Phi_{\mathscr{P}_a}$ | $\mathbf{R}S_A\circ\mathbf{R}a_*$, page 3 |
| $\mathbf{R}\Delta_Z$ | Dualizing functor of $Z$, page 3 |
| $[\alpha]$ | Point in $\widehat{A}$ associated to the line bundle $\alpha$, page 2 |
| $V^i_a(X,\mathscr{F})$, $V^i_a(\mathscr{F})$ | $i$-th cohomological support locus, page 4 |
| $\mathrm{gv}_a(\mathscr{F})$ | Generic vanishing index of $\mathscr{F}$, page 4 |
| $\mathrm{gv}_a(X)$ | Generic vanishing index of $X$, page 4 |
| $GV_k$ | $\mathscr{F}$ such that $\mathrm{gv}_a(\mathscr{F})\geq k$, page 4 |
| $\mathcal{K}_X$ | Kummer variety associated to $X$, page 9 |
| $R_{\mathscr{A}}$ | Sections ring associated to the line bundle $\mathscr{A}$, page 9 |
| $S_{\mathscr{A}}$ | Symmetric algebra of $H^0(Z\mathscr{A})$, page 9 |
| $N_p$ | Property of syzygies of projective varieties, page 10 |
| $N^r_0$, $N^r_p$ | Property of syzygies of projective varieties, page 10 |
| $\pi_X$ | Quotient map $X\to\mathcal{K}_X$, page 11 |
| $M_L$ | Kernel of the evaluation map of $L$, page 13 |
| $\mathcal{K}_X$ | Kummer variety associated to $X$, page 15 |
| $\pi_X : X\longrightarrow\mathcal{K}_X$ | Quotient map, page 15 |
| $M_W$ | Kernel of the evaluation map $W\otimes\mathscr{O}_Z\to L$, page 14 |
| $\psi_L$ | Mumford normalized isomorphism, page 16 |
| $-\mathbf{1}_X$ | Involution on $X$ given by $x \mapsto -x$, page 17 |
| $[-1]_{\mathscr{L}}$ | Involution induced in cohomology by $i_X$ and the linearization on $\mathscr{L}$, page 17 |
| $H^0(X,\, \mathscr{L})^{\pm}$ | the eigenspace associated to the eigenvalue $\pm 1$ of the involution $[-1]_{\mathscr{L}}$, page 17 |
| $m^+_\alpha$ | A multiplication map of sections, page 20 |
| $\xi$ | Isogeny $X \times X \to X \times X$ given by $(p_1 + p_2,\ p_1 - p_2)$, page 20 |
| $\Phi^\alpha_\beta$ | An isomorphism, page 20 |
| $\Phi_\beta$ | $\Phi^\alpha_\beta$ with $[\beta] \simeq [\alpha^\vee]$, page 21 |
| $W_n$ | $H^0(X,\, \mathscr{A}^{2n})^+$ , page 27 |
| $M_{W_n}$ | $\pi^*_X M_{A^{\otimes n}}$, page 27 |
| $\mathrm{exc}(\mu)$ | Exceptional locus of $\mu$, page 38 |
| $\mathscr{J}(c \cdot \lvert V \rvert)$ | Multiplier ideal of the linear series $V$, page 39 |
| $\mathscr{J}(c \lvert\lvert D \rvert\rvert)$ | Asymptotyc multiplier ideal associated to $D$, page 39 |
| $a_X$ | Albanese map of $X$, page 48 |
| $P_m(X)$ | $m$-th plurigenus of $X$, page 48 |
| $\mathcal{H}_X$ | A special line bundle, page 51 |

# INDEX